\documentclass[review,onefignum,onetabnum]{siamart250211}

\usepackage[backend=biber, style=numeric, maxnames=99]{biblatex}
\usepackage{comment}
\usepackage{amssymb}
\nolinenumbers

\usepackage{algorithm}
\usepackage{algorithmicx}
\usepackage{algpseudocode}

\usepackage{booktabs}

\usepackage{subcaption}

\usepackage{diagbox}

\usepackage{float}   
\usepackage{caption} 

\newtheorem{assumption}[theorem]{Assumption}

\newtheorem{remark}[theorem]{Remark}
\usepackage{xcolor}

\ifpdf
  \DeclareGraphicsExtensions{.eps,.pdf,.png,.jpg}
\else
  \DeclareGraphicsExtensions{.eps}
\fi

\headers{A fast particle-field method for 3D Chemotaxis systems}{Jingyuan Hu, Zhongjian Wang, Jack Xin and Zhiwen Zhang}

\title{A fast stochastic interacting particle-field method for 3D parabolic-parabolic Chemotaxis  systems: numerical algorithms and error analysis}

\author{Jingyuan Hu\thanks{Department of Mathematics, The University of Hong Kong, Pokfulam Road, Hong Kong SAR,
P.R.China. hujy@connect.hku.hk.}
\and Zhongjian Wang\thanks{Division of Mathematical Sciences, School of Physical and Mathematical Sciences, Nanyang Technological University, 21 Nanyang Link, Singapore 637371. zhongjian.wang@ntu.edu.sg}.
\and Jack Xin\thanks{Corresponding author. Department of Mathematics, University of California at Irvine, Irvine, CA 92697, USA.
jack.xin@uci.edu}.
\and Zhiwen Zhang\thanks{Corresponding author. Department of Mathematics, The University of Hong Kong, Pokfulam
Road, Hong Kong SAR, P.R.China. Materials Innovation Institute for Life Sciences and Energy
(MILES), HKU-SIRI, Shenzhen, P.R. China. zhangzw@hku.hk.}
}

\usepackage{amsopn}

\makeatletter
\newcommand*{\addFileDependency}[1]{
  \typeout{(#1)}
  \@addtofilelist{#1}
  \IfFileExists{#1}{}{\typeout{No file #1.}}
}
\makeatother

\date{}
\begin{document}
\maketitle

\begin{abstract}
In this paper, we develop a novel numerical framework, namely the stochastic interacting particle-field method with particle-in-cell acceleration (SIPF-PIC), for the efficient simulation of the three-dimensional (3D) parabolic-parabolic Keller-Segel (KS) systems. The SIPF-PIC method integrates Lagrangian particle dynamics with spectral field solvers by leveraging localized particle-grid interpolations and fast Fourier transform (FFT) techniques. For $P$ particles and $H$ Fourier modes per spatial dimension, the SIPF-PIC method achieves a computational complexity of $O(P + H^3 \log H)$ per time step, a significant improvement over the original SIPF method (proposed in \cite{SIPF1}), which has a computational complexity of $O(PH^3)$, while preserving numerical accuracy. Moreover, we carry out a rigorous error analysis for the proposed method and establish rigorous corresponding error estimates. Finally, we present numerical experiments to validate the convergence order and demonstrate the computational efficiency of the SIPF-PIC method. To the best of our knowledge, this work represents the first numerical approach that can resolve complex 3D blowup dynamics beyond single-point collapse, including the formation of ring-shaped singular structures.


\end{abstract}




\section{Introduction}
The Keller-Segel (KS) system \cite{KS} describes chemotaxis-\\driven aggregation processes, particularly modeling the dynamics of the slime mold amoeba Dictyostelium discoideum in response to the gradients of the chemoattractant. Here, we study the parabolic-parabolic KS system on a bounded spatial domain $\Omega \subset \mathbb{R}^d$:
\begin{equation}
\left\{
\begin{aligned}
\partial_t \rho &= \nabla \cdot \left( \mu \nabla\rho - \chi \rho \nabla c \right), \\
\epsilon \partial_t c &= \Delta c - k^2 c + \rho, \quad \mathbf{x} \in \Omega,\ t \in [0,T],
\end{aligned}
\right.
\label{eq:para_system}
\end{equation}
where $\rho(\mathbf{x},t)$ is the density of the particles, $c(\mathbf{x},t)$ is the concentration of the chemoattractant, and $\chi, \mu, \epsilon, k\,>0$ are physical constants. Formally taking $\epsilon = 0$ yields a parabolic-elliptic KS system that models instantaneous chemical equilibration.

The KS system demonstrates interdisciplinary utility, with applications spanning biological, ecological, and medical domains. In biology, these models describe cell aggregation and migration in bacterial colonies and cancer cells driven by chemotaxis \cite{Transport}. In ecology, they analyze population dynamics shaped by chemotactic strategies, such as algal bloom dynamics and resource-seeking behaviors \cite{crossref2}. In medical research, the KS system offers insights into both physiological mechanisms (e.g., tissue repair and immune response) and pathological processes (e.g., cancer metastasis) \cite{crossref3}.

Due to the nonlinear and potentially singular behavior of KS systems, especially in the presence of blow-up phenomena \cite{nagai1995blow, herrero1998self}, numerical methods have become essential for analyzing their solutions. Mesh-based numerical methods are among the most widely used approaches for solving the KS system. Chertock et al. \cite{pos_preserving} proposed a high-order hybrid finite volume and finite difference scheme with positivity-preserving properties for two-dimensional (2D) problems. Shen et al. \cite{shen2020unconditionally} proposed an energy dissipation and bound-preserving scheme that is not restricted to specific spatial discretizations. Chen et al. \cite{chen2022error} developed a fully-discrete finite element scheme for the 2D parabolic-elliptic case, and demonstrated that the scheme exhibits finite-time blow-up under assumptions analogous to those of the continuous case. Other works include a positivity-preserving and asymptotic-preserving semi-discrete scheme \cite{LiuJian-Guo2018Paap}, and a saturation concentration setting \cite{hillen2001global} which prevents blow-up and yields spiky solutions. Beyond these 2D cases, Gnanasekaran et al. \cite{3Dblowup} simulated a 3D blow-up solution in a cube with a finite difference scheme in a short time.
 
In addition to the Eulerian discretization methods mentioned above, there have been notable advancements in the Lagrangian framework for the KS system \eqref{eq:para_system} and related equations. The parabolic-elliptic case is comparatively simpler, characterized by a memoryless particle interaction system~\cite{velazquez2004point} with singularities arising when a pair of particles approach closely. Several particle methods seek to circumvent the singularity; Havskovéč and Ševčovič \cite{havskovec2009stochastic} employed a regularized interaction potential, and Liu et al. \cite{liu2017random} developed a random particle blob method with a mollified kernel. Convergence studies for these approaches have been conducted in \cite{havskovec2011convergence, mischler2013kac} and \cite{liu2019propagation}, with mean-field equations serving as a key tool \cite{mischler2013kac}. Other works include the computation of 2D advective parabolic-elliptic KS systems with passive flow \cite{khan2015global}, and the reduction of computational costs via the random batch method for particle interactions \cite{random_batch}.

Existing Lagrangian frameworks for parabolic-elliptic cases cannot be directly extended to fully parabolic systems, as the chemo-attractant concentration no longer achieves rapid local equilibration and the drift terms depend on historical particle positions \cite{chen2022mean, fournier2023particle}. Stevens \cite{stevens2000} developed an $N$-particle system for the fully parabolic case, representing the chemo-attractant concentration with particles subject to probabilistic creation and annihilation. However, this method also suffers from the accumulation of per-step complexity over time. To address this limitation, Wang et al. \cite{SIPF1} developed a stochastic interacting particle-field (SIPF) method that computes the chemo-attractant concentration through spectral decomposition and maintains only the current-step particle positions, which ensures that the per-step complexity remains constant over time.

However, we have to point out that in the SIPF method \cite{SIPF1}, each particle of $\rho$ interacts with every Fourier mode of $c$. This requires the computation of the trigonometric sum $\frac{1}{P}\sum_{p=1}^P e^{-\frac{2\pi i \mathbf{q}}{L}\cdot \widehat{\mathbf{X}}_p}$ for all wavenumbers $\mathbf{q}$, where $\widehat{\mathbf{X}}_p$ for $1\leq p \leq P$ are the positions of the particles (see Algorithm \ref{Alg:SIPF}). For 3D systems with $H$ Fourier modes per dimension, this results in an $O(PH^3)$ computational complexity per timestep. Such a computational bottleneck severely restricts the applicability of SIPF for large-scale 3D KS systems, particularly in scenarios requiring high-resolution and long-time simulations.

In this paper, to reduce the complexity, we propose the SIPF-PIC method (see Algorithm \ref{SIPF-PIC}) that integrates localized particle-grid operations with spectral field solvers accelerated by the fast Fourier transform (FFT). This hybrid method confines particle interactions to adjacent grid nodes while resolving global frequency-domain couplings via FFT, yielding a computational complexity of $O\left(P + H^3 \log H\right)$ per timestep, as presented in Theorem \ref{thm:overall_complexity}. The bidirectional interaction mechanism operates in each time step through two stages: particle-to-grid projection (see Algorithm \ref{Alg:P2G}) transfers particle-derived data to the computational grid, while grid-to-particle interpolation (see Algorithm \ref{Alg:G2P}) distributes field solutions from the grid back to individual particles; details are provided in Section \ref{sec:PIC}. The particle-to-grid projection in the SIPF-PIC method, which numerically realizes the trigonometric summation that is also known as the unequally spaced fast Fourier transform (USFFT) \cite{USFFT}, has been employed in Fourier-basis particle simulations \cite{USFFTS}. The grid-to-particle interpolation is employed in classical particle-in-cell methods \cite{alma99,TskhakayaD.2007TPM} and does not rely on spectral decompositions. \begin{color}{black}Some deterministic particle methods \cite{coppa96,cho18} also rely on particle-in-cell interactions to couple some computational particles with a physical grid. Although these methods provide properties of mass conservation, low-dissipation, and computationally efficient transport, repeated particle-grid transfers may trigger ringing instability, making grid-scale high-frequency numerical perturbations difficult to eliminate. In contrast, the SIPF-PIC method is based on a temporal discretization of stochastic particle dynamics, in which random diffusive displacements rapidly decorrelate high-frequency components and prevent their coherent accumulation.\end{color}

Our algorithm reduces the particle-grid interaction complexity from $O\left(P H^3\right)$ to $O\left(P + H^3 \log H\right)$, thus enabling the simulation of 3D KS systems with significantly increased particle numbers $(P>10^6)$ and finer grids $(H=256)$, resulting in improved numerical accuracy. \begin{color}{black}Furthermore, in Section~\ref{subsec:os}, we provide an operator-splitting interpretation of the algorithm in addition to the original basis-function perspective. This formulation shows that the method naturally preserves not only the positivity of the particle-based density $\hat{\rho}$ but also that of $\hat{c}$. It also clarifies that the essential role of FFT is to efficiently evaluate a single diffusion step for the periodic spatial domain and boundary conditions. Therefore, the proposed algorithm is applicable to domains of arbitrary geometry and to general boundary conditions, provided that an efficient one-step diffusion operator is available for the domain under consideration.\end{color}

We also conduct a systematic and in-depth study on the error analysis for the proposed SIPF-PIC method. The convergence of our method is theoretically guaranteed by the analysis in Section \ref{sec:Error}. Our main result, presented in Theorem \ref{thm:overall_error}, shows the convergence of the numerical solution \((\hat{\rho},\hat{c})\) obtained by the SIPF-PIC method to the reference solution \((\rho,c)\), under certain mild assumptions on $\rho$, $c$, and $\hat{c}$ (see Assumption \ref{as}). Specifically, the error in the particle distribution, measured by the Wasserstein distance \( \mathcal{W}_1 \), scales as \(O\left(H^{-\frac{16}{13}} + H^{\frac{4}{13}}P^{-\frac{1}{2}}\right)\) for a fixed computational time $T$. Furthermore, under specific parameter combinations, the assumption on \(\hat{c}\) can be relaxed (see Section~\ref{subsec:c_hat_relax}), which guarantees the unconditional stability of our algorithm. The proofs of these theorems rely on the Brownian coupling technique (see Definition \ref{def:br_cp}), which compares the particles to i.i.d. samples from the reference distribution with the same Brownian path.

To complement the theoretical analysis above, we also present numerical experiments to verify the convergence rate, computational complexity, and consistency with respect to $H$. We verify the critical mass \(8 \pi\) and the evolution of the second moment in a standard two-dimensional parabolic-elliptic system \eqref{eq:para_system}, where the dimension \(d=2\) and the parameters \((k, \epsilon, \mu, \chi) = (0, 0, 1, 1)\). For 3D parabolic-parabolic KS systems with \(k,\epsilon>0\), the critical mass lacks a precise characterization due to its dependence on the initial density distribution \(\rho(\mathbf{x}, 0)\). We thus validate the accuracy of the SIPF-PIC method under radially symmetric initial data, where the reference solution can be obtained by reducing the system to 1D and running a finite difference method on a very fine mesh.

Moreover, the SIPF-PIC method successfully reproduces complex blow-up patterns beyond the classical single-point radial type blow-up. To the best of our knowledge, this work presents the first numerical approach capable of resolving complex 3D blow-up dynamics beyond the single-point collapse regime, marking a key advance in the numerical simulation of such singular phenomena. We first investigate a non-radially symmetric configuration consisting of four uniform density balls located at the vertices of a tetrahedron, and then accurately resolve the ring-shaped blow-up pattern proposed in a recent paper \cite{ring_ref}, where mass concentrates on an evolving ring \(x_1^2 + x_2^2 = r(t)^2\) with a time-dependent radius \(r(t)\). In contrast, the original SIPF method, with parameter scalings ($P=10^4, H=24$), lacks sufficient resolution to capture the ring-shaped blow-up pattern. Specifically, the SIPF-PIC method successfully captures the line-density blow-up, which is characterized by mass concentration along evolving spatial structures rather than isolated points, and further captures its instability under symmetry-breaking perturbations \cite{collapsing}. These complex collapse mechanisms deviate fundamentally from classical delta-type blow-up dynamics, and the SIPF-PIC method, by leveraging large particle numbers and high-resolution grids, reliably captures their evolutionary processes in simulations. These results validate the performance of the SIPF-PIC method in addressing challenging 3D blow-up problems that were previously intractable by using existing methods.

\begin{color}{black}
The proposed SIPF-PIC method both inherits key ingredients from the previous SIPF method and introduces several substantial differences. As in the SIPF method, the particle positions are updated through a time discretization of the underlying stochastic differential equation (SDE), while a collection of Fourier basis functions, viewed as eigenfunctions of the diffusion operator under periodic boundary conditions, is used to evolve the concentration field. In the SIPF method, both the evaluation of the particle drift and the one-step deposition require direct interactions between all particles and all Fourier modes, resulting in a computational complexity of $O(PH^3)$. In contrast, the SIPF-PIC method employs local particle-in-cell interactions between each particle and its neighboring physical grid points, followed by FFT-based transformations between the physical and spectral grids, which reduces the overall complexity to $O(P+H^3\log H)$. In addition, we provide an operator-splitting interpretation formulated entirely in physical space. This viewpoint naturally reveals the positivity-preserving properties of both the density and the concentration and indicates how the method can be extended to non-periodic boundary conditions and geometrically complex domains whenever an efficient one-step diffusion operator is available. Another new contribution of the present work is a rigorous error analysis, which establishes the convergence of the SIPF-PIC method, and rules out the ringing instabilities or artificial blow-up generated solely by the numerical discretization.  
\end{color}

The rest of the paper is organized as follows. Section \ref{sec:num_method} presents theoretical finite-time blow-up results for the KS system \eqref{eq:para_system}, and introduces the SIPF method, a hybrid particle-field method derived through temporal discretization of the system. In Section \ref{sec:PIC}, we formalize the interpolation operators governing bidirectional particle-grid coupling, including particle-to-grid discretization and grid-to-particle interpolation, followed by the detailed formulation of the SIPF-PIC method. Section \ref{sec:Error} rigorously establishes a comprehensive error analysis framework that quantifies the approximation errors from various sources in the SIPF-PIC method. Section \ref{sec:Num} validates the method through numerical experiments, including convergence tests, performance benchmarks, and the visualization of 3D blow-up patterns. Finally, a summary of our findings and future research directions is given in Section \ref{sec:conclusion}.

\section{Preliminary results}
\label{sec:num_method}
In this section, we present preliminary analytical and numerical results in the study of parabolic-parabolic KS system.
\subsection{Analytical results of KS system}\label{subsec:review}
Recall the parabolic-parabolic KS system:
\begin{equation}
\left\{
\begin{aligned}
&\rho_t = \nabla \cdot \left( \mu \nabla\rho - \chi \rho \nabla c \right), \\
&\epsilon c_t = \Delta c - k^2 c + \rho, && \mathbf{x} \in \Omega,\ t \in [0,T],
\end{aligned}
\right.
\label{eq:system}
\end{equation}
where \( \rho(\mathbf{x},t) \) denotes the particle density, \( c(\mathbf{x},t) \) the chemoattractant concentration, with \( \chi, \mu > 0 \) (chemotactic sensitivity and diffusivity) and \( \epsilon, k > 0 \) (relaxation and degradation rates).
\paragraph{Finite-Time Blow-Up Behaviors}
For the parabolic-elliptic case \eqref{eq:system} with \(k = \epsilon = 0\) and \(\mu = 1\), the spatial dimension \(d\) governs the scaling properties of blow-up dynamics. In two dimensions (\(d = 2\)), the system is \(L^1\)-critical: the scale transformation \(
\rho_\lambda(\mathbf{x},t) = \frac{1}{\lambda^2} \rho\left(\frac{\mathbf{x}}{\lambda}, \frac{t}{\lambda^2}\right)  
\) maintains solution to system \eqref{eq:system}, and preserves the total mass \(\int_{\mathbb{R}^2} \rho(\mathbf{x},t) \mathrm{d}\mathbf{x}\). There is a critical mass threshold \( M_{\text{crit}} = \frac{8\pi}{\chi} \), independent of the initial density configuration. Subcritical initial distributions satisfying \( \int_{\mathbb{R}^2} \rho(\mathbf{x},0)\,\mathrm{d}\mathbf{x} < M_{\text{crit}} \) admit global-in-time solutions \cite{dolbeault2004optimal}, whereas supercritical initial masses with \( \int_{\mathbb{R}^2} \rho(\mathbf{x},0)\,\mathrm{d}\mathbf{x} > M_{\text{crit}} \) exhibit finite-time blow-up. In dimensions \(d \geq 3\), the system becomes mass-supercritical, since the scaling \(\rho_\lambda(\mathbf{x},t) = \frac{1}{\lambda^d} \rho\left(\frac{\mathbf{x}}{\lambda}, \frac{t}{\lambda^2}\right)\) preserves the \(L^{d/2}\) norm of \(\rho\) instead of \(L^1\). If the \(L^{d/2}\) norm of the initial density configuration is small enough, i.e., \(\|\rho(\mathbf{x},0)\|_{L^{d/2}} \leq N(d,\chi)\) for some dimension-dependent \(N(d,\chi)\), the system admits a global-in-time solution with a decreasing \(L^{d/2}\) norm \cite{corrias2004global}, that is, \(\|\rho(\mathbf{x},t)\|_{L^{d/2}} \leq \|\rho(\mathbf{x},0)\|_{L^{d/2}}\) for all \(t \geq 0\). Conversely, finite-time blow-up may occur \cite{collapsing} for an arbitrary initial mass \(M_0 = \|\rho(\mathbf{x},0)\|_{L^1}\).  

The blow-up is classified as Type I if \[\limsup_{t \to T^-} (T - t)\|\rho(t)\|_{L^\infty(\mathbb{R}^d)} < C\] for some constant \(C > 0\). Otherwise, the blow-up is classified as Type II. Type II blow-up requires \(d \geq 3\) and corresponds to singularities that do not collapse into a single-point delta measure, instead forming spatially distributed concentration patterns such as ring-shaped structures, as shown in \cite{ring_ref}.

\paragraph{Critical Mass}
For the 2D parabolic-elliptic case in \eqref{eq:system} with parameters \((k, \epsilon, \mu,\\ \chi) = (0, 0, 1, 1)\), the critical mass threshold \(M_{\text{crit}} = 8\pi\) governs the formation of finite-time singularities. Solutions with initial mass \(M_0  > 8\pi\) collapse into singular measures within finite time \(T^* < \infty\); while subcritical masses \(M_0 < 8\pi\) exhibit dispersion. The transition is characterized by the evolution of the second moment \cite{Transport}:  
\begin{equation}
\label{eq:second_moment}
\frac{\mathrm{d}}{\mathrm{d}t}m_2(t) = 4M_0\left(1 - \frac{M_0}{8\pi}\right), 
\end{equation}
where
\begin{equation}
\quad m_2(t) := \int_{\mathbb{R}^2} |\mathbf{x}|^2\rho(t,\mathbf{x})\mathrm{d}\mathbf{x}.
\end{equation}
In other words, \(m_2(t)\) grows linearly for \(M_0 < 8\pi\), and decays linearly for \(M_0 > 8\pi\). 
For \(M_0 > 8\pi\), we conclude that \(m_2(t)\) should become negative in finite time, which is impossible since \(\rho\) is nonnegative. Therefore, a blow-up occurs before this contradiction arises. For further details on the critical mass analysis in the 2D parabolic-elliptic case, see \cite{collapsing}.

Unlike the simple 2D parabolic-elliptic KS system, where the total mass is the only factor that determines the aggregation behaviors, there are no similar results in more complex cases \cite{he2022fast}. Numerical experiments \cite{SIPF1} demonstrate that in 3D parabolic-parabolic KS systems, the critical mass depends on the initial distribution of \(\rho\).

\subsection{SIPF Method}
In this section, we briefly introduce the Stochastic Interacting Particle Field (SIPF) method developed in \cite{SIPF1}. The SIPF method is a hybrid particle-spectral approach, which approximates the solutions to \eqref{eq:system} through coupled particle dynamics and field equations, avoiding historical dependence in drift computations while maintaining bounded operational complexity.

Throughout this section, the system is restricted to a domain \( \Omega = \left[-\frac{L}{2}, \frac{L}{2}\right]^3 \) with periodic boundary conditions. In the SIPF method for solving the parabolic-parabolic KS system~\eqref{eq:system}, the chemical concentration \(\hat{c}(\mathbf{x}, n\tau) \) is represented in frequency space using a Fourier series expansion, while the particle density \( \hat{\rho}(\mathbf{x}, n\tau) \) is represented in physical space by some moving particles. The temporal domain \([0, T]\) is partitioned into $N_T = \left\lceil \frac{T}{\tau} \right\rceil$ equal intervals of size \( \tau \), where the ceiling function $\lceil x \rceil$ denotes the smallest integer greater than or equal to $x$.

At each time step \(t_n = n\tau \), the concentration is represented as:
\begin{equation}
\label{eq:U}
\hat{c}^{(n)}(\mathbf{x}):=\hat{c}(\mathbf{x}, n\tau) = \sum_{\mathbf{q} \in U} \hat{\alpha}_{\mathbf{q}}^{(n)} e^{i \frac{2\pi}{L} \mathbf{q} \cdot \mathbf{x}},
\end{equation}
where \( \hat{\alpha}_{\mathbf{q}}^{(n)} \) are the time-dependent Fourier coefficients for \(\hat{c}\), and the frequency index set \(U := \left\{-\frac{H}{2}, \ldots, \frac{H}{2}\right\}^3\). The particle density is represented as:
\begin{equation}
\label{eq:def_rho_hat}
    \hat{\rho}^{(n)}(\mathbf{x}) = \frac{M_0}{P} \sum_{p=1}^P \delta\left(\mathbf{x} - \widehat{\mathbf{X}}_p^{(n)}\right),
\end{equation}
where \(\delta(\cdot)\) is the Dirac delta function, \(M_0 = \int_{\Omega} \rho^{(0)}(\mathbf{x})\,\mathrm{d}\mathbf{x}\) is the initial mass, and \((\widehat{\mathbf{X}}_p^{(n)})_{1\leq p \leq P}\) are positions of particles at time $t_n$.

\paragraph{Updating chemical concentration \(\hat{c}\)} 
\label{par:updating_c}
We discretize the governing equation for \(c\) in \eqref{eq:system} using an implicit Euler scheme:
\begin{equation}
\label{eq:c_discrete}
\frac{\epsilon(\hat{c}^{(n)}-\hat{c}^{(n-1)})}{\tau} = \Delta \hat{c}^{(n)} -k^2 \hat{c}^{(n)} +\hat{\rho}^{(n-1)}.
\end{equation}
Applying the Fourier transform to \eqref{eq:c_discrete} yields:
\begin{equation}
\epsilon \frac{\hat{\alpha}_{\mathbf{q}}^{(n)} - \hat{\alpha}_{\mathbf{q}}^{(n-1)}}{\tau} = -\left(\frac{4\pi^2}{L^2} |\mathbf{q}|^2 + k^2\right) \hat{\alpha}_{\mathbf{q}}^{(n)} + \hat{\beta}_{\mathbf{q}}^{(n-1)},
\label{eq:implicit_alpha}
\end{equation}
where the Fourier coefficients for the density \(\hat{\rho}\) at time $t=(n-1)\tau$, denoted by \(\hat{\beta}_{\mathbf{q}}^{(n-1)}\), are computed as:
\begin{equation}
\label{eq:beta_hat_comp}
\begin{aligned}
\hat{\beta}_{\mathbf{q}}^{(n-1)} = \frac{1}{L^3} \int_{\Omega} \hat{\rho}^{(n-1)}(\mathbf{x}) e^{-i \frac{2\pi}{L} \mathbf{q} \cdot \mathbf{x}} \, \mathrm{d}\mathbf{x}=\frac{M_0}{PL^3} \sum_{p=1}^{P} e^{-\frac{2\pi i}{L}\mathbf{q} \cdot \widehat{\mathbf{X}}_p^{(n-1)}}.
\end{aligned}
\end{equation}
The equation \eqref{eq:implicit_alpha} simplifies to the update rule of \(\hat{\alpha}_{\mathbf{q}}^{(n)}\):
\begin{equation}
\hat{\alpha}_{\mathbf{q}}^{(n)} = \frac{1}{1 + \frac{\tau}{\epsilon}\left(\frac{4\pi^2}{L^2} |\mathbf{q}|^2 + k^2\right)} \hat{\alpha}_{\mathbf{q}}^{(n-1)} + \frac{1}{\frac{4\pi^2}{L^2} |\mathbf{q}|^2 + k^2 + \frac{\epsilon}{\tau}} \hat{\beta}_{\mathbf{q}}^{(n-1)}.
\label{eq:alpha_hat_update}
\end{equation}

\paragraph{Updating particle density \(\rho\)} \label{subsec:particle_simulation}
The governing equation for \(\rho\) in \eqref{eq:system} is the Fokker-Planck equation of the following SDE, where \(\mathbf{W}_{t} = \sqrt{2 \mu} \mathbf{B}_{t}\) and \(\mathbf{B}_{t}\) is a standard Brownian motion:
\begin{equation}
\label{eq:SDE}
\mathrm{d}\mathbf{X} = \chi \nabla_{\mathbf{X}} c (\mathbf{X},t)\mathrm{d}t + \mathrm{d}\mathbf{W}_{t}.
\end{equation}
We apply the Euler-Maruyama scheme:  
\begin{equation}
\widehat{\mathbf{X}}_p^{(n)} = \widehat{\mathbf{X}}_p^{(n-1)} + \chi \tau \nabla \hat{c}^{(n-1)}(\widehat{\mathbf{X}}_p^{(n-1)}) + \mathbf{W}_p^{(n)},
\label{eq:rho_hat_update}
\end{equation}
where \(\mathbf{W}_p^{(n)} \sim \mathcal{N}(0,\, 2\mu\tau \mathbf{I}_3) \) are independent 3D Wiener increments, \( \mathbf{I}_3 \) is the \( 3 \times 3 \) identity matrix, and the drift \(\nabla \hat{c}^{(n-1)}(\widehat{\mathbf{X}}_p^{(n-1)})\) is computed by
\begin{equation}
\label{eq:gradc}
\nabla \hat{c}^{(n-1)}(\widehat{\mathbf{X}}_p^{(n-1)}) = \sum_{\mathbf{q} \in U}  \mathbf{q} \frac{2\pi i}{L} \hat{\alpha}^{(n-1)}_{\mathbf{q}} e^{\frac{2\pi i}{L} \mathbf{q} \cdot \widehat{\mathbf{X}}_p^{(n-1)}} .
\end{equation}

We summarize the SIPF method as Algorithm \ref{Alg:SIPF}. We can find that the computational complexity per timestep is $O\left( P H^3 \right)$, where operations \eqref{eq:beta_hat_comp} and \eqref{eq:gradc} are the bottlenecks.

\begin{algorithm}
\caption{The SIPF method developed in \cite{SIPF1}}\label{Alg:SIPF}
\begin{algorithmic}[1]
\State \textbf{Initialization:}
\State \(\hat{\alpha}^{(0)} \gets \mathbf{0}_{H \times H \times H}\). \Comment{Initialize 3D zero tensor}
\State Sample \(\widehat{\mathbf{X}}_p^{(0)} \sim \text{i.i.d. } \mathbf{X}^{(0)}\). \Comment{With density \(\frac{\rho}{M_0}\) at \(t=0\)}
\Statex
\State \textbf{Update:}
\For{\(n = 1, 2, \dots,N_T\)}
    \State \textbf{Update \(\hat{\alpha}^{(n)}_{H \times H \times H}\):}
    \State \(\hat{\alpha}_{\mathbf{q}}^{(n)} \gets \frac{1}{1 + \frac{\tau}{\epsilon}\left(\frac{4\pi^2}{L^2}|\mathbf{q}|^2 + k^2\right)} \hat{\alpha}_{\mathbf{q}}^{(n-1)} + \frac{1}{\frac{4\pi^2}{L^2}|\mathbf{q}|^2 + k^2 + \frac{\epsilon}{\tau}} \hat{\beta}_{\mathbf{q}}^{(n-1)}\),
    \State where \(\hat{\beta}_{\mathbf{q}}^{(n-1)} = \frac{M_0}{PL^3} \sum_{p=1}^{P} e^{-\frac{2\pi i}{L}\mathbf{q} \cdot \widehat{\mathbf{X}}_p^{(n-1)}}\).
    \Statex
    \State \textbf{Update \(\widehat{\mathbf{X}}_p^{(n)}\):}
    \For{\(p = 1\) \textbf{to} \(P\)}
        \State Compute \(\mathbf{Z}_p^{(n)} \gets \widehat{\mathbf{X}}_p^{(n-1)} + \chi \tau \nabla \hat{c}^{(n-1)}(\widehat{\mathbf{X}}_p^{(n-1)})\).
        \State Sample \(\mathbf{W}_p^{(n)} \sim \mathcal{N}(0, 2\mu \tau \mathbf{I}_3)\). \Comment{Independent noise}
        \State Update \(\widehat{\mathbf{X}}_p^{(n)} \gets \mathbf{Z}_p^{(n)} + \mathbf{W}_p^{(n)}\).
    \EndFor
    \State \hspace{0.5cm} where \(\nabla \hat{c}^{(n-1)}(\widehat{\mathbf{X}}_p^{(n-1)}) = \sum_{\mathbf{q} \in U}  \mathbf{q} \frac{2\pi i}{L} \hat{\alpha}^{(n-1)}_{ \mathbf{q}} e^{\frac{2\pi i}{L} \mathbf{q} \cdot \widehat{\mathbf{X}}_p^{(n-1)}}\). 
\EndFor
\end{algorithmic}
\end{algorithm}

\section{SIPF-PIC Method}
Building on this computational bottleneck, we now introduce particle-in-cell operators as a pathway to a more scalable algorithm. The core idea is to leverage a structured Eulerian representation to reduce the complexity of particle-frequency interactions in the SIPF method, while retaining the Lagrangian flexibility of particles for transport and sampling.

\textcolor{black}{
In this section, we describe the algorithmic details of the proposed method and analyze its key properties, including computational complexity, the errors associated with the introduced particle-in-cell operators, the mass conservation property, and the positivity-preserving property afforded by the newly introduced operator-splitting framework.}
\label{sec:PIC}
\subsection{Particle-in-cell interpolation framework}\label{subsec:interpolation}
The drawback of the original SIPF method (Algorithm~\ref{Alg:SIPF}) is the particle-frequency interaction with $O(PH^3)$ computational complexity. Our new framework is designed to break through this computational bottleneck.

The efficient interaction between the physical domain and the frequency domain is based on the FFT. We interpolate the particles using only a few nearby grid points, while the grid as a whole interacts with the frequency domain. This approach avoids the high complexity of interactions between each particle and each frequency.

In FFT, we require physical domain and frequency domain grids of the same size. We define the physical grid as \(U \times \frac{L}{H} = \{ -\frac{L}{2}, \ldots, \frac{L}{2} - \frac{L}{H}\}^3\), where the index set \(U\) is described in \eqref{eq:U}.

Consider a smooth target function \( f: \mathbb{R}^3 \to \mathbb{R} \); \(\mathbf{X} \in \Omega\) is an arbitrary coordinate. A conventional particle-in-cell (PIC) method \cite{alma99} represents \(f(\mathbf{X})\) as a combination of the function values at grid points:
\begin{equation}
\label{eq:full}
f(\mathbf{X}) \approx \sum_{\mathbf{q} \in U} R\left(\mathbf{X}, \mathbf{q}\frac{L}{H}\right)f\left(\mathbf{q}\frac{L}{H}\right),
\end{equation}
where the particle-grid interaction kernel $R\left(\mathbf{X}, \mathbf{q}\frac{L}{H}\right)$ satisfies \(\sum_{\mathbf{q}\in U}R\left(\mathbf{X}, \mathbf{q}\frac{L}{H}\right)=1\) for every $\mathbf{X} \in \Omega$. Furthermore, the kernel \(R(\mathbf{X},\cdot)\) is locally supported, meaning that for a given $\mathbf{X} \in \Omega$, only a few nearby grid points contribute to the approximation.

For example, in a trilinear interpolation with a second-order truncation error, the kernel is represented as
\begin{equation}
R\left(\mathbf{X}, \mathbf{q}\frac{L}{H}\right) = \prod_{1\leq j \leq 3}\max\left(1-\frac{H}{L}\left|X_j - q_j \frac{L}{H}\right|,0\right).
\end{equation}
In particular, $R(\mathbf{X},\cdot)$ is nonzero only for the eight vertices of the grid cell containing $\mathbf{X}$. We denote the second-order basis \(\Gamma_2 = \{0,1\}^3\) and the relative grid positions \(\mathbf{X}_{{\gamma}} = \left( \left\lfloor \frac{H}{L} \mathbf{X} \right\rfloor + {\gamma} \right) \times \frac{L}{H}\), where the notation $\lfloor \cdot \rfloor$ denotes the floor function applied component-wise. The interpolation \eqref{eq:full} can be compactly expressed as
\begin{equation}
\label{eq:2ord}
f(\mathbf{X}) \approx \sum_{\boldsymbol{\gamma} \in \Gamma_2} F_{{\gamma}}^{(2)} (\boldsymbol{\lambda})f(\mathbf{X}_{{\gamma}}),
\end{equation}
where the weights \(F_{j_1 j_2 j_3}^{(2)}(\boldsymbol{\lambda}) = \prod_{s=1}^{3} \left( 1 - j_s + (2j_s - 1) \lambda_s \right)\) for \((j_1,j_2,j_3)\in \Gamma_2\), and $\boldsymbol{\lambda}=\boldsymbol{\lambda}(\mathbf{X}):=\frac{H}{L}\mathbf{X}-\left\lfloor \frac{H}{L} \mathbf{X} \right\rfloor \in [0,1)^3$ is the relative location of $\mathbf{X}$ in the grid cell containing $\mathbf{X}$. For points near the boundary under periodic boundary conditions, the indices are wrapped around the grid edges.

\begin{figure}[h]
    \centering
    \includegraphics[width=0.42\linewidth]{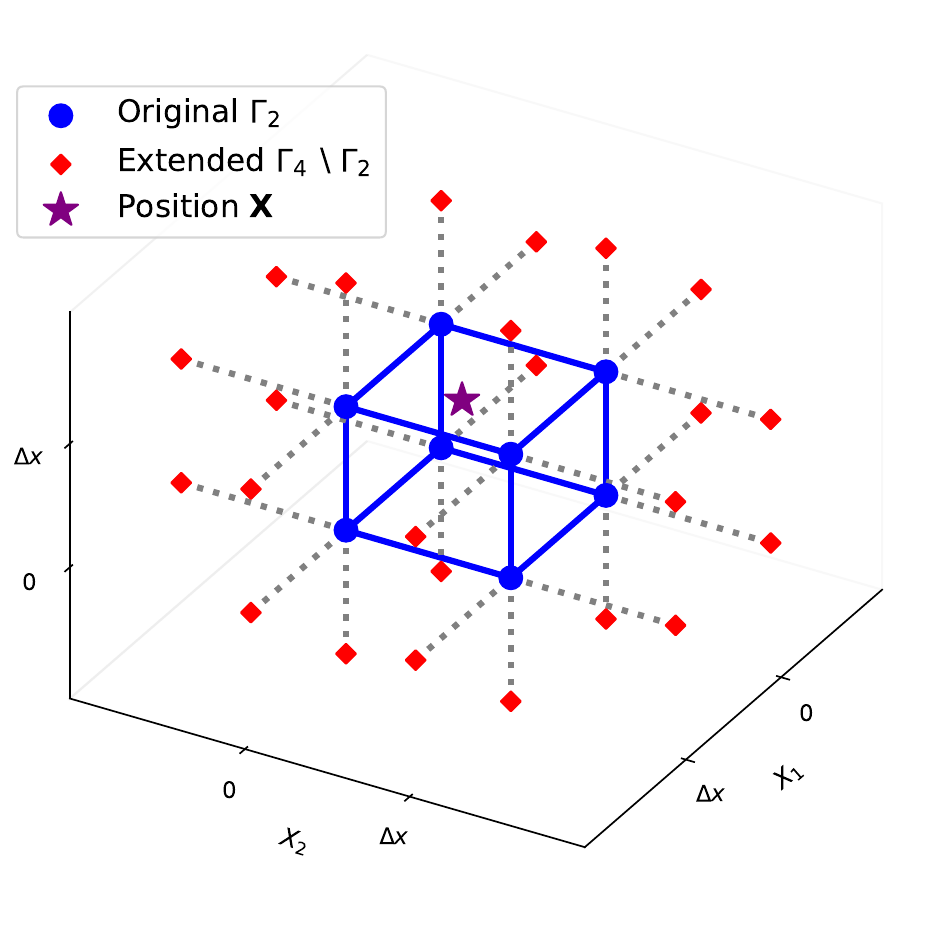}
    \caption{Extended vertices for fourth-order.}
    \label{Extended}
\end{figure}

To achieve higher precision, we also introduce the fourth-order basis\[\Gamma_4 = \Gamma_2 \cup \{ -1,2 \} \times \{0,1\} \times \{0,1\} \cup \{0,1\} \times \{-1,2\} \times \{0,1\} \cup \{0,1\} \times \{0,1\} \times \{-1,2\},\]
which enables interpolation \cite{tricubic} exhibiting a fourth-order truncation error, as shown in Figure~\ref{Extended}.
The complete fourth-order interpolation formula becomes as follows:
\begin{equation}
\label{eq:4ord}
f(\mathbf{X}) \approx \sum_{\gamma \in \Gamma_4} F_{\gamma}^{(4)} (\boldsymbol{\lambda})f(\mathbf{X}_{\gamma}),
\end{equation}
where the weights \(F_\gamma^{(4)}(\boldsymbol{\lambda})\) are shown below, exhibiting cubic symmetry:
\begin{equation}
\begin{aligned}
&F_{2,j_2,j_3}^{(4)}(\boldsymbol{\lambda}) = -\frac{(1 + \lambda_1)(1 - \lambda_1) \lambda_1 \left( 1 - j_2 + (2j_2 - 1) \lambda_2 \right)\left( 1 - j_3 + (2j_3 - 1) \lambda_3 \right)}{6},\\
&F_{-1,j_2,j_3}^{(4)}(\boldsymbol{\lambda}) = -\frac{(2 - \lambda_1)\lambda_1(1 - \lambda_1)\left( 1 - j_2 + (2j_2 - 1) \lambda_2 \right)\left( 1 - j_3 + (2j_3 - 1) \lambda_3 \right)}{6},\\
&\cdots\\
&F_{j_1j_2j_3}^{(4)}(\boldsymbol{\lambda}) = \prod_{s=1}^{3} \left( 1 - j_s + (2j_s - 1) \lambda_s \right) \left( 1 + \sum_{k=1}^{3} \frac{\lambda_k (1 - \lambda_k)}{2} \right),(j_1,j_2,j_3)\in \{0,1\}^3.
\end{aligned}
\end{equation}

We define particle-to-grid discretization as follows.
\begin{definition}
\label{def:particle2grid}
Note the set of weighted empirical measures \(\mathcal{M}_D(\Omega):=\\ \left\{\sum_{p=1}^{P} w_p \delta(\mathbf{x}-\widehat{\mathbf{X}}_p):P \in \mathbb{N},w_p \in \mathbb{R}^+,\widehat{\mathbf{X}}_p \in \Omega\right\}\). The discretization operator \( I:\\ \mathcal{M}_D(\Omega) \to \mathbb{R}^{H\times H\times H} \) maps each empirical measure onto a three-dimensional grid partitioning the spatial domain $\Omega$ into $H^3$ equal cells, defined as:
\begin{equation}
    I(\sum_{p=1}^{P} w_p \delta(\mathbf{x}-\widehat{\mathbf{X}}_p)) = \sum_{p=1}^{P} w_p \sum_{\gamma \in \Gamma} F_{p,\gamma} \delta(\mathbf{x}-\widehat{\mathbf{X}}_{p,\gamma}),
\end{equation}
where \( \Gamma \) is the interpolation basis index set from Section~\ref{subsec:interpolation}, \( F_{p,\gamma}:=F_\gamma (\boldsymbol{\lambda}(\mathbf{X}_p))\) and \(\widehat{\mathbf{X}}_{p,\gamma}:=\left( \left\lfloor \frac{H}{L} \mathbf{X}_p \right\rfloor + {\gamma} \right) \times \frac{L}{H}\) are the corresponding coefficients and grid points in \eqref{eq:2ord} or \eqref{eq:4ord} for \(\widehat{\mathbf{X}}_p\).
\end{definition}


Similarly, we define grid-to-particle interpolation as follows.
\begin{definition}
\label{def:grid2particle}
The interpolation operator \( I^*: \mathbb{R}^{H\times H\times H} \to C(\Omega) \) constructs continuous approximations via:
\begin{equation}
    (I^* f)(\mathbf{X}) = \sum_{\gamma \in \Gamma} F_\gamma (\boldsymbol{\lambda}(\mathbf{X}))f(\mathbf{X}_\gamma),
\end{equation} 
where \( f \) is sampled at grid points \( \mathbf{X}_\gamma \).
\end{definition}

The operator \(I^*\) produces piecewise polynomial functions that depend only on the nodal values, satisfying the adjoint relation \( \langle I\rho, f\rangle = \langle \rho, I^*f\rangle \) for \(\rho \in \mathcal{M}(\Omega)\) and \(f\) bounded.
\begin{color}{black}
\begin{remark}
\label{remark:smoothness}
Since the particle-to-grid operator is applied to the relatively smooth trigonometric function
\(e^{\frac{2\pi i}{L}\mathbf{q}\cdot\mathbf{x}}\),
whereas the grid-to-particle operator is applied to the less smooth quantity $\nabla \hat{c}(\mathbf{x})$, we employ the fourth-order particle-to-grid operator and the second-order grid-to-particle operator. This combination is also adopted in the error analysis presented in Section~\ref{sec:Error}.

Since fourth-order interpolation inevitably involves negative interpolation weights, a smoothing procedure near the density boundary \cite{doi:10.1137/16M105962X,NATARAJAN2021114001} is required to guarantee the positivity-preserving property stated in Section~\ref{subsec:os}. This property is not needed in the error analysis of Section~\ref{sec:Error}; therefore, the smoothing procedure is omitted there.    
\end{remark}
\end{color}

For completeness, we consider the following convention of the discrete Fourier transform
\(\mathcal{F}_H\), which admits the following definition:
\begin{equation}
\mathcal{F}_H f\left(\mathbf{y}\right) := L^{-3} \int_{\mathbf{x} \in \Omega} e^{-i \mathbf{x} \cdot \mathbf{y}} f(\mathbf{x}) \, \mathrm{d}\mathbf{x},
\end{equation}
where \(\mathbf{y} = \mathbf{q} \frac{2\pi}{L}\) is a frequency vector.
The inverse Fourier transform operator \(\mathcal{F}_H^{-1}\) is then given by:
\begin{equation}
\mathcal{F}_H^{-1} \hat{\alpha}(\mathbf{x}) = \sum_{\mathbf{q} \in U} \hat{\alpha}_{\mathbf{q}} e^{i \frac{2\pi}{L} \mathbf{q} \cdot \mathbf{x}}.
\end{equation}

The Fourier coefficients \(\hat{\alpha}_{\mathbf{q}}\) and \(\hat{\beta}_{\mathbf{q}}\) defined in \eqref{eq:U} and \eqref{eq:beta_hat_comp} are specific evaluations:
\(
\hat{\alpha}_{\mathbf{q}} = \mathcal{F}_H \hat{c}\left(\mathbf{q}\frac{2\pi}{L}\right), \quad 
\hat{\beta}_{\mathbf{q}} = \mathcal{F}_H \hat{\rho}\left(\mathbf{q}\frac{2\pi}{L}\right).
\)
The composition \(\mathcal{F}_H^{-1} \circ \mathcal{F}_H\) forms a projection onto the space of trigonometric polynomials \(\mathcal{T}_H:= \left\{ f : f = \sum_{\mathbf{q} \in U} \hat{\alpha}_{\mathbf{q}} e^{i \frac{2\pi}{L} \mathbf{q} \cdot \mathbf{x}} \right\}\). For any function \( f \in \mathcal{T}_H \), it holds that
\(
\mathcal{F}_H^{-1} \circ \mathcal{F}_H (f) = f.
\)
Additionally, for \(\mathbf{q}, \mathbf{q'} \in U\) with \(\mathbf{q} \neq \mathbf{q'}\), the following properties are satisfied:
\begin{equation}
  \begin{aligned}
    \mathcal{F}_H^{-1} \circ \mathcal{F}_H \left(\delta \left(\mathbf{x}-\mathbf{q} \frac{L}{H} \right) \right) \left(\mathbf{q} \frac{L}{H} \right) &= \left( \frac{H}{L} \right)^3, \\
    \mathcal{F}_H^{-1} \circ \mathcal{F}_H \left(\delta \left(\mathbf{x}-\mathbf{q} \frac{L}{H} \right) \right) \left(\mathbf{q'} \frac{L}{H} \right) &= 0.
  \end{aligned}
\end{equation}

\subsection{Computational complexity and algorithmic details}
\label{subsec:complexity}
Consider \( P \)\\ equally weighted particles at positions \( \{\widehat{\mathbf{X}}_p\}_{p=1}^P \subset \Omega\) with total mass $M_0$. The empirical measure is shown as:
\begin{equation}
    \label{eq:emp_measure}
    \hat{\rho}(\mathbf{x}) = \frac{M_0}{P}\sum_{p=1}^P \delta(\mathbf{x} - \widehat{\mathbf{X}}_p).
\end{equation}
We set the computational grid \( \mathcal{L} = \{\mathbf{X}_\gamma = \frac{L}{H} \mathbf{q} : \mathbf{q}  \in U\} \) with spacing \( h = \frac{L}{H} \). We perform Algorithm \ref{Alg:P2G} to compute \( \mathcal{F}_H I(\hat{\rho})(\mathbf{y}) \) for all \( \mathbf{y} \in U \times \frac{2\pi}{L} \), where the discretization operator $I$ is specified in Definition \ref{def:particle2grid}.

\begin{algorithm}
\caption{Particle-to-Grid discretization}
\label{Alg:P2G}
\begin{algorithmic}[1]
    \State Initialize zero tensor \(\phi[\mathbf{q}]:= \phi_\mathbf{q}\) for all \(\mathbf{q} \in U\)
    \For{\(p = 1\) \textbf{to} \(P\)}
        \For{\(\gamma \in \Gamma\)}
        \Comment{Defined in Section~\ref{subsec:interpolation}}
            \State Find grid index: \(\mathbf{q}_{p,\gamma} = \mathbf{X}_{p,\gamma} \cdot \frac{H}{L} \) 
            \State Calculate interpolation weight: \(F_{p,\gamma}\)
            \State Update density: \(\phi[\mathbf{q}_{p,\gamma}] \leftarrow \phi[\mathbf{q}_{p,\gamma}] + \frac{M_0 F_{p,\gamma}}{P}\)
        \EndFor
    \EndFor
    \State Compute \(\mathcal{F}_H \phi\) via Fast Fourier Transform
\end{algorithmic}
\end{algorithm}
\begin{lemma} With $H$ Fourier modes and $P$ particles, 
the complexity of particle-to-grid discretization in Algorithm  \ref{Alg:P2G} is \( O(H^3 \log H + P) \). 
\end{lemma}
\begin{proof}
The particle-to-grid transfer costs \( O(P) \) (each particle contributes to \(|\Gamma|\) neighbors), and the FFT on the \(H \times H \times H\) grid costs \( O(H^3 \log H) \).
\end{proof}

For the inverse operation, given the values of a function \( f \) on the grid \( f\left(U \cdot \frac{L}{H}\right) \), we perform Algorithm \ref{Alg:G2P} to compute \( I^* f(\mathbf{X}_1), \ldots, I^* f(\mathbf{X}_P) \), where the interpolation operator $I^*$ is specified in Definition \ref{def:grid2particle}.
\begin{algorithm}
\caption{Grid-to-Particle interpolation}
\label{Alg:G2P}
\begin{algorithmic}[1]
    \State Input tensor \(f[\mathbf{q}]:= f_{\mathbf{q}}\) defined on grid indices \(\mathbf{q} \in U\)
    \For{\(p = 1\) \textbf{to} \(P\)}
        \State Initialize particle value: \(f_p \leftarrow 0\)
        \For{\(\gamma \in \Gamma\)} \Comment{Defined in Section~\ref{subsec:interpolation}}
            \State Find grid index: \(\mathbf{q}_{p,\gamma} = \mathbf{X}_{p,\gamma} \cdot \frac{H}{L}\)
            \State Obtain interpolation weight: \(F_{p,\gamma}\)
            \State Accumulate contribution: \(f_p \leftarrow f_p + F_{p,\gamma} \cdot f[\mathbf{q}_{p,\gamma}]\)
        \EndFor
        \State Store interpolated value: \(I^*  f(\mathbf{X}_p) \leftarrow f_p\)
    \EndFor
\end{algorithmic}
\end{algorithm}

\begin{lemma}
With $P$ particles, the complexity of the grid-to-particle interpolation in Algorithm \ref{Alg:G2P} is \( O(P) \).
\end{lemma}
\begin{proof}
The grid-to-particle transfer costs \( O(P) \) (each particle interacts with a fixed number of \(|\Gamma|\) neighbors).
\end{proof}

The computational complexity of the original method, dominated by all-to-all particle-frequency interactions scaling as \(O(P H^3)\), is reduced to \(O(P + H^3 \log H)\) via the PIC acceleration presented in Section \ref{subsec:interpolation}. This acceleration involves local grid-projected density interpolation and gradient reconstruction, limiting particle interactions to neighboring grid points while handling frequency components globally. Periodic boundary conditions are imposed on the boundary grids. \begin{color}{black} To guarantee the numerical convergence of Algorithm~\ref{SIPF-PIC} (see Section~\ref{sec:Error}), we apply the Gaussian spectral filter $\Phi_{\mathbf q}=\exp\left(-\frac{2\pi ^2|\mathbf q|^2}{H_0^2}\right)$ to the Fourier coefficients. This operation is equivalent to convolving the particle distribution with the Gaussian kernel $\mathcal N\left(0,\frac{L^2}{H_0^2}\mathbf I_3\right)$ where \( \mathbf{I}_3 \) is the \( 3 \times 3 \) identity matrix, which smooths the distribution over the length scale corresponding to a grid of resolution $H_0$.

For the particle update over a single time step, we employ a symmetric
half-diffusion--advection--half-diffusion splitting. This choice is motivated
by the error analysis in Section~\ref{sec:Error}, which requires the pre-diffusion step in
order to close the Grönwall estimate pathwise. Applying the same symmetric
splitting to the concentration field may also improve the temporal
convergence order. In the present work, however, we adopt this formulation
primarily for theoretical stability rather than for temporal order
enhancement.
\end{color}

\begin{algorithm}
\caption{SIPF-PIC method}\label{SIPF-PIC}
\begin{algorithmic}[1]
\State \textbf{Initialization:}
\State $\hat{\alpha}_{\mathbf{q}}^{(0)} \gets \mathbf{0}_{H \times H \times H}$ \Comment{3D zero tensor (indexed by $\mathbf{q} \in U$)}
\State Sample $\widehat{\mathbf{X}}_p^{(0)} \sim \text{i.i.d. } \mathbf{X}^{(0)}$
\Statex
\State \textbf{Update:}
\For{$n = 1, 2, \dots,N_T$}
    \State \textbf{Update $\hat{\alpha}_{\mathbf{q}}^{(n)}$:}
    \State $\hat{\alpha}_{\mathbf{q}}^{(n)} \gets \frac{1}{1 + \frac{\tau}{\epsilon}\left(\frac{4\pi^2}{L^2}|\mathbf{q}|^2 + k^2\right)} \hat{\alpha}_{\mathbf{q}}^{(n-1)} + \frac{1}{\frac{4\pi^2}{L^2}|\mathbf{q}|^2 + k^2 + \frac{\epsilon}{\tau}} \mathcal{F}_H I(\hat{\rho}^{(n-1)})\left(\mathbf{q} \frac{2\pi}{L}\right)$
    \State where $\hat{\rho}^{(n-1)}(\mathbf{x}) = \sum_{p=1}^{P} \frac{M_0}{P} \delta(\mathbf{x} - \widehat{\mathbf{X}}_p^{(n-1)})$
    \State and $\mathcal{F}_H I(\hat{\rho}^{(n-1)})\left(\mathbf{q} \frac{2\pi}{L}\right)$ computed by Alg. \ref{Alg:P2G}
    \Statex
    \State \textbf{Update $\widehat{\mathbf{X}}_p^{(n)}$:}
    \For{$p = 1$ \textbf{to} $P$}
        \State Sample $\mathbf{W}_p^{-(n)},\mathbf{W}_p^{+(n)} \sim \mathcal{N}(0, \mu \tau \mathbf{I}_3)$ \Comment{Independent 3D increments}
        \State Compute $\mathbf{Z}_p^{(n)} \gets \widehat{\mathbf{X}}_p^{(n-1)}+\mathbf{W}_p^{-(n)} + \chi \tau I^*(\nabla \hat{c}^{(n-1)})(\widehat{\mathbf{X}}_p^{(n-1)}+\mathbf{W}_p^{-(n)})$
        \State Update $\widehat{\mathbf{X}}_p^{(n)} \gets \mathbf{Z}_p^{(n)} + \mathbf{W}_p^{+(n)}$ \Comment{Periodic setting}
    \EndFor
    \State with $\nabla \hat{c}^{(n-1)}\left( \mathbf{p} \frac{L}{H} \right) = \mathcal{F}_H^{-1}\left(\frac{2\pi i}{L}\mathbf{q} \hat{\alpha}_{\mathbf{q}}^{(n-1)} \Phi_\mathbf{q}\right)\left( \mathbf{p} \frac{L}{H} \right)$
    \Statex \Comment{Physical grid values computed all at once via FFT}
    \State where $\Phi_{\mathbf q}=\exp\left(-\frac{2\pi ^2|\mathbf q|^2}{H_0^2}\right)$\Comment{Gaussian spectral filter}
    \State and $I^*(\nabla \hat{c}^{(n-1)})$ computed by Alg. \ref{Alg:G2P}
\EndFor
\end{algorithmic}
\end{algorithm}

\begin{theorem}
\label{thm:overall_complexity}
The single-step complexity of Algorithm~\ref{SIPF-PIC} is \( O \left(H^3 \log H + P \right) \).
\end{theorem}
\begin{proof}
 Each iteration costs \( O(H^3 \log H + P) \): spectral updates require FFT with complexity \( O(H^3 \log H) \) and particle-to-grid transfer with \( O(P) \), while particle updates require inverse FFT with \( O(H^3 \log H) \) and gradient operations with \( O(P) \).
\end{proof}

\subsection{Bounds of interpolation error}

Recall the empirical measure \eqref{eq:emp_measure} and the discretization operator $I$ defined in Definition \ref{def:particle2grid}. We claim the following error bounds for $I$ under a Fourier transform.

\begin{lemma}
\label{Thm:P2G}
For any particle distribution \(\{\mathbf{X}_p\}_{p=1}^P \subset \Omega\), the empirical density \(\hat{\rho}^{(n-1)}\) is defined by \eqref{eq:def_rho_hat}. let the frequencies be \( \mathbf{y} = \mathbf{q} \frac{2\pi}{L} \) for indices \( \mathbf{q} \in U\). The interpolation errors satisfy:
\begin{align}
|(\mathcal{F}_H I_2(\hat{\rho}) - \mathcal{F}_H \hat{\rho})(\mathbf{y})| &\leq M_0 C_2 L^{-3} |\mathbf{y}|^2 \left(\frac{L}{H}\right)^2, \\
|(\mathcal{F}_H I_4(\hat{\rho}) - \mathcal{F}_H \hat{\rho})(\mathbf{y})| &\leq M_0 C_4 L^{-3} |\mathbf{y}|^4 \left(\frac{L}{H}\right)^4,
\end{align}
where \( C_{\eta} \) depends on the interpolation order \( \eta \) with $\eta=2, 4$.
\end{lemma}

\begin{proof}
Let \( g_\mathbf{y}(\mathbf{x}) = e^{-i\mathbf{y} \cdot \mathbf{x}} \). Then we have the Fourier transform 
\begin{equation}
\begin{aligned}
\label{eq:FH}
\mathcal{F}_H \hat{\rho}(\mathbf{y}) = \frac{M_0}{PL^3} \sum_{p=1}^{P} g_\mathbf{y}(\mathbf{X}_p),\quad
\mathcal{F}_H I(\hat{\rho})(\mathbf{y}) = \frac{M_0}{PL^3} \sum_{p=1}^{P} \sum_{\gamma \in \Gamma} F_{p,\gamma} g_\mathbf{y}(\mathbf{X}_{p, \gamma}).
\end{aligned}
\end{equation}
 Using a \(d\)-th order Taylor expansion with Lagrange remainder at \(\mathbf{X}\), for multi-index \(\mathbf{k} = (k_1, k_2, k_3)\) and \(D_{\mathbf{k}} = \partial^{|\mathbf{k}|}/\partial x_1^{k_1} \partial x_2^{k_2} \partial x_3^{k_3}\), we have
\begin{equation}
g_\mathbf{y}(\mathbf{x}) = \sum_{|\mathbf{k}| \leq \eta-1} \frac{D_{\mathbf{k}} g_\mathbf{y}(\mathbf{X})}{\mathbf{k}!} (\mathbf{x} - \mathbf{X})^{\mathbf{k}} + \sum_{|\mathbf{k}| = \eta} \frac{D_{\mathbf{k}} g_\mathbf{y}(\mathbf{X} + \xi)}{\mathbf{k}!} (\mathbf{x} - \mathbf{X})^{\mathbf{k}},
\end{equation}
where \(\xi\) lies between \(\mathbf{X}\) and \(\mathbf{x}\), and \(\mathbf{k}! = k_1! k_2! k_3!\).
Throughout this paper, the notation $|\mathbf{k}|=k_1+k_2+k_3$ means the order of a multi-index $\mathbf{k}$, which is different from the norm $|\mathbf{x}|$ for a vector $\mathbf{x}$ in Euclidean space.

The exactness of \(I_\eta\) for polynomials of degrees \(0\) to \( \eta-1\) implies cancellation of lower-order terms, in other words,
\(
\sum_{\gamma \in \Gamma} \overline{F_{\gamma}} f(\mathbf{X}_\gamma) = f(\mathbf{X})
\)
holds for all polynomials \( f \) of degree at most \( \eta-1 \), where $\overline{F_{\gamma}}:=F_{\gamma}(\boldsymbol{\lambda}(\mathbf{X}))$, yielding: 
\begin{equation}
-g_\mathbf{y}(\mathbf{X}) + \sum_{\gamma \in \Gamma} \overline{F_{\gamma}} g_\mathbf{y}(\mathbf{X}_\gamma) = \sum_{|\mathbf{k}| = \eta} \sum_{\gamma \in \Gamma} \overline{F_{\gamma}} \frac{D_{\mathbf{k}} g_\mathbf{y}(\mathbf{X} + \xi(\mathbf{X}, \gamma))}{\mathbf{k}!} (\mathbf{X}_\gamma - \mathbf{X})^{\mathbf{k}}.
\end{equation}

Since \( g_\mathbf{y}(\mathbf{x}) = e^{-i\mathbf{y} \cdot \mathbf{x}} \), we have 
\(
|D_{\mathbf{k}} g_\mathbf{y}(\mathbf{X} + \xi(\mathbf{X}, \gamma))| \leq |\mathbf{y}|^\eta,
\)
for all \(\mathbf{k},\mathbf{X},\xi, \gamma\).

For the case of \( \eta = 2 \) and \( \Gamma =\Gamma_2 \), we have \( \overline{F_{\gamma}} \geq 0 \), \( \sum_{\gamma \in \Gamma} \overline{F_{\gamma}} = 1 \) and
\(
(\mathbf{X}_\gamma - \mathbf{X})^{\mathbf{k}} \leq \left( \frac{L}{H} \right)^2.
\)
Thus, we establish the second-order error bound
\begin{equation}
    \left|-g_\mathbf{y}(\mathbf{X}) + \sum_{\gamma \in \Gamma} \overline{F_{\gamma}} g_\mathbf{y}(\mathbf{X}_\gamma)\right| \leq \sum_{|\mathbf{k}| = 2} \frac{|\mathbf{y}|^2}{\mathbf{k}!} \left( \frac{L}{H} \right)^2 = \frac{9}{2} |\mathbf{y}|^2 \left( \frac{L}{H} \right)^2.
    \label{eq:PIC2}
\end{equation}
For the case of \( \eta = 4 \), let \(\Gamma = \Gamma_4\). The coefficients in the outer set \( \Gamma \setminus \Gamma_2 \) are negative, while the coefficients in the inner set \( \Gamma_2 \) are positive.

We only need to handle the following
\begin{equation}
\begin{aligned}
\left| \sum_{\gamma \in \Gamma} \overline{F_{\gamma}} (\mathbf{X}_\gamma - \mathbf{X})^{\mathbf{k}} \right| 
&\leq \left| \sum_{\gamma \in \Gamma_2} \overline{F_{\gamma}} (\mathbf{X}_\gamma - \mathbf{X})^{\mathbf{k}} \right| 
   + \left| \sum_{\gamma \in \Gamma \setminus \Gamma_2} \overline{F_{\gamma}} (\mathbf{X}_\gamma - \mathbf{X})^{\mathbf{k}} \right| \\
&\leq \left( \frac{L}{H} \right)^4 \left( \sum_{\gamma \in \Gamma_2} \overline{F_{\gamma}} 
   -2^4 \sum_{\gamma \in \Gamma \setminus \Gamma_2} \overline{F_{\gamma}} \right) \\
&= \left( \frac{L}{H} \right)^4 \left( 1 -17 \sum_{\gamma \in \Gamma \setminus \Gamma_2} \overline{F_{\gamma}} \right).
\end{aligned}
\end{equation}
We partition the 24 outer grid points into 8 clusters based on their nearest inner vertices. For instance, the cluster containing \( \mathbf{X}_{211}, \mathbf{X}_{121}, \) and \( \mathbf{X}_{112} \) shares adjacency with the inner vertex \( \mathbf{X}_{111} \).
For each cluster, the coefficient combination satisfies
\begin{equation}
-(\overline{F}_{211} + \overline{F}_{121} + \overline{F}_{112}) = \lambda_1 \lambda_2 \lambda_3 \left( \frac{3 - \lambda_1^2 - \lambda_2^2 - \lambda_3^2}{6} \right) \leq \frac{1}{2} \lambda_1 \lambda_2 \lambda_3.
\end{equation}
Aggregating over all 8 clusters yields
\begin{equation}
-\sum_{\gamma \in \Gamma \setminus \Gamma_2} \overline{F_{\gamma}} \leq \frac{1}{2} \left( \lambda_1 \lambda_2 \lambda_3 + \cdots + (1 - \lambda_1)(1 - \lambda_2)(1 - \lambda_3) \right) = \frac{1}{2}. 
\end{equation}
Combining these results, we establish the fourth-order error bound
\begin{equation}
\label{eq:PIC4}
\left|-g_\mathbf{y}(\mathbf{X}) + \sum_{\gamma \in \Gamma} \overline{F_{\gamma}} g_\mathbf{y}(\mathbf{X}_\gamma) \right| \leq \frac{19}{2} \sum_{|\mathbf{k}| = 4} \frac{|\mathbf{y}|^4}{\mathbf{k}!} \left( \frac{L}{H} \right)^4 = \frac{513}{16} |\mathbf{y}|^4 \left( \frac{L}{H} \right)^4.
\end{equation}
The overall spectral error consequently satisfies
\begin{equation}
\begin{aligned}
\left|\mathcal{F}_H I(\hat{\rho})(\mathbf{y}) - \mathcal{F}_H \hat{\rho}(\mathbf{y}) \right| 
&= \left| \frac{M_0}{P L^3} \left( \sum_{p=1}^{P} \sum_{\gamma \in \Gamma} F_{p,\gamma} g_\mathbf{y}(\mathbf{X}_{p, \gamma}) - \sum_{p=1}^{P} g_\mathbf{y}(\mathbf{X}_p) \right) \right| \\
&\leq \frac{M_0}{PL^3} \cdot P \cdot C_\eta |\mathbf{y}|^\eta \left( \frac{L}{H} \right)^\eta \\
&= C_\eta M_0 L^{-3} |\mathbf{y}|^\eta \left( \frac{L}{H} \right)^\eta,
\end{aligned}
\end{equation}
where \( \eta = 2 \) (with \( C_2 = \frac{9}{2} \)) and \( \eta = 4 \) (with \( C_4 = \frac{513}{16} \)) correspond to second-order and fourth-order interpolation, respectively.
\end{proof}

As a direct consequence of Lemma \ref{Thm:P2G} (see equations \eqref{eq:PIC2} and \eqref{eq:PIC4}), we present the error bounds for the operation $I^*$ specified in Definition \ref{def:grid2particle} as follows.

\begin{lemma}[Error bounds of grid-to-particle interpolation]
\label{Thm:G2P}
Assume \(f \in \\C^2(\Omega)\). Given the values of \( f \) on the grid \( f\left(U \cdot \frac{L}{H}\right) \), for each \( p = 1,..., P \), the error of the second-order interpolation scheme at position \( \mathbf{X}_p \) satisfies 
\begin{equation}
|I^*_2 f(\mathbf{X}_p) - f(\mathbf{X}_p)| \leq K_2 \left(\frac{L}{H}\right)^2;
\end{equation}
if we further assume \(f \in C^4(\Omega)\), the error of the fourth-order interpolation scheme satisfies
\begin{equation}
|I^*_4 f(\mathbf{X}_p) - f(\mathbf{X}_p)| \leq K_4 \left(\frac{L}{H}\right)^4,
\end{equation}
where \(K_2 = \frac{9}{2} \max_{|\mathbf{k}|=2} \|D_\mathbf{k} f\|_\infty\) is associated with the second derivative of \( f \), and \(K_4 = \frac{513}{16} \max_{|\mathbf{k}|=4} \|D_\mathbf{k} f\|_\infty\) is associated with the fourth derivative of \( f \).
\end{lemma}
\begin{color}{black}
\subsection{Operator-Splitting Formulation}
\label{subsec:os}
The original formulation of the SIPF method (Algorithm \ref{Alg:SIPF}) is obtained by truncating the Fourier basis functions and therefore lacks a direct interpretation in physical space. We now reformulate the SIPF-PIC method (Algorithm \ref{SIPF-PIC}) from an operator-splitting perspective. Besides providing a natural positivity-preserving property of both $\hat{\rho}$ and $\hat{c}$, this formulation also extends naturally to more general computational domains and boundary conditions.

In Algorithm~\ref{SIPF-PIC}, the Fourier coefficients $\hat{\alpha}$ of the concentration $\hat{c}$ are updated according to
\begin{equation}
\hat{\alpha}_{\mathbf q}^{(n+1)}
=
\hat{\kappa}_{\tau,\mathbf q}
\left(
\hat{\alpha}_{\mathbf q}^{(n)}
+
\tau \hat{\beta}_{\mathbf q}^{(n)}
\right),
\end{equation}
where $\hat{\kappa}_{\tau,\mathbf q}= \exp \left(-\frac{\tau}{\epsilon}\left(\left|\mathbf q\frac{2\pi}{L}\right|^2+k^2\right)\right)$ is the attenuation factor of the Fourier mode $\mathbf q$ over one diffusion step of size $\tau$, and $\hat{\beta}_{\mathbf q}^{(n)}=\mathcal{F}_H\!\left(I\hat{\rho}^{(n)}\right)$ is the Fourier coefficient of the source term. Here, $I$ denotes the discretization operator defined in Definition~\ref{def:particle2grid}.
\begin{remark}
This attenuation factor differs slightly from that obtained by the implicit Euler discretization \(
\kappa_{\tau,\mathbf q}
=
\left(
1+\frac{\tau}{\epsilon}
\left(
\left|\mathbf q\frac{2\pi}{L}\right|^2+k^2
\right)
\right)^{-1}\). Nevertheless, in the presence of the spectral regularization scale
$H_0\sim\tau^{-1/2}$, the two propagators exhibit nearly identical numerical behavior, since their discrepancy is confined to the high-frequency modes that are already exponentially suppressed by the regularization. As the error analysis in Section~\ref{sec:Error} does not rely on strict positivity preservation, we do not pursue this comparison further.    
\end{remark}

Equivalently, in the physical space the update can be written as
\begin{equation}
\hat{c}^{(n+1)}
=
\widehat{G}_{\tau}
*
\left(
\hat{c}^{(n)}
+
\tau I\hat{\rho}^{(n)}
\right),
\end{equation}
where $\widehat{G}_{\tau}=\mathcal{F}_H^{-1}\hat{\kappa}_{\tau}$ is the diffusion kernel with exponential attenuation associated with one time step of size $\tau$. Therefore, provided that both $\widehat{G}_{\tau}$ and $I\hat{\rho}^{(n)}$ are nonnegative everywhere, the update is naturally positivity-preserving.

The nonnegativity of $I\hat{\rho}^{(n)}$ has already been discussed in Remark~\ref{remark:smoothness}. We now turn to the positivity of the diffusion kernel with exponential attenuation $\widehat{G}_{\tau}=\mathcal{F}_H^{-1}\hat{\kappa}_{\tau}$. The following theorem shows that when the diffusion length within one time step is larger than the grid spacing, namely $\tau \gtrsim \left(\frac{L}{H}\right)^{2-\sigma}$ for some $\sigma>0$, there exists a nonnegative kernel $\widetilde{G}_{\tau}$ that differs from $\widehat{G}_{\tau}$ only by an exponentially small amount. Consequently, replacing $\widehat{G}_{\tau}$ by $\widetilde{G}_{\tau}$ introduces an error that is exponentially small compared with the polynomial discretization error of the algorithm. Similar exponential tail bounds for Fourier-grid approximations of Gaussian kernels have been derived by Barnett, Greengard, and Rachh~\cite{uniform}, where the truncation error of the squared-exponential kernel is shown to decay exponentially with the frequency cut-off.
\begin{theorem}
\label{theorem:positive_diffusion_kernel}
There exists a kernel
$\widetilde{G}_{\tau}\in\mathbb{R}^{U_H}$ whose components are strictly positive and whose discrete Fourier coefficients
$\widetilde{\kappa}_{\tau}=\mathcal{F}_H\widetilde{G}_{\tau}$ satisfy
\begin{equation}
\left|
\widetilde{\kappa}_{\tau,\mathbf q}
-
\hat{\kappa}_{\tau,\mathbf q}
\right|
\leq
4\exp\left(
-\frac{\pi^2\tau H^2}{\epsilon L^2}
\right),
\qquad \mathbf q\in U_H,
\end{equation}
provided that $\tau H^2L^{-2}$ is sufficiently large, i.e., $\tau\gtrsim (L/H)^{2-\sigma}$ for some $\sigma>0$.
\end{theorem}

\begin{proof}
For each $\mathbf q\in U_H$, define the aliased Fourier coefficients by
\begin{equation}
\widetilde{\kappa}_{\tau,\mathbf q}
=
\sum_{\mathbf n\in\mathbb Z^3}
\hat{\kappa}_{\tau,\mathbf q+\mathbf nH}
=
\sum_{\mathbf n\in\mathbb Z^3}
\exp\left(
-\frac{\tau}{\epsilon}
\left(
\left|
\frac{2\pi}{L}
(\mathbf q+\mathbf nH)
\right|^2
+k^2
\right)
\right),
\end{equation}
and let
$\widetilde{G}_{\tau}
=
\mathcal{F}_H^{-1}\widetilde{\kappa}_{\tau}$.
By the Poisson summation formula, the inverse discrete Fourier transform of
$\widetilde{\kappa}_{\tau}$ is obtained by sampling the periodization of the
corresponding Gaussian heat kernel. More precisely, for every grid point
$\mathbf x_{\mathbf p}=\mathbf pL/H$,
\begin{equation}
\widetilde{G}_{\tau}
\left(
\mathbf p\frac{L}{H}
\right)
=
\sum_{\mathbf n\in\mathbb Z^3}
G_{\tau}
\left(
\mathbf p\frac{L}{H}+\mathbf nL
\right),
\end{equation}
up to the normalization convention of $\mathcal{F}_H$, where
\begin{equation}
G_{\tau}(\mathbf x)
=
\exp\left(-\frac{k^2\tau}{\epsilon}\right)
\left(\frac{\epsilon}{4\pi\tau}\right)^{\frac{3}{2}}
\exp\left(
-\frac{\epsilon|\mathbf x|^2}{4\tau}
\right).
\end{equation}
Since every term in this sum is strictly positive, all components of
$\widetilde{G}_{\tau}$ are strictly positive.

It remains to estimate the difference between the aliased coefficients and
their principal components. For every $\mathbf q\in U_H$,
\begin{equation}
\begin{aligned}
\left|
\widetilde{\kappa}_{\tau,\mathbf q}
-
\hat{\kappa}_{\tau,\mathbf q}
\right|
&=
\sum_{\mathbf n\in\mathbb Z^3\setminus\{\mathbf0\}}
\exp\left(
-\frac{\tau}{\epsilon}
\left(
\left|
\frac{2\pi}{L}
(\mathbf q+\mathbf nH)
\right|^2
+k^2
\right)
\right)
\\
&\leq
\sum_{\mathbf n\in\mathbb Z^3\setminus\{\mathbf0\}}
\exp\left(
-\frac{\pi^2\tau H^2}{\epsilon L^2}
|\mathbf n|^2
\right).
\end{aligned}
\end{equation}
Writing $a=\pi^2\tau H^2/(\epsilon L^2)$, the remaining Gaussian lattice
tail satisfies
\begin{equation}
\sum_{\mathbf n\in\mathbb Z^3\setminus\{\mathbf0\}}
e^{-a|\mathbf n|^2}
\leq 4e^{-a}
\end{equation}
whenever $a$ is sufficiently large. This proves the stated estimate.

Finally, if
$\tau\gtrsim(L/H)^{2-\sigma}$, then
$\tau H^2/L^2\gtrsim H^\sigma$, and hence
\begin{equation}
\left|
\widetilde{\kappa}_{\tau,\mathbf q}
-
\hat{\kappa}_{\tau,\mathbf q}
\right|
\leq
4e^{-CH^\sigma}
\end{equation}
for some constant $C>0$. Since $e^{-CH^\sigma} \lesssim H^{-r}$ for every
$r>0$, this perturbation is asymptotically smaller than any polynomial
discretization error. Consequently, from an equivalent perspective, the explicitly computable kernel $\widehat{G}_{\tau}$ combined with the positivity-preserving interpolation operator $I$ introduces only an exponentially small violation of positivity in the computed concentration $\hat{c}$.
\end{proof}
Since the mass of each particle remains unchanged throughout the simulation and no particles are absorbed under periodic boundary conditions, mass conservation follows immediately:
\begin{equation}
\int_{\Omega}\hat{\rho}^{(n)}(\mathbf{x})\,\mathrm{d}\mathbf{x}
=
\int_{\Omega}I\hat{\rho}^{(n)}(\mathbf{x})\,\mathrm{d}\mathbf{x}
=
M_0.
\end{equation}
Consequently, the total concentration is given by
\begin{equation}
\int_{\Omega}\hat{c}^{(n)}(\mathbf{x})\,\mathrm{d}\mathbf{x}
=
e^{-k^2n\tau}
\int_{\Omega}\hat{c}^{(0)}(\mathbf{x})\,\mathrm{d}\mathbf{x}
+
\tau M_0
\sum_{s=1}^{n}
e^{-k^2s\tau}.
\end{equation}

For more general computational domains and boundary conditions, the underlying continuous model may itself no longer satisfy mass conservation. Mass may also change because of growth or death processes in cell populations, which appear as reaction terms in the continuous model. To remain consistent with such models, the particles may undergo boundary treatments that alter their masses, while growth and death processes can be incorporated by updating the particle masses according to the corresponding reaction terms. Nevertheless, the proposed algorithm can still be interpreted from the operator-splitting perspective:
\begin{equation}
\hat{c}^{(n+1)}
=
\mathcal{D}_{\tau}
\!\left(
\mathcal{R}_{\tau}
\!\left(
\hat{c}^{(n)},
I\hat{\rho}^{(n)}
\right)
\right),
\end{equation}
where
$\mathcal{R}_{\tau}(f_1,f_2)=f_1+\tau f_2$
is the reaction operator, and
$\mathcal{D}_{\tau}$ denotes the attenuated diffusion operator.

Under this formulation, the computational complexity of each time step is
$O(P+D(H))$,
where $D(H)$ denotes the cost of performing one diffusion step on a grid of resolution $H$. In particular, under periodic boundary conditions,
$D(H) \sim H^3\log H$.

Furthermore, for general computational domains and boundary conditions, if the attenuated diffusion operator $\mathcal{D}_{\tau}$ is positivity-preserving, then the resulting operator-splitting algorithm remains positivity-preserving. The convergence of the SIPF-PIC method with respect to the particle number $P$ and the grid resolution $H$ under general boundary conditions, however, is beyond the scope of the present work.
\end{color}
\section{Convergence of the SIPF-PIC method}
\label{sec:Error}
In this section, we present the convergence result of the proposed SIPF-PIC method (Algorithm \ref{SIPF-PIC}) under three-dimensional periodic boundary conditions, which is stated as Theorem \ref{thm:overall_error}. Specifically, we introduce the following Brownian coupling technique to construct an auxiliary random variable $\widetilde{\mathbf{X}}$ that bridges the ground truth random variable $\mathbf{X}$ and the empirical random variable $\widehat{\mathbf{X}}$ maintained in the algorithm.

\begin{definition}
\label{def:br_cp}
The auxiliary particles \( \widetilde{\mathbf{X}}_p^{(n)} \) (for \( 1 \leq p \leq P \)) satisfy the following update process, \begin{color}{black}where \(\mathbf{W}_p^{(n)}=\mathbf{W}_p^{-(n)}+\mathbf{W}_p^{+(n)}\) are independent increments,\end{color} the same as those driving the corresponding particles mentioned in Algorithm~\ref{SIPF-PIC}:
\begin{equation}
\begin{aligned}
\widetilde{\mathbf{X}}_p^{(0)} &= \mathbf{X}_p^{(0)}, \\
\widetilde{\mathbf{X}}_p^{(n)} &= \widetilde{\mathbf{X}}_p^{(n-1)} + \chi \tau \nabla c(\widetilde{\mathbf{X}}_p^{(n-1)}+\mathbf{W}_p^{-(n)},(n-1)\tau) + \mathbf{W}_p^{(n)}.
\end{aligned}
\end{equation}
\( \widetilde{\mathbf{X}}^{(n)} \) denotes the discrete random variable consisting of particles \( \widetilde{\mathbf{X}}_1^{(n)}, \ldots, \widetilde{\mathbf{X}}_P^{(n)} \), where each particle is sampled with  probability \( \frac{1}{P} \).

The auxiliary discrete density \( \tilde{\rho}^{(n)} \) and its Fourier coefficients $\tilde{\beta}^{(n)}$ are defined as 
\begin{equation}
\tilde{\rho}^{(n)}(\mathbf{x}) = \frac{M_0}{P} \sum\limits_{p=1}^{P} \delta(\mathbf{x} - \widetilde{\mathbf{X}}_p^{(n)}),
\end{equation}
\begin{equation}
\tilde{\beta}_\mathbf{q}^{(n)} = \frac{M_0}{P} L^{-3} \sum_{p=1}^P e^{-i \widetilde{\mathbf{X}}_p^{(n)} \cdot \frac{2\pi}{L} \mathbf{q}}.
\end{equation}
Furthermore, the auxiliary concentration $\tilde{c}^{(n)}$ and its Fourier coefficients $\tilde{\alpha}^{(n)}$ are updated via Algorithm \ref{SIPF-PIC}, \textcolor{black}{where the spectral filter $\Phi_{\mathbf q}=\exp\left(-\frac{2\pi ^2|\mathbf q|^2}{H_0^2}\right)$:}
\begin{equation}
\label{eq:alpha_tilde_update}
\tilde{\alpha}^{(n)}_{\mathbf{q}} = \kappa_{\tau,\mathbf q} \tilde{\alpha}^{(n-1)}_{\mathbf{q}} + \tau \kappa_{\tau,\mathbf q} \tilde{\beta}^{(n-1)}_{\mathbf{q}},\,
\tilde{c}^{(n)} (\mathbf{x})  = \sum_{\mathbf{q} \in U}  \tilde{\alpha}_{\mathbf{q}}^{(n)} \Phi_{\mathbf q} e^{ \frac{2\pi i}{L} \mathbf{q} \cdot \mathbf{x}}.
\end{equation}
The semi-discrete implicit reference concentration $c^{(n)}$ and its Fourier coefficients $\alpha^{(n)}$ are updated as follows:
\begin{equation}
{\alpha}^{(n)}_{\mathbf{q}} = \kappa_{\tau,\mathbf q} {\alpha}^{(n-1)}_{\mathbf{q}} + \tau\kappa_{\tau,\mathbf q}{\beta}^{(n-1)}_{\mathbf{q}}, \text{where}\,\,
{\beta}_{\mathbf{q}}^{(n-1)} = \frac{1}{L^3} \int_{\Omega} {\rho}^{(n-1)}(\mathbf{x}) e^{-i \frac{2\pi}{L} \mathbf{q} \cdot \mathbf{x}} \, \mathrm{d}\mathbf{x},
\end{equation}
while the corresponding random variable \(\mathbf{X}^{(n)}\) with distribution \(\rho ^{(n)}\) is updated via an Euler-Maruyama semi-discretization over the time interval:
\begin{equation}
\mathbf{X}^{(n)} = \mathbf{X}^{(n-1)} + \chi \tau \nabla {c}(\mathbf{X}^{(n-1)}+ \mathbf{W}^{-(n)},(n-1)\tau) + \mathbf{W}^{(n)}.
\label{eq:rho_update}
\end{equation}
\end{definition}

In other words, \( \widehat{\mathbf{X}}_p \) and \( \widetilde{\mathbf{X}}_p \) have the same initial values and diffuse using the same Brownian process, but \( \widehat{\mathbf{X}}_p \) uses \( I^*(\nabla \hat{c}) \) maintained by the algorithm to calculate the drift, while \( \widetilde{\mathbf{X}}_p \) uses the ground truth \( \nabla c \). Since the reference gradient \( \nabla c \) is underlying and does not depend on randomness \(\widetilde{\mathbf{X}}_p^{(0)}\) or \(\mathbf{W}_p^{(n)}\), the coordinates \( \widetilde{\mathbf{X}}_p^{(n)} \) are i.i.d. with respect to \( \mathbf{X}^{(n)} \). Briefly, \( \widetilde{\mathbf{X}}_p \) are independent branches of \(\mathbf{X}\). Also note that $\tilde{c}$ is used only as an intermediate quantity for estimation, and $\nabla \tilde{c}$ does not participate in the update of $\widetilde{\mathbf{X}}_p$.

\begin{remark}
\label{def_random_space}
The probability space in the algorithm is generated as follows: the initial values \((\mathbf{X}_p^{(0)})_{1 \leq p \leq P} \), and the diffusion terms \( \mathbf{W}_p^{(n)}=\mathbf{W}_p^{-(n)}+\mathbf{W}_p^{+(n)} \) in the \( n \)-th step for \( 1 \leq n \leq N_T \). When discussing probabilities or using the expectation notation $\mathbb{E}$, we always operate within this space. The updates of particles $\widehat{\mathbf{X}}_p^{(n)}$ and auxiliary particles $\widetilde{\mathbf{X}}_p^{(n)}$ depend on this shared probability space, which constitutes the essence of the Brownian coupling technique.
\end{remark}
\begin{color}{black}
\begin{definition}
\label{def:mean_abs_error}
The element-wise deviation \(\mathbf{e}_p^{(n)}\) is defined as $\mathbf{e}_p^{(n)} := \widehat{\mathbf{X}}_p^{(n)} - \widetilde{\mathbf{X}}_p^{(n)}$. The mean square error \( d^{(n)} \) is defined as 
\begin{equation}
\label{eq:avg_distance}
d^{(n)} = \frac{1}{P} \sum\limits_{p=1}^{P} |\mathbf{e}_p^{(n)}|^2.
\end{equation}
\end{definition}

\begin{remark}
The quantity \( d^{(n)} \) serves as an upper bound for the \( \mathcal{W}_2 \) distance between the empirical distribution \( \widehat{\mathbf{X}}^{(n)} \) and the auxiliary distribution \( \widetilde{\mathbf{X}}^{(n)} \).
\end{remark}
\end{color}

To ensure the convergence of Algorithm \ref{SIPF-PIC}, we impose the following assumption, which is intended to hold only in the non-blow-up regime. The applicability of the proposed algorithm in the blow-up regime is discussed separately in Remark~\ref{remark:blowup}.
\begin{assumption}
\label{as}
The solution \(\rho, c\) to \eqref{eq:para_system} is periodic over the box\\ \(\Omega=\left[-\frac{L}{2},\frac{L}{2}\right]^3\), that is, for \(q_1,q_2,q_3 \in \mathbf{Z}\),
\begin{equation}
\begin{aligned}
\rho(x_1+q_1 L,x_2+q_2 L,x_3 + q_3 L) &= \rho(x_1,x_2,x_3), \\
c(x_1+q_1 L,x_2+q_2 L,x_3 + q_3 L)    &= c(x_1,x_2,x_3).
\end{aligned}
\label{eq:periodic}
\end{equation}
The function \(\rho \in H^2(\Omega) \cap L^\infty(\Omega)\) , i.e., 
\begin{equation}
\begin{aligned}
\int_{\mathbf{x} \in \Omega}&\left(|\rho(\mathbf{x},t)|^2+ |\nabla_\mathbf{x}\rho(\mathbf{x},t)|^2+|D_\mathbf{x}^2\rho(\mathbf{x},t)|^2\right)\, \mathrm{d}\mathbf{x} \leq N_2^2, \\ 
&|\rho(\mathbf{x},t)|\leq N_0, \quad \forall t \in \left[0,T\right],
\end{aligned}
\end{equation}

The first to third order spatial derivatives of \( c \) and \( \hat{c} \) are bounded, i.e., for any spatial multi-index \( \mathbf{k} \) such that \( |\mathbf{k}| = \eta \) where \(\eta=1,2,3\), it holds that
\begin{equation}
\label{eq:c_smoothness}
\max(|D_{\mathbf{k}} c(\mathbf{x},t)|, |D_{\mathbf{k}} \hat{c}(\mathbf{x},t)|) \leq M_\eta \quad \forall \mathbf{x} \in \left[-\frac{L}{2}, \frac{L}{2}\right]^3,t \in \left[0,T\right].
\end{equation}
The temporal derivative of \(\nabla c\) is bounded, that is,
\begin{equation}
\left| \frac{\partial \nabla c(\mathbf{x},t)}{\partial t} \right| \leq J_1\quad \forall \mathbf{x} \in \left[-\frac{L}{2}, \frac{L}{2}\right]^3,t \in \left[0,T\right].
\end{equation}
\end{assumption}

\subsection{Main convergence results}\label{subsec:main_convergence_proof}
We outline the structure of the argument. We first introduce several error measures required for the subsequent analysis. We then present a series of lemmas that establish bounds for the individual components of the error dynamics, before proving Theorem \ref{thm:overall_error}, which constitutes the main convergence result for Algorithm \ref{SIPF-PIC}. Finally, in Section \ref{subsec:c_hat_relax}, we develop an argument showing that the smoothness assumption on $\hat{c}$ imposed in Assumption \ref{as} can be naturally guaranteed under certain conditions.

We begin by introducing the error measures that will be used throughout the subsequent analysis.

\begin{definition}
The frequency error $\epsilon_\mathbf{q}^{(n)}$ is defined as
\begin{equation}
\epsilon_\mathbf{q}^{(n)}=\tilde{\alpha}_\mathbf{q}^{(n)}-\hat{\alpha}_\mathbf{q}^{(n)}.
\end{equation}
\begin{color}{black}The global velocity error $g^{(n)}$ is defined as
\begin{equation}
g^{(n)} = \left\|
\nabla c^{(n-1)}
-
\nabla\hat{c}^{(n-1)}
\right\|_2^2 = \int_\Omega \left\|
\nabla c^{(n-1)} (\mathbf{x})
-
\nabla\hat{c}^{(n-1)}(\mathbf{x})
\right\|_2^2\,\mathrm{d}\mathbf{x},
\end{equation}
and is decomposed as
\begin{equation}
g_0^{(n)} = \left\|
\nabla c^{(n-1)}
-
\nabla\tilde{c}^{(n-1)}
\right\|_2^2,g_1^{(n)} = \left\|
\nabla \tilde{c}^{(n-1)}
-
\nabla\hat{c}^{(n-1)}
\right\|_2^2.
\end{equation}
The mollified physical displacement error $h^{(n)}$ is defined as
\begin{equation}
\label{eq:hs}
h^{(n)} = \left\|G_0 \ast \left(\frac{1}{P}\sum_{p=1}^P\mathbf{e}_p^{(n)}\delta(\mathbf{x} - \widetilde{\mathbf{X}}_p^{(n)})\right)\right\|_2,
\end{equation}
where $G_0$ is the density of the Gaussian kernel $\mathcal N\left(0,\frac{L^2}{H_0^2}\mathbf I_3\right)$. Here, the \(2\)-norm of a vector-valued $(\mathbb{R}^3\to\mathbb{R}^3)$ function is understood as the sum of the \(2\)-norms of its three components.
\end{color}
\end{definition}

The single-step growth of $d^{(n)}$ is estimated as follows.
\begin{lemma}
\label{lem:Gronwall2}
The mean square error \(d^{(n)}\) is controlled by \(d^{(n-1)}\)
\begin{color}{black}and the global velocity error \(g^{(n-1)}\)\end{color}:
\begin{equation}
\label{eq:e03}
\begin{aligned}
\mathbb{E}[d^{(n)}]
&\leq
\left(1+\tau+3\chi^2(\tau+\tau^2)M_2^2\right)
\mathbb{E}[d^{(n-1)}] \\
&\quad
+3\chi^2(\tau+\tau^2)
\left(
M_3^2\left(\frac{L}{H}\right)^4
+\frac{N_0}{M_0}
\mathbb{E}\left[g^{(n-1)}\right]
\right),
\end{aligned}
\end{equation}
where \(M_2\), \(M_3\), and \(N_0\) are defined in
Assumption \ref{as}, and
\(\nabla c^{(n-1)}=\nabla c(\cdot,(n-1)\tau)\).
\end{lemma}

\begin{proof}
Define the discrete transport mappings
\[
T^{(n-1)}(\mathbf{X})
=
\mathbf{X}+\chi\tau\nabla c^{(n-1)}(\mathbf{X}),
\]
\[
\overline{T}^{(n-1)}(\mathbf{X})
=
\mathbf{X}+\chi\tau\nabla\hat{c}^{(n-1)}(\mathbf{X}),
\]
and
\[
\widehat{T}^{(n-1)}(\mathbf{X})
=
\mathbf{X}+\chi\tau I^*\nabla\hat{c}^{(n-1)}(\mathbf{X}),
\]
respectively. Since the updates of the solutions
\(\mathbf{X}^{(n)}\) and \(\widehat{\mathbf{X}}^{(n)}\)
use the same diffusion \(\mathbf{W}^{-(n)}\) and \(\mathbf{W}^{+(n)}\), we have
\begin{equation}
\label{eq:one_step_particles}
\begin{aligned}
d^{(n)}
&=
\frac{1}{P}\sum_{p=1}^{P}
\left|
T^{(n-1)}(\widetilde{\mathbf{X}}_p^{(n-1)}+\mathbf{W}_p^{-(n)})
-
\widehat{T}^{(n-1)}(\widehat{\mathbf{X}}_p^{(n-1)}+\mathbf{W}_p^{-(n)})
\right|^2.
\end{aligned}
\end{equation}
For notational simplicity, set
\begin{equation}
\label{eq:v_p}
\begin{aligned}
\mathbf{v}_p^{(n-1)}
&=
\chi\Bigl(
\nabla c^{(n-1)}(\widetilde{\mathbf{X}}_p^{(n-1)}+\mathbf{W}_p^{-(n)})
-
I^*\nabla\hat{c}^{(n-1)}
(\widehat{\mathbf{X}}_p^{(n-1)}+\mathbf{W}_p^{-(n)})
\Bigr).
\end{aligned}
\end{equation}
Then the one-step particle error can be written as
\[
\mathbf{e}_p^{(n)}
=
\mathbf{e}_p^{(n-1)}
+
\tau\mathbf{v}_p^{(n-1)}.
\]
Using
\[
|\mathbf{e}+\tau\mathbf{v}|^2
\leq
(1+\tau)|\mathbf{e}|^2
+
(\tau+\tau^2)|\mathbf{v}|^2,
\]
we obtain
\begin{equation}
\label{eq:initial_square_decomposition}
\begin{aligned}
d^{(n)}
&\leq
(1+\tau)d^{(n-1)}
+
\frac{\tau+\tau^2}{P}
\sum_{p=1}^{P}
|\mathbf{v}_p^{(n-1)}|^2.
\end{aligned}
\end{equation}

To estimate the last term, we decompose the velocity error as
\begin{equation}
\begin{aligned}
\frac{1}{\chi}\mathbf{v}_p^{(n-1)}
&=
\Bigl(
\nabla c^{(n-1)}(\widetilde{\mathbf{X}}_p^{(n-1)}+\mathbf{W}_p^{-(n)})
-
\nabla\hat{c}^{(n-1)}
(\widetilde{\mathbf{X}}_p^{(n-1)}+\mathbf{W}_p^{-(n)})
\Bigr)\\
&\quad+
\Bigl(
\nabla\hat{c}^{(n-1)}
(\widetilde{\mathbf{X}}_p^{(n-1)}+\mathbf{W}_p^{-(n)})
-
\nabla\hat{c}^{(n-1)}
(\widehat{\mathbf{X}}_p^{(n-1)}+\mathbf{W}_p^{-(n)})
\Bigr)\\
&\quad+
\Bigl(
\nabla\hat{c}^{(n-1)}
(\widehat{\mathbf{X}}_p^{(n-1)}+\mathbf{W}_p^{-(n)})
-
I^*\nabla\hat{c}^{(n-1)}
(\widehat{\mathbf{X}}_p^{(n-1)}+\mathbf{W}_p^{-(n)})
\Bigr).
\end{aligned}
\end{equation}
It follows from
\[
|\mathbf{a}+\mathbf{b}+\mathbf{c}|^2
\leq
3\left(|\mathbf{a}|^2+|\mathbf{b}|^2+|\mathbf{c}|^2\right)
\]
that
\begin{equation}
\label{eq:velocity_error_decomposition}
\begin{aligned}
\frac{1}{P}\sum_{p=1}^{P}
|\mathbf{v}_p^{(n-1)}|^2
&\leq 3\chi^2
\left(
A_1^{(n-1)}+A_2^{(n-1)}+A_3^{(n-1)}
\right),
\end{aligned}
\end{equation}
where
\begin{equation}
\begin{aligned}
A_1^{(n-1)}
&=
\frac{1}{P}\sum_{p=1}^{P}
\left|
\nabla c^{(n-1)}(\widetilde{\mathbf{X}}_p^{(n-1)}+\mathbf{W}_p^{-(n)})
-
\nabla\hat{c}^{(n-1)}
(\widetilde{\mathbf{X}}_p^{(n-1)}+\mathbf{W}_p^{-(n)})
\right|^2,\\
A_2^{(n-1)}
&=
\frac{1}{P}\sum_{p=1}^{P}
\left|
\nabla\hat{c}^{(n-1)}
(\widetilde{\mathbf{X}}_p^{(n-1)}+\mathbf{W}_p^{-(n)})
-
\nabla\hat{c}^{(n-1)}
(\widehat{\mathbf{X}}_p^{(n-1)}+\mathbf{W}_p^{-(n)})
\right|^2,\\
A_3^{(n-1)}
&=
\frac{1}{P}\sum_{p=1}^{P}
\left|
\nabla\hat{c}^{(n-1)}
(\widehat{\mathbf{X}}_p^{(n-1)}+\mathbf{W}_p^{-(n)})
-
I^*\nabla\hat{c}^{(n-1)}
(\widehat{\mathbf{X}}_p^{(n-1)}+\mathbf{W}_p^{-(n)})
\right|^2.
\end{aligned}
\end{equation}

For the estimate of \(A_1^{(n-1)}\), note that
\(\widetilde{\mathbf{X}}_p^{(n-1)}+\mathbf{W}_p^{-(n)}\) are identically
distributed with density \(\rho^{(n-1/2)}/M_0:=(\rho^{(n-1)} \ast G_{\tau/2})/M_0\). Thus,
using the \(L^\infty\) bound
\(\|\rho^{(n-1/2)}\|_\infty\leq\|\rho^{(n-1)}\|_\infty\leq N_0\), we have
\begin{equation}
\label{eq:e13}
\begin{aligned}
\mathbb{E}[A_1^{(n-1)}]
&=
\frac{1}{M_0}
\mathbb{E}\left[
\int_\Omega
\left|
\nabla c^{(n-1)}(\mathbf{x})
-
\nabla\hat{c}^{(n-1)}(\mathbf{x})
\right|^2
\rho^{(n-1/2)}(\mathbf{x})\,\mathrm{d}\mathbf{x}
\right]\\
&\leq
\frac{N_0}{M_0}
\mathbb{E}\left[
\left\|
\nabla c^{(n-1)}
-
\nabla\hat{c}^{(n-1)}
\right\|_2^2
\right]=
\frac{N_0}{M_0}
\mathbb{E}\left[g^{(n-1)}\right].
\end{aligned}
\end{equation}

For the estimate of \(A_2^{(n-1)}\), we employ
Assumption \ref{as}, which states that the second derivatives
of \(\hat{c}\) are bounded by \(M_2\). By the mean value theorem,
\begin{equation}
\label{eq:e23}
\begin{aligned}
A_2^{(n-1)}
&\leq
\frac{1}{P}\sum_{p=1}^{P}
M_2^2
\left|
\widetilde{\mathbf{X}}_p^{(n-1)}
-
\widehat{\mathbf{X}}_p^{(n-1)}
\right|^2\\
&=
M_2^2d^{(n-1)}.
\end{aligned}
\end{equation}

Finally, by Lemma \ref{Thm:G2P}, the interpolation error
\(A_3^{(n-1)}\) admits the following bound:
\begin{equation}
\label{eq:e33}
\begin{aligned}
A_3^{(n-1)}
&=
\frac{1}{P}\sum_{p=1}^{P}
\left|
\nabla\hat{c}^{(n-1)}
(\widehat{\mathbf{X}}_p^{(n-1)}+\mathbf{W}_p^{-(n)})
-
I^*\nabla\hat{c}^{(n-1)}
(\widehat{\mathbf{X}}_p^{(n-1)}+\mathbf{W}_p^{-(n)})
\right|^2\\
&\leq
\frac{1}{P}\sum_{p=1}^{P}
M_3^2\left(\frac{L}{H}\right)^4=
M_3^2\left(\frac{L}{H}\right)^4.
\end{aligned}
\end{equation}

Combining \eqref{eq:initial_square_decomposition} and
\eqref{eq:velocity_error_decomposition} with
\eqref{eq:e13}, \eqref{eq:e23}, and \eqref{eq:e33}, and
then taking expectations, we obtain
\[
\begin{aligned}
\mathbb{E}[d^{(n)}]
&\leq
(1+\tau)\mathbb{E}[d^{(n-1)}]\\
&\quad+
3\chi^2(\tau+\tau^2)
\left(
\frac{N_0}{M_0}
\mathbb{E}\left[g^{(n-1)}\right]
+
M_2^2\mathbb{E}[d^{(n-1)}]
+
M_3^2\left(\frac{L}{H}\right)^4
\right),
\end{aligned}
\]
which is precisely the target estimate \eqref{eq:e03}.
\end{proof}

The discrepancy between the auxiliary and empirical gradient fields is estimated as follows.
\begin{lemma}
\label{lem:grad_c}
The $2$-norm of the difference between the auxiliary gradient \( \nabla \tilde{c}\) and the empirical gradient \(\nabla \hat{c}\) over the domain $\Omega$ is bounded by
\begin{equation}
\label{eq:gradc_bound}
g_1^{(n)} 
\leq \frac{6M_0^2 \tau}{\epsilon} \sum_{s=0}^{n-1} h^{(s)} + 2 \cdot (2\pi)^6 M_0^2 C^2 L^{-1} H_0^9 H^{-8},
\end{equation}
where the mollified physical displacement error $h^{(s)}$ is defined in \eqref{eq:hs}, and the interpolation error constant is  \(C = \frac{513}{16}\) for the fourth-order particle-to-grid interpolation.
\end{lemma}

\begin{proof}
We use Parseval's theorem to establish the required bound:
\begin{equation}
    g_1^{(n)}=\int_\Omega \left| \nabla \tilde{c}^{(n)}(\mathbf{x}) - \nabla \hat{c}^{(n)}(\mathbf{x}) \right|^2 \mathrm{d}\mathbf{x} = L^3 \sum_{\mathbf{q} \in U} \left| \mathbf{q} \frac{2\pi}{L} \Phi_{\mathbf q} \cdot \epsilon^{(n)}_{\mathbf{q}} \right|^2,
  \label{eq:energy_low_freq}
\end{equation}
where the frequency index set and the spectral filter defined as \(U := \left\{-\frac{H}{2}, \ldots, \frac{H}{2}\right\}^3,\\\Phi_{\mathbf q}=\exp\left(-\frac{2\pi ^2|\mathbf q|^2}{H_0^2}\right)\).

The dynamics of spectral error $\epsilon_{\mathbf{q}}^{(n)}$ defined in Definition \ref{def:mean_abs_error} follow directly from the update formulas \eqref{eq:alpha_hat_update} and \eqref{eq:alpha_tilde_update}. For any frequency index $\mathbf{q} \in U$, 
\begin{equation}
\label{eq:e_decomp}
\begin{aligned}
\epsilon^{(n)}_{\mathbf{q}} 
&= \kappa_{\tau,\mathbf q} \epsilon^{(n-1)}_{\mathbf{q}} + \tau \kappa_{\tau,\mathbf q} 
   \biggl( (\tilde{\beta}^{(n-1)}_{\mathbf{q}} - \hat{\beta}^{(n-1)}_{\mathbf{q}})  + ( \hat{\beta}^{(n-1)}_{\mathbf{q}} - \overline{\beta}^{(n-1)}_{\mathbf{q}} )\biggr),
\end{aligned}
\end{equation}
where $\tilde{\beta}_{\mathbf{q}} := \mathcal{F}_H\tilde{\rho}(\mathbf{q}\frac{2\pi}{L})$,  $\hat{\beta}_{\mathbf{q}} := \mathcal{F}_H\hat{\rho}(\mathbf{q}\frac{2\pi}{L})$, and $\overline{\beta}_{\mathbf{q}} := \mathcal{F}_H I\left( \hat{\rho}\right)(\mathbf{q}\frac{2\pi}{L}) $ denote the discrete Fourier coefficients of the respective density fields.

By Lemma \ref{Thm:P2G}, we have the following bound for the first term:
\begin{equation}
|\hat{\beta}^{(n)}_{\mathbf{q}} - \overline{\beta}^{(n)}_{\mathbf{q}}| \leq C M_0 \left|\mathbf{q} \frac{2\pi}{L}\right|^4 \left( \frac{L}{H} \right)^4 L^{-3}.
\end{equation}

Recall the displacement definition \( \mathbf{e}_p^{(n)} := \widehat{\mathbf{X}}_p^{(n)} - \widetilde{\mathbf{X}}_p^{(n)} \). We have the following bound for the second term: 
\begin{equation}
\bigl| \tilde{\beta}_{\mathbf{q}}^{(n)} - \hat{\beta}_{\mathbf{q}}^{(n)} \bigr| 
\leq \frac{M_0}{P L^3}\left| \sum_{p=1}^P e^{-i\widetilde{\mathbf{X}}_p^{(n)}\mathbf{q}\frac{2\pi}{L}} \mathbf{e}_p^{(n)}\right| \left|\mathbf{q}\frac{2\pi}{L}\right|.
\end{equation}

We next define the spectral displacement error 
\(\nu_{\mathbf{q}}^{(n)} := \frac{1}{P L^3}\sum_{p=1}^P e^{-i\widetilde{\mathbf{X}}_p^{(n)}\mathbf{q}\frac{2\pi}{L}} \mathbf{e}_p^{(n)}\), and the interpolation errors $\hat{\epsilon}_{\mathbf{q}}^{(n)}$ and displacement-induced errors $\check{\epsilon}_{\mathbf{q}}^{(n)}$ as
\begin{equation}
\begin{aligned}
&\hat{\epsilon}_{\mathbf{q}}^{(0)}:=0,\hat{\epsilon}_{\mathbf{q}}^{(n)}:=\kappa_{\tau,\mathbf q}\hat{\epsilon}^{(n-1)}_{\mathbf{q}}+ \tau \kappa_{\tau,\mathbf q}
  ( \hat{\beta}^{(n-1)}_{\mathbf{q}} - \overline{\beta}^{(n-1)}_{\mathbf{q}} ),\\
&\check{\epsilon}_{\mathbf{q}}^{(0)}:=0,\check{\epsilon}_{\mathbf{q}}^{(n)}:=\kappa_{\tau,\mathbf q} \check{\epsilon}^{(n-1)}_{\mathbf{q}}+ \tau \kappa_{\tau,\mathbf q}
  ( \tilde{\beta}^{(n-1)}_{\mathbf{q}} - \hat{\beta}^{(n-1)}_{\mathbf{q}} ).
\end{aligned}
\end{equation}
By \eqref{eq:e_decomp} we have $|\epsilon_{\mathbf{q}}^{(n)}| \leq |\hat{\epsilon}_{\mathbf{q}}^{(n)}| + |\check{\epsilon}_{\mathbf{q}}^{(n)}|.$

These two error terms admit the following estimates, respectively
\begin{equation}
\label{eq:e1}
|\hat{\epsilon}_{\mathbf{q}}^{(n)}| \leq \kappa_{\tau,\mathbf q} |\hat{\epsilon}_{\mathbf{q}}^{(n-1)}| 
+ M_0 L^{-3} \tau \kappa_{\tau,\mathbf q} \cdot  C \left|\mathbf{q}\frac{2\pi}{L}\right|^4 \left( \frac{L}{H} \right)^4,
\end{equation}
and
\begin{equation}
\label{eq:e2}
|\check{\epsilon}_{\mathbf{q}}^{(n)}| \leq \kappa_{\tau,\mathbf q}|\check{\epsilon}_{\mathbf{q}}^{(n-1)}| 
+ M_0 \tau \kappa_{\tau,\mathbf q} \cdot  \left|\mathbf{q}\frac{2\pi}{L}\right|\left| \nu_{\mathbf{q}}^{(n-1)} \right|.
\end{equation}

For \eqref{eq:e1}, summing the decaying geometric series yields
\begin{equation}
\label{eq:sum_decay}
\begin{aligned}
|\hat{\epsilon}_{\mathbf{q}}^{(n)}| 
&\leq \frac{1 + \frac{\tau}{\epsilon}\left(|\mathbf{q}\frac{2\pi}{L}|^2 + k^2\right)}{\frac{\tau}{\epsilon}\left(|\mathbf{q}\frac{2\pi}{L}|^2 + k^2\right)} 
\cdot \frac{M_0 L^{-3}}{|\mathbf{q}\frac{2\pi}{L}|^2 + k^2 + \frac{\epsilon}{\tau}} 
\cdot  C \left|\mathbf{q}\frac{2\pi}{L}\right|^4 \left( \frac{L}{H} \right)^4 \\
&= \frac{M_0 L^{-3}}{|\mathbf{q}\frac{2\pi}{L}|^2 + k^2} 
\cdot  C \left|\mathbf{q}\frac{2\pi}{L}\right|^4 \left( \frac{L}{H} \right)^4.
\end{aligned}
\end{equation}

For \eqref{eq:e2}, we derive the estimate as follows:
\begin{equation}
\begin{aligned}
|\check{\epsilon}_{\mathbf{q}}^{(n)}|^2 &\leq \frac{M_0^2 \tau^2}{\epsilon^2} \left( \left|\nu_{\mathbf{q}}^{(n-1)}\right| + \kappa_{\tau,\mathbf q} \left|\nu_{\mathbf{q}}^{(n-2)}\right| + \cdots + \kappa_{\tau,\mathbf q}^{n-1} \left|\nu_{\mathbf{q}}^{(0)}\right| \right)^2\\
&\leq \frac{M_0^2 \tau^2}{\epsilon^2} \left(1 + \kappa_{\tau,\mathbf q} ^2 + \cdots + \kappa_{\tau,\mathbf q}^{2n-2} \right)\left(\left|\nu_{\mathbf{q}}^{(n-1)}\right|^2 + \left|\nu_{\mathbf{q}}^{(n-2)}\right|^2 + \cdots +  \left|\nu_{\mathbf{q}}^{(0)}\right|^2 \right)\\
&\leq \frac{M_0^2 \tau^2}{\epsilon^2} \cdot \frac{1}{1-\kappa_{\tau,\mathbf q} ^2}\left(\left|\nu_{\mathbf{q}}^{(n-1)}\right|^2 + \left|\nu_{\mathbf{q}}^{(n-2)}\right|^2 + \cdots +  \left|\nu_{\mathbf{q}}^{(0)}\right|^2 \right)\\
&\leq \frac{3M_0^2 \tau}{\epsilon} \cdot \frac{1}{|\mathbf{q}\frac{2\pi}{L}|^2 + k^2}\left(\left|\nu_{\mathbf{q}}^{(n-1)}\right|^2 + \left|\nu_{\mathbf{q}}^{(n-2)}\right|^2 + \cdots +  \left|\nu_{\mathbf{q}}^{(0)}\right|^2 \right).
\end{aligned}
\end{equation}

Applying the above spectral error bounds and using Parseval's identity \eqref{eq:energy_low_freq}, we translate the frequency-domain estimates back to the physical space to proceed with the proof. Specifically, we can obtain that 
\begin{equation}
\sum_{\mathbf{q} \in U} \left| \mathbf{q} \frac{2\pi}{L} \Phi_{\mathbf q}\cdot \epsilon_{\mathbf{q}}^{(n)} \right|^2 \leq 2\sum_{\mathbf{q} \in U} \left| \mathbf{q} \frac{2\pi}{L}\Phi_{\mathbf q} \cdot \hat{\epsilon}_{\mathbf{q}}^{(n)} \right|^2 + 2\sum_{\mathbf{q} \in U} \left| \mathbf{q} \frac{2\pi}{L} \Phi_{\mathbf q}\cdot \check{\epsilon}_{\mathbf{q}}^{(n)} \right|^2,
\label{eq:energy_sum_low_freq}
\end{equation}
\begin{equation}
  \begin{aligned}
    \sum_{\mathbf{q} \in U} \left| \mathbf{q} \frac{2\pi}{L} \Phi_{\mathbf q}\cdot \hat{\epsilon}_{\mathbf{q}}^{(n)} \right|^2 
    &\leq M_0^2 L^{-6} \sum_{\mathbf{q} \in U} \left(\frac{C\left|\mathbf{q} \frac{2\pi}{L}\right|^5}{\left|\mathbf{q} \frac{2\pi}{L}\right|^2 + k^2} \left( \frac{L}{H} \right)^4 e^{-\frac{2\pi ^2|\mathbf q|^2}{H_0^2}}\right) ^2 \\
    &\leq H_0^3 \cdot M_0^2 C^2 L^{-6} \left( H_0 \frac{2\pi}{L} \right)^6 \left( \frac{L}{H} \right)^8,
  \end{aligned}
  \label{eq:energy_sum_low_freq1}
\end{equation}
and 
\begin{equation}
  \begin{aligned}
    \sum_{\mathbf{q} \in U} \left| \mathbf{q} \frac{2\pi}{L}\Phi_{\mathbf q} \cdot \check{\epsilon}_{\mathbf{q}}^{(n)} \right|^2 
    &\leq \frac{3M_0^2 \tau}{\epsilon} \sum_{\mathbf{q} \in U} \frac{|\mathbf{q}\frac{2\pi}{L}|^2}{|\mathbf{q}\frac{2\pi}{L}|^2 + k^2} \sum_{s=0}^{n-1} \left|\Phi_{\mathbf q}\nu_{\mathbf{q}}^{(s)}\right|^2 \\
    &\leq \frac{3M_0^2 \tau}{\epsilon} \sum_{s=0}^{n-1} \sum_{\mathbf{q} \in U}  \left|\Phi_{\mathbf q}\nu_{\mathbf{q}}^{(s)}\right|^2\\
    &\leq \frac{3M_0^2 \tau}{\epsilon L^3} \sum_{s=0}^{n-1} h^{(s)}.
  \end{aligned}
  \label{eq:energy_sum_low_freq2}
\end{equation}
Substituting \eqref{eq:energy_sum_low_freq}-\eqref{eq:energy_sum_low_freq2} into \eqref{eq:energy_low_freq} yields the desired result.
\end{proof}

The discrepancy between the reference and auxiliary gradient fields is estimated as follows.
\begin{lemma}
\label{lem:grad_c2}
The $2$-norm of the difference between the reference gradient \( \nabla c\) and the auxiliary gradient \( \nabla \tilde{c}\) is bounded by
\begin{equation}
\label{eq:gradc_bound2}
\mathbb{E} \left[ g_0^{(n)} \right] \leq \frac{M_0^2 H_0}{PL} + L^4 M_3^2 {H}_0^{-4},
\end{equation}
where \(M_3\) is defined in Assumption \ref{as}, and the underlying probability space is specified in Remark \ref{def_random_space}.
\end{lemma}

\begin{proof}
Note $\Phi_{\mathbf q}=0$ for $\mathbf q \in \mathbb{Z}^3 \setminus U$. By Parseval's theorem, we have
\begin{equation}
  \begin{aligned}
    g_0^{(n)}&=\int_\Omega \left| \nabla c^{(n)}(\mathbf{x}) - \nabla \tilde{c}^{(n)}(\mathbf{x}) \right|^2 \mathrm{d}\mathbf{x} \\
    &= L^3 \biggl( \sum_{\mathbf{q} \in \mathbb{Z}^3}\left|(1-\Phi_{\mathbf q})\mathbf{q}\frac{2\pi}{L} \alpha^{(n)}_{\mathbf{q}} \right|^2  + \sum_{\mathbf{q} \in U} \left| \Phi_{\mathbf q}\mathbf{q}\frac{2\pi}{L} (\alpha^{(n)}_{\mathbf{q}} - \tilde{\alpha}^{(n)}_{\mathbf{q}}) \right|^2 \biggr).
  \end{aligned}
  \label{eq:random}
\end{equation}
Let the first part
$S_1^{(n)}=\sum_{\mathbf q\in\mathbb Z^3}
\left|
(1-\Phi_{\mathbf q})
\mathbf q\frac{2\pi}{L}
\alpha_{\mathbf q}^{(n)}
\right|^2.$ For $\mathbf q\in U$, 
$1-\Phi_{\mathbf q}\leq\frac{L^2}{2H_0^2}
\left|\mathbf q\frac{2\pi}{L}\right|^2$, and for $\mathbf q\notin U$, we have $|\mathbf q|\geq H/2$ and hence
$
\left|\mathbf q\frac{2\pi}{L}\right|^2
\leq
\frac{L^4}{\pi^4H^4}
\left|\mathbf q\frac{2\pi}{L}\right|^6.$
Splitting the summation over $U$ and its complement gives
\begin{equation}
S_1^{(n)}
\leq
L^4
\left(
\frac{1}{4H_0^4}
+
\frac{1}{\pi^4H^4}
\right)
\sum_{\mathbf q\in\mathbb Z^3}
\left|\mathbf q\frac{2\pi}{L}\right|^6
|\alpha_{\mathbf q}^{(n)}|^2.
\end{equation}
By Parseval's identity,
\begin{equation}
\begin{aligned}
\sum_{\mathbf q\in\mathbb Z^3}
\left|\mathbf q\frac{2\pi}{L}\right|^6
|\alpha_{\mathbf q}^{(n)}|^2
\leq
M_3^2.
\end{aligned}
\end{equation}
Consequently,
\begin{equation}
S_1^{(n)}
\leq
M_3^2
\left[
\frac{1}{4}\left(\frac{L}{H_0}\right)^4
+
\frac{1}{\pi^4}\left(\frac{L}{H}\right)^4
\right].
\end{equation}
For $H\geq H_0$, it follows that
\begin{equation}
\label{eq:high_freq}
S_1^{(n)}
\leq
M_3^2\left(\frac{L}{H_0}\right)^4.
\end{equation}

For the second part, the independence of the particles \( \widetilde{\mathbf{X}}_1^{(n)}, \ldots, \widetilde{\mathbf{X}}_P^{(n)} \sim \mathbf{X}^{(n)} \) implies that 
\begin{equation}
\begin{aligned}
\mathbb{E}\left| \beta^{(n)}_{\mathbf{q}} - \tilde{\beta}^{(n)}_{\mathbf{q}} \right|^2 
&= \frac{1}{P} \mathrm{Var}\!\left(M_0 L^{-3}e^{-i\mathbf{q}\frac{2\pi}{L} \widetilde{\mathbf{X}}^{(n)}}\right) \\
&= \frac{M_0^2}{P}L^{-6}\left(1 - \left|\beta^{(n)}_{2\mathbf{q}}\right|^2\right) \\
&\leq \frac{M_0^2}{PL^6}.
\end{aligned}
\end{equation}

Similar to \eqref{eq:sum_decay}, applying the attenuation factor \(\frac{1}{|\mathbf{q}\frac{2\pi}{L}|^2 + k^2}\) and taking the expectation over the sum yields
\begin{equation}
\begin{aligned}
\mathbb{E}\left[\sum_{\mathbf{q} \in U} \left| \mathbf{q}\frac{2\pi}{L}\Phi_{\mathbf q}(\alpha^{(n)}_{\mathbf{q}} - \tilde{\alpha}^{(n)}_{\mathbf{q}}) \right|^2 \right]
&\leq \sum_{\mathbf{q} \in U} \frac{|\mathbf{q}\frac{2\pi}{L}\Phi_{\mathbf q}|^2}{(|\mathbf{q}\frac{2\pi}{L}|^2 + k^2)^2} \cdot \frac{M_0^2}{PL^6}\\
&\leq \frac{M_0^2 H_0}{P L^4}.
\end{aligned}
\label{eq:random_low_freq}
\end{equation}

The last inequality follows from a standard lattice-sum estimate. Grouping the three-dimensional lattice points into shells gives
\(\sum_{\mathbf q\neq\mathbf 0}|\mathbf q|^{-2}e^{-4\pi^2|\mathbf q|^2/H_0^2}\leq H_0\).

Combining the high-frequency error bound \eqref{eq:high_freq} and the low-frequency variance estimate \eqref{eq:random_low_freq} via the Parseval decomposition \eqref{eq:random}, we establish the desired final error bound.
\end{proof}

\begin{color}{black}
The single-step growth of the mollified physical displacement error $h^{(n)}$ is estimated as follows.
\begin{lemma}
\label{lemma:h_bound}
Under the conditions \(\mu\tau \geq \frac{L^2}{H_0^2}\) and \(P\gtrsim H_0^5\), the mollified physical displacement error $h^{(n)}$ is controlled by $h^{(n-1)}$, $d^{(n-1)}$ and $g^{(n-1)}$, where $\mu$ is the diffusion coefficient of particles in system \eqref{eq:para_system}:
\begin{equation}
\begin{aligned}
\mathbb{E}[h^{(n)}] &\leq \mathbb{E}[h^{(n-1)}](1+(1+3\chi^2M_2) \tau) +8\chi^2 N_2^2 M_3^2 \left(\frac{L}{H}\right)^4\\
&+ \tau \left(\left(3\mu \frac{H_0^5}{PL^5}+8\chi ^2 N_2^2M_2^2\right)\mathbb{E}[d^{(n-1)}]  +4 \chi^2N_0^2\mathbb{E}[g^{(n-1)}]\right).
\end{aligned}
\end{equation}
\end{lemma}
\begin{proof}
Recall $h^{(n)} = \left\|G_0 \ast \left(\sum_{p=1}^P\mathbf{e}_p^{(n)}\delta(\mathbf{x} - \widetilde{\mathbf{X}}_p^{(n)})\right)\right\|_2$. Define the middle step
\begin{equation}
h^{(n-\frac{1}{2})} = \left\|G_0 \ast \left(\sum_{p=1}^P\mathbf{e}_p^{(n-1)}\delta(\mathbf{x} - \widetilde{\mathbf{X}}_p^{(n)})\right)\right\|_2,
\end{equation}
we have $h^{(n)} - h^{(n-1)} = (h^{(n-\frac{1}{2})}-h^{(n-1)})+(h^{(n)}-h^{(n-\frac{1}{2})})$.

For the first part, recall the discrete transport mapping \(T^{(n-1)}(\mathbf{X})
=\\
\mathbf{X}+\chi\tau\nabla c^{(n-1)}(\mathbf{X})\). Note the distribution-valued error field
\begin{equation}
\begin{aligned}
\mathcal{E}^{(n)}&:=\frac{1}{P}\sum_{p=1}^P\mathbf{e}_p^{(n)}\delta(\mathbf{x} - \widetilde{\mathbf{X}}_p^{(n)}),\mathcal{E}^{(n-\frac{1}{2})}:=\frac{1}{P}\sum_{p=1}^P\mathbf{e}_p^{(n-1)}\delta(\mathbf{x} - \widetilde{\mathbf{X}}_p^{(n)}),\\
\mathcal{E}^{(n-\frac{3}{4})}&:=\frac{1}{P}\sum_{p=1}^P\mathbf{e}_p^{(n-1)}\delta\left(\mathbf{x} - T^{(n-1)}(\widetilde{\mathbf{X}}_p^{(n-1)}+\mathbf{W}_p^{-(n)})\right).
\end{aligned}
\end{equation}
To avoid confusion, throughout the proof of this lemma, we use
\(K_\sigma := \operatorname{Law}\bigl(\mathcal{N}(0,\sigma I_3)\bigr)\)
to denote the three-dimensional Gaussian kernel with covariance matrix
$\sigma I_3$. With this convention, the convolution of two such kernels
satisfies
\(K_\sigma * K_{\sigma'} = K_{\sigma+\sigma'}.\)
The Gaussian mollifier introduced in the algorithm is given by
\(G_0 = K_{{L^2}/{H_0^2}}\). Simple computations show that \(\|K_\sigma\|_2^2 \sim \sigma^{-\frac{3}{2}},\|\nabla K_\sigma\|_2^2 \sim \sigma^{-\frac{5}{2}}\) for \(\sigma \to 0\).

By the update scheme of \(\widetilde{\mathbf{X}}_p^{(n)}\), we have
\begin{equation}
\mathbb{E}\left[G_0 \ast \mathcal{E}^{(n-\frac{1}{2})}\right] = K_{\frac{L^2}{H_0^2}+\mu \tau} \ast T_\#^{(n-1)}(K_{\mu \tau} \ast \mathcal{E}^{(n-1)} ).
\end{equation}
Independence of $\mathbf{W}_p^{+(n)}$ gives the exact identity (notice that $d^{(n-1)}\frac{1}{P}\sum_{p=1}^P |\mathbf{e}_p^{(n-1)}|^2$)
\begin{equation}\label{eq:post-variance-exact}
\begin{aligned}
&\mathbb{E}\Bigl[
  \bigl\|K_{\frac{L^2}{H_0^2}} * \mathcal{E}^{(n-\frac{1}{2})} - K_{\frac{L^2}{H_0^2}+\mu\tau} * \mathcal{E}^{(n-\frac{3}{4})}\bigr\|_2^2
  \,\Big|\,
  \mathbf{W}_1^{-(n)},
  \ldots,
  \mathbf{W}_P^{-(n)}
\Bigr]\\
&=
\frac{d^{(n-1)}}{P}
\left(
\bigl\|K_{\frac{L^2}{H_0^2}}\bigr\|_2^2
-
\bigl\|K_{\frac{L^2}{H_0^2}+\mu \tau}\bigr\|_2^2
\right).
\end{aligned}
\end{equation}
By Parseval's identity, we have the difference of $L^2$ norm
\begin{equation}
\begin{aligned}
\bigl\|K_{\frac{L^2}{H_0^2}}\bigr\|_2^2
-
\bigl\|K_{\frac{L^2}{H_0^2}+\mu \tau}\bigr\|_2^2 &= L^3 \sum_{\mathbf{q} \in \mathbb{Z}^3}e^{-\frac{4\pi^2|\mathbf{q}|^2}{H_0^2}}\left(1-e^{-\frac{4\pi^2|\mathbf{q}|^2 \mu \tau}{L^2}}\right)\\
&\leq L^3 \sum_{\mathbf{q} \in \mathbb{Z}^3} \frac{4\pi^2|\mathbf{q}|^2 \mu \tau}{L^2}e^{-\frac{4\pi^2|\mathbf{q}|^2}{H_0^2}}\\
&=\mu \tau \bigl\|\nabla K_{\frac{L^2}{H_0^2}}\bigr\|_2^2 \leq \mu \tau \frac{H_0^5}{L^5}.
\end{aligned}
\end{equation}
To control the variance generated by the independent increments $\mathbf{W}_p^{-(n)}$, define the
random kernels $\Psi_p^{(n-1)}(\mathbf{x};\mathbf{w}):=K_{\frac{L^2}{H_0^2}+\mu \tau}(\mathbf{x}-T^{(n-1)}(\widetilde{\mathbf{X}}_p^{(n-1)}+\mathbf{w}))$ and the average kernels $\overline{\Psi}_p^{(n-1)}(\mathbf{x};\mathbf{w})=\mathbb{E}\left[\Psi_p^{(n-1)}(\mathbf{x};\mathbf{w})\big|\mathbf{w}\right].$ We have the exact identity
\begin{equation}
\begin{aligned}
&\mathbb{E}\Bigl[\bigl\|
  K_{\frac{L^2}{H_0^2}+\mu\tau} * \mathcal{E}^{(n-\frac{3}{4})}-K_{\frac{L^2}{H_0^2}+\mu \tau} \ast T_\#^{(n-1)}(K_{\mu \tau} \ast \mathcal{E}^{(n-1)})\bigr\|_2^2
\Bigr]\\
&=\frac{1}{P^2}\sum_{p=1}^P |\mathbf{e}_p^{(n-1)}|^2\mathbb{E}\left[\left\|\Psi_p^{(n-1)}(\cdot;\mathbf{W}_p^{-(n)})-\overline{\Psi}_p^{(n-1)}(\cdot;\mathbf{W}_p^{-(n)})\right\|_2^2\right].
\end{aligned}
\end{equation}
Directly apply the Lipschitz bound on the push-forward map, and notice that the second derivative bound $M_2$ in Assumption~\ref{as} is independent of the timestep $\tau$, we have
\begin{equation}
\begin{aligned}
\mathbb{E}\left[\left\|\Psi_p^{(n-1)}(\cdot;\mathbf{W}_p^{-(n)})-\overline{\Psi}_p^{(n-1)}(\cdot;\mathbf{W}_p^{-(n)})\right\|_2^2\right]&\leq \|\nabla T^{(n-1)}\|_\infty^2 \mathbb{E}[\mathbf{W}_p^{-(n)}]\|\nabla K_{\frac{L^2}{H_0^2}+\mu\tau}\|_2^2\\
&\leq  (1+\tau M_2)^2 \mu \tau \left(\frac{L^2}{H_0^2}+\mu\tau\right)^{-\frac{5}{2}}\\
&\leq 2 \mu \tau \frac{H_0^5}{L^5}.
\end{aligned}
\end{equation}
Finally, the compression rate of the push-forward map is bounded by the second derivative of $c$:
\begin{equation}
\left\|T_\#^{(n-1)}(K_{\mu \tau} \ast \mathcal{E}^{(n-1)})\right\|_2^2 \leq (1+3\chi^2\tau M_2) \left\|K_{\mu \tau} \ast \mathcal{E}^{(n-1)}\right\|_2^2.
\end{equation}
Under the condition \(\mu \tau \geq \frac{L^2}{H_0^2}\), combining all the estimates gives
\begin{equation}
\begin{aligned}
&\mathbb{E}\left[\left\|G_0 \ast \mathcal{E}^{(n-\frac{1}{2})}\right\|_2^2\right] \\
&\leq \left\|K_{\frac{L^2}{H_0^2}+\mu \tau} \ast T_\#^{(n-1)}(K_{\mu \tau} \ast \mathcal{E}^{(n-1)} )\right\|_2^2+3d^{(n-1)}\mu\tau \frac{H_0^5}{PL^5}\\
&\leq \left\|T_\#^{(n-1)}(K_{\mu \tau} \ast \mathcal{E}^{(n-1)} )\right\|_2^2+3d^{(n-1)}\mu\tau \frac{H_0^5}{PL^5}\\
&\leq (1+3\tau M_2)\left\|K_{\mu \tau} \ast \mathcal{E}^{(n-1)}\right\|_2^2+3d^{(n-1)}\mu\tau \frac{H_0^5}{PL^5}\\
&\leq (1+3\tau M_2)\left\|G_0 \ast \mathcal{E}^{(n-1)}\right\|_2^2+3d^{(n-1)}\mu\tau \frac{H_0^5}{PL^5}.
\end{aligned}
\end{equation}

For the second part, recall the definition \eqref{eq:v_p} of the single-step velocity error $\mathbf{v}_p^{(n-1)}$. Using \(
|\mathbf{e}+\tau\mathbf{v}|^2
\leq
(1+\tau)|\mathbf{e}|^2
+
(\tau+\tau^2)|\mathbf{v}|^2
\), we have the following estimate:
\begin{equation}
\begin{aligned}
\mathbb{E}\left[\left\|G_0 \ast \mathcal{E}^{(n)}\right\|_2^2\right] \leq& (1+\tau)\mathbb{E}\left[\left\|G_0 \ast \mathcal{E}^{(n-\frac{1}{2})}\right\|_2^2\right]\\
&+ (\tau+\tau^2)\mathbb{E}\left[\left\|\frac{1}{P}\sum_{p=1}^P\mathbf{v}_p^{(n-1)}G_0(\mathbf{x} - \widetilde{\mathbf{X}}_p^{(n)})\right\|_2^2\right].
\end{aligned}
\end{equation}
Lemma~\ref{lem:Gronwall2} gives
\begin{equation}
\left|\mathbf{v}_p^{(n-1)}-\chi \left(\nabla c^{(n-1)}-\nabla \hat{c}^{(n-1)}\right)(\widetilde{\mathbf{X}}_p^{(n-1)})\right| \leq \chi \left(M_2|\mathbf{e}_p^{(n-1)}|+M_3\left(\frac{L}{H}\right)^2\right).
\end{equation}
Under the condition \(P\gtrsim H_0^5\), the nominal part is controlled by
\begin{equation}
\begin{aligned}
&\mathbb{E}\left[\left\|\frac{1}{P}\sum_{p=1}^P\chi \left(\nabla c^{(n-1)}-\nabla \hat{c}^{(n-1)}\right)(\widetilde{\mathbf{X}}_p^{(n-1)})\,G_0(\mathbf{x} - \widetilde{\mathbf{X}}_p^{(n)})\right\|_2^2\right]
\\&\leq\chi^2\left\|\nabla c^{(n-1)}-\nabla \hat{c}^{(n-1)}\right\|_2^2 \left(N_0^2+\frac{H_0^3}{PL^3}\right) \leq 2\chi^2 N_0^2 g^{(n-1)}.
\end{aligned}
\end{equation}

The disturbance part is controlled by
\begin{equation}
\begin{aligned}
&\mathbb{E}\left[\left\|\frac{1}{P}\sum_{p=1}^P\chi \left(M_2|\mathbf{e}_p^{(n-1)}|+M_3\left(\frac{L}{H}\right)^2\right)G_0(\mathbf{x} - \widetilde{\mathbf{X}}_p^{(n)})\right\|_2^2\right]
\\&\leq\chi^2\left(2\frac{M_2^2}{P}\sum_{p=1}^P|\mathbf{e}_p^{(n-1)}|^2+2M_3^2\left(\frac{L}{H}\right)^4\right) \left(N_2^2+\frac{H_0^3}{PL^3}\right)
\\&\leq 4 \chi^2 N_2^2 \left(M_2^2 d^{(n-1)}+M_3^2\left(\frac{L}{H}\right)^4\right).
\end{aligned}
\end{equation}
Combining all the estimates gets the target bound.
\end{proof}
\end{color}

The temporal error of the Strang-type semi-discretization over the time interval is stated as follows (see \textcite[p.~15-16]{milstein2004stochastic}).
\begin{lemma}
\label{lemma:temporal}
Let \(\{\mathbf{X}_p\}_{1 \leq p \leq P}\) be \(P\) independent realizations of the SDE \eqref{eq:SDE}, and let the Brownian increments \(\mathbf{W}_p\) satisfy
\begin{equation}
\mathbf{W}_p(n\tau) - \mathbf{W}_p((n-1)\tau) = \mathbf{W}_p^{(n)}. 
\end{equation}
The Strang-type semi-discretization exhibits first-order strong convergence for additive noise, i.e., when the diffusion term is independent of the solution process. More precisely, for \(t = n\tau \leq T\) with integer \(n\geq 1\),
\begin{equation}
\left(\mathbb{E}\bigl[ | \mathbf{X}_p(t) - \widetilde{\mathbf{X}}_p^{(n)} |^2 \bigr]\right)^\frac{1}{2}\leq \left( C_0 (t^2 +t)\right)^\frac{1}{2} \tau,   
\end{equation}
where \(C_0 = C_0(J_1,M_1,M_2,M_3)\) is a positive constant independent of the time step \(\tau\) and \(t\).
\end{lemma}

The main convergence result is presented as follows.
\begin{theorem}
\label{thm:overall_error}
For the SIPF-PIC method with fourth-order particle-to-grid and second-order grid-to-particle interpolation, the mean absolute error at the final time T, defined as $w^{(n)}:=\left(\frac{1}{P}\sum_{1\leq p \leq P}|\mathbf{X}_p(n \tau)-\widehat{\mathbf{X}}_p^{(n)}|^2\right)^\frac{1}{2}$, satisfies the following growth bound, \begin{color}{black}under  Assumption~\ref{as} and the conditions \(\mu\tau \geq \frac{L^2}{H_0^2}\) and \(P \geq H_0^5(1+L^{-5})\)\end{color}:
\begin{equation}
\label{eq:main_results_thm}
\begin{aligned}
\mathbb{E}\left[w^{\left(N_T\right)}\right] &\leq C_2 e^{C_1T}
\Bigg(
C_3 L^2 H^{-2}
+C_4L^{-\frac{1}{2}}H_0^\frac{9}{2}H^{-4} \\
&\hspace{3.2cm}
+C_5L^2H_0^{-2}
+C_6L^{-\frac{1}{2}}P^{-\frac{1}{2}}H_0^\frac{1}{2}+C_0 \tau
\Bigg),
\end{aligned}
\end{equation}
where \(C_0\) to \(C_6\) are positive constants derived in the proof, independent of the domain size \(L\), grid resolution \(H\), filter resolution \(H_0\), number of particles \(P\), and time step \(\tau\). \begin{color}{black}More specifically, the term $C_3L^2H^{-2}$ arises from the second-order interpolation error of $\nabla\hat{c}$, while $C_4L^{-\frac{1}{2}}H_0^{\frac{9}{2}}H^{-4}$ results from the summation of all frequencies of the fourth-order interpolation error of $e^{\frac{2\pi i}{L}\mathbf{q}\cdot\mathbf{x}}$. The term $C_5L^2H_0^{-2}$ represents the error introduced by the Gaussian filter with cut-off frequency $H_0$, and $C_6L^{-\frac{1}{2}}P^{-\frac{1}{2}}H_0^{\frac{1}{2}}$ is the statistical sampling error associated with the i.i.d.\ auxiliary particles.\end{color}
\end{theorem}

\begin{proof}
Define the total error energy \(D^{(n)}:=\mathbb{E}[d^{(n)}]+\mathbb{E}[h^{(n)}]\). First, we employ the discrete Gronwall inequality to derive an estimate for \(D^{(n)}\). From Lemmas~\ref{lem:Gronwall2} and \ref{lemma:h_bound}, we have
\begin{equation}
\label{eq:D_part}
D^{(n)}
\leq
(1+C_1^*\tau)D^{(n-1)}
+C_2^*\tau g^{(n-1)}
+C_3^*\tau\left(\frac{L}{H}\right)^4,
\end{equation}
where the constants
\(C_1^*=1+3\mu+M_2^2(3+8\chi^2N_2^2)\),
\(C_2^*=3M_0^{-1}N_0\chi^2+4N_0^2\chi^2\), and
\(C_3^*=(3+8N_2^2)\chi^2M_3^2\)
are independent of the computational parameters: domain size \(L\), grid resolution \(H\), filter resolution \(H_0\), number of particles \(P\), and time step \(\tau\).

From Lemmas~\ref{lem:grad_c} and \ref{lem:grad_c2}, we have
\begin{equation}
\label{eq:G_part}
\mathbb{E}[g^{(n)}]
\leq
\frac{6M_0^2}{\epsilon}\tau
\sum_{s=0}^{n-1}D^{(s)}
+C_4^*L^{-1}H_0^9H^{-8}
+C_5^*L^4H_0^{-4}
+C_6^*L^{-1}P^{-1}H_0,
\end{equation}
where
\(C_4^*=\left(\frac{513}{16}\right)^2 2^7\pi^6M_0^2\),
\(C_5^*=M_0^2\), and \(C_6^*=M_3^2\)
are independent of the computational parameters.
Let \(
S^{(n)}:=\tau\sum_{s=0}^{n-1}D^{(s)}.
\) Substituting \eqref{eq:G_part} into \eqref{eq:D_part} gives
\begin{equation}
\label{eq:D_combined}
\begin{aligned}
D^{(n)}
\leq{}&
(1+C_1^*\tau)D^{(n-1)}
+\frac{6C_2^*M_0^2}{\epsilon}\,
\tau^2\sum_{s=0}^{n-2}D^{(s)} \\
&+\tau C_3^*\left(\frac{L}{H}\right)^4
+\tau C_2^*C_4^*L^{-1}H_0^9H^{-8}
+\tau C_2^*C_5^*L^4H_0^{-4} \\
&+\tau C_2^*C_6^*L^{-1}P^{-1}H_0.
\end{aligned}
\end{equation}

Equivalently, using
\[
S^{(n)}=S^{(n-1)}+\tau D^{(n-1)},
\]
we obtain
\begin{equation}
\label{eq:D_S_coupled}
\begin{aligned}
D^{(n)}
\leq{}&
(1+C_1^*\tau)D^{(n-1)}
+\frac{6C_2^*M_0^2}{\epsilon}\,
\tau S^{(n-1)}
+\tau R,
\end{aligned}
\end{equation}
where
\begin{equation}
\label{eq:R_definition}
\begin{aligned}
R:={}&
C_3^*\left(\frac{L}{H}\right)^4
+C_2^*C_4^*L^{-1}H_0^9H^{-8}
+C_2^*C_5^*L^4H_0^{-4} \\
&+C_2^*C_6^*L^{-1}P^{-1}H_0.
\end{aligned}
\end{equation}

Applying the discrete Gronwall inequality to the coupled quantities
\(D^{(n)}\) and \(S^{(n)}\) and notice that $D^{(0)}=0$, we obtain
\begin{equation}
\label{eq:D_gronwall}
D^{(n)}
\leq
\frac{e^{\lambda t_n}-1}{\lambda}
\Bigg(
C_3^*\left(\frac{L}{H}\right)^4
+C_2^*C_4^*L^{-1}H_0^9H^{-8}
+C_2^*C_5^*L^4H_0^{-4}
+C_2^*C_6^*L^{-1}P^{-1}H_0
\Bigg),
\end{equation}
where
\[
\lambda
=
\frac{
C_1^*
+\sqrt{(C_1^*)^2+24C_2^*M_0^2/\epsilon}
}{2}.
\]

Since the initial errors vanish, \(D^{(0)}=0\). Therefore, setting
\(t_{N_T}=N_T\tau=T\) yields
\begin{equation}
\label{eq:D_final}
\begin{aligned}
D^{(N_T)}
\leq{}&
C_2^2 e^{2C_1T}
\Bigg(
C_3^2\left(\frac{L}{H}\right)^4
+C_4^2L^{-1}H_0^9H^{-8} \\
&\hspace{3.2cm}
+C_5^2L^4H_0^{-4}
+C_6^2L^{-1}P^{-1}H_0
\Bigg),
\end{aligned}
\end{equation}
where the constants are independent of the computational parameters:
\begin{equation}
\begin{aligned}
2C_1
&=
\frac{
C_1^*
+\sqrt{(C_1^*)^2+24C_2^*M_0^2/\epsilon}
}{2}+2,
&
C_2&=(2C_1-2)^{-1/2}+2,\\
C_3^2&=C_3^*,
&
C_4^2&=C_2^*C_4^*,\\
C_5^2&=C_2^*C_5^*,
&
C_6^2&=C_2^*C_6^*.
\end{aligned}
\end{equation}

Finally, from the triangle inequality of \(L^2\) distance,
\begin{equation}\label{eq:convergence_estimate}
\begin{aligned}
\mathbb{E}[w^{(N_T)}]
    &\leq \sqrt{\mathbb{E}[d^{(N_T)}]} + 
    \Biggl( \frac{1}{P} \sum_{p=1}^{P} \bigl| \mathbf{X}_p(T) - \widetilde{\mathbf{X}}_p^{(N_T)} \bigr|^2 \Biggr)^{\!\frac{1}{2}}\\  &\leq \sqrt{D^{(N_T)}} + 
    \Biggl( \frac{1}{P} \sum_{p=1}^{P} \bigl| \mathbf{X}_p(T) - \widetilde{\mathbf{X}}_p^{(N_T)} \bigr|^2 \Biggr)^{\!\frac{1}{2}},    
\end{aligned}
\end{equation}
combining with the temporal error bound from Lemma~\ref{lemma:temporal} (notice that $T^2+T <2e^T$ for $T>0$), we obtain the claimed convergence bound in  \eqref{eq:main_results_thm}.
\end{proof}

\begin{corollary}
\label{cor:optimized_error}
By optimizing the filter resolution \(H_0 = H^{\frac{8}{13}}\), the error bound simplifies to
\begin{equation}
\begin{aligned}
\mathbb{E}\left[w^{\left(N_T\right)}\right] &\leq C_2 e^{C_1T}
\Bigg(
\left(C_4L^{-\frac{1}{2}}+C_5L^2\right)H^{-\frac{16}{13}} \\
&\hspace{3.2cm}
+C_6L^{-\frac{1}{2}}P^{-\frac{1}{2}}H^\frac{4}{13}+C_0 \tau
\Bigg).
\label{eq:main_results_opt}
\end{aligned}
\end{equation}
\end{corollary}

\subsection{Relaxation of the smoothness assumption}
\label{subsec:c_hat_relax}
In this section, we establish several conclusions showing that, under specific parameter combinations, the smoothness of $\hat{c}$ in Assumption \ref{as} can be naturally guaranteed by the algorithm itself.

Since the proof of Theorem~\ref{thm:overall_error} yields $D^{(n)}\lesssim H_0^{-4}$ for the optimal parameters \(H_0 \sim H^\frac{8}{13}\), the following lemma implies that $\hat{c}$ naturally exhibits third-order convergence over finite time.

\begin{lemma}
\label{lemma:c_hat_2norm_reg}
Let the computational parameters satisfy \( P \geq H^{\frac{40}{13}}(1+L^{-5}) \) and \( H_0 = H^{\frac{8}{13}} \). Under the assumption that the total error energy \(D^{(s)} \leq Q H_0^{-4}\) for \(1 \leq s \leq n\) and constant $Q$ independent of \(\tau,H,H_0,P\), the algorithm-derived \(\hat{c}^{(n)}\) has a square-integrable third derivative:
\begin{equation}
\max_{|\mathbf{k}|= 3}\| D_{\mathbf{k}} \hat{c}^{(n)} \|_2 \leq 2L^\frac{3}{2}M_3^*+\overline{C}_1Q+\overline{C}_2,
\end{equation}
where $M_3^* := \max_{|\mathbf{k}|= 3}
\left\|D_{\mathbf{k}}c(\mathbf{x},t)\right\|_\infty$ is the ground-truth smoothness bound of the third order, and $\overline{C}_1,\overline{C}_2$ is independent of \(\tau,H,H_0,P\).
\end{lemma}
\begin{proof}
We show that $\|D_{\mathbf{k}} \hat{c}\|_2$ is controlled by $\|D_{\mathbf{k}} c\|_2$. By Parseval's theorem, 
\begin{equation}
\begin{aligned}
\sum_{|\mathbf{k}|=3}& \int_{\Omega} |D_{\mathbf{k}} \hat{c}^{(n)}( \mathbf{x})|^2 \, \mathrm{d}\mathbf{x} = L^3\sum_{|\mathbf{k}|=3} \sum_{\mathbf{q}\in U} \left(\mathbf{q}\frac{2\pi}{L} \right)^{2\mathbf{k}}\left|\Phi_{\mathbf q}\hat{\alpha}^{(n)}_{\mathbf{q}}\right|^2\\
&\leq 2L^3\sum_{|\mathbf{k}|=3} \sum_{\mathbf{q}\in U} \left(\mathbf{q}\frac{2\pi}{L} \right)^{2\mathbf{k}}\left(\left|\Phi_{\mathbf q}\alpha^{(n)}_{\mathbf{q}}\right|^2 + \left|\Phi_{\mathbf q}\hat{\alpha}^{(n)}_{\mathbf{q}}-\Phi_{\mathbf q}\alpha^{(n)}_{\mathbf{q}}\right|^2 \right)\\
&\leq 2\sum_{|\mathbf{k}|=3} \int_{\Omega} |D_{\mathbf{k}} c^{(n)}(\mathbf{x})|^2 \, \mathrm{d}\mathbf{x} + 20L^3 \sum_{\mathbf{q}\in U} \left|\mathbf{q}\frac{2\pi}{L} \right|^6\left|\Phi_{\mathbf q}\hat{\alpha}^{(n)}_{\mathbf{q}}-\Phi_{\mathbf q}\alpha^{(n)}_{\mathbf{q}}\right|^2.
\end{aligned}
\end{equation}

Similar to Lemmas \ref{lem:grad_c} and \ref{lem:grad_c2}, we obtain the following decomposition estimate:
\begin{equation}
  \begin{aligned}
    \sum_{\mathbf{q} \in U} \left| \mathbf{q} \frac{2\pi}{L}\right|^6\left| \Phi_{\mathbf q}\cdot \hat{\epsilon}_{\mathbf{q}}^{(n)} \right|^2 
    &\leq M_0^2 L^{-6} \sum_{\mathbf{q} \in U} \left(\frac{C\left|\mathbf{q} \frac{2\pi}{L}\right|^{11}}{\left|\mathbf{q} \frac{2\pi}{L}\right|^2 + k^2} \left( \frac{L}{H} \right)^4 e^{-\frac{2\pi ^2|\mathbf q|^2}{H_0^2}}\right) ^2 \\
    &\leq H_0^3 \cdot M_0^2 C^2 L^{-6} \left( H_0 \frac{2\pi}{L} \right)^{10} \left( \frac{L}{H} \right)^8,
  \end{aligned}
  \label{eq:energy_sum_low_freq13}
\end{equation}

\begin{equation}
  \begin{aligned}
    &\sum_{\mathbf{q} \in U} \left| \mathbf{q} \frac{2\pi}{L}\right|^6\left|\Phi_{\mathbf q} \cdot \check{\epsilon}_{\mathbf{q}}^{(n)} \right|^2 
    \leq \frac{M_0^2 \tau}{\epsilon} \sum_{\mathbf{q} \in U} \frac{|\mathbf{q}\frac{2\pi}{L}|^6}{|\mathbf{q}\frac{2\pi}{L}|^2 + k^2} \sum_{s=0}^{n-1} \left|\Phi_{\mathbf q}\nu_{\mathbf{q}}^{(s)}\right|^2 \\
    &\leq \frac{M_0^2 \tau}{\epsilon} \sum_{s=0}^{n-1} \sum_{\mathbf{q} \in U} \left| \mathbf{q} \frac{2\pi}{L}\right|^4 \left|\Phi_{\mathbf q}\nu_{\mathbf{q}}^{(s)}\right|^2\leq \frac{M_0^2 \tau}{\epsilon L^3} \sum_{s=0}^{n-1} \|D^2  G_0 \ast\mathcal{E}^{(s)}\|_2^2\\
    &\leq \frac{M_0^2 \tau H_0^4}{\epsilon L^7} \sum_{s=0}^{n-1} \|G_0 \ast\mathcal{E}^{(s)}\|_2^2 =\frac{M_0^2 \tau H_0^4}{\epsilon L^7} \sum_{s=0}^{n-1} h^{(s)},
  \end{aligned}
  \label{eq:energy_sum_low_freq23}
\end{equation}
\begin{equation}
\begin{aligned}
\mathbb{E}\left[\sum_{\mathbf{q} \in U} \left| \mathbf{q}\frac{2\pi}{L}\right|^6\left|\Phi_{\mathbf q}(\alpha^{(n)}_{\mathbf{q}} - \tilde{\alpha}^{(n)}_{\mathbf{q}}) \right|^2 \right]
&\leq \sum_{\mathbf{q} \in U} \frac{|\mathbf{q}\frac{2\pi}{L}|^6|\Phi_{\mathbf q}|^2}{(|\mathbf{q}\frac{2\pi}{L}|^2 + k^2)^2} \cdot \frac{M_0^2}{PL^6}\\
&\leq \frac{M_0^2 H_0^5}{P L^8}.
\end{aligned}
\label{eq:random_low_freq3}
\end{equation}

Note that none of these conclusions requires any smoothness assumption on $\hat{c}^{(n)}$, and $\|D_{\mathbf{k}} c^{(n)}\|_2 \leq L^{3/2} \|D_{\mathbf{k}} c^{(n)}\|_\infty$. Under our assumptions on the computational parameters, we have $H_0^{13}H^{-8}=1,H_0^5 P^{-1} \leq 1$. Taking $\overline{C}_1=20M_0^2 T \epsilon^{-1}L^{-4},\overline{C}_2= (\frac{513}{16})^2 (2\pi)^{10}M_0^2 L^{-5}$ then yields the desired estimate.
\end{proof}

The third-order regularity assumption on $\hat{c}$ can be removed by directly incorporating the weighted error term $H_0^{-4}\max_{|\mathbf{k}|=3}\lVert D_{\mathbf{k}}(\hat{c}^{(n)}-c^{(n)})\rVert_2$ into the Gronwall estimate. However, replacing the $L^\infty$-norm of the third derivative of $\hat{c}$ by its $L^2$-norm in the estimates of Lemmas~\ref{lem:Gronwall2} and \ref{lemma:h_bound} would be technically involved. We therefore state the following corollary, which shows that, under more restrictive conditions of the computational parameters, the third derivative of $\hat{c}$ is unconditionally bounded in the $L^\infty$-norm. Consequently, the assumption on $\hat{c}$ in Assumption~\ref{as} can be removed under these parameters.

\begin{corollary}
\label{cor:c_hat_2norm_reg}
Let the computational parameters satisfy \( P \geq H^4(1+L^{-5}) \), \( H_0 \leq H^{\frac{1}{2}} \), the underlying $c^{(n)}$ satisfies $\max_{|\mathbf{k}|= 5}\| D_{\mathbf{k}} c^{(n)}\|_\infty \leq M_5^*$. The algorithm-derived \(\hat{c}^{(n)}\) has a bounded third derivative:
\begin{equation}
\max_{|\mathbf{k}|= 3}\| D_{\mathbf{k}} \hat{c}^{(n)}\|_\infty \leq M_3^*+\overline{C}_0,
\end{equation}
where $M_3^* := \max_{|\mathbf{k}|= 3}
\left\|D_{\mathbf{k}}c(\mathbf{x},t)\right\|_\infty$ is the ground-truth smoothness bound of the third order, and $\overline{C}_0$ is independent of \(\tau,H,H_0,P\).   
\end{corollary}
\begin{proof}
Under the fifth-order smoothness assumption of $c$, the source term of the error \eqref{eq:high_freq} becomes $S_1^{(n)}\leq M_5^2\left(\frac{L}{H_0}\right)^8$. Substituting this into the error accumulation estimate in \eqref{eq:main_results_thm}, and imposing the additional parametric condition $H\geq H_0^2$ and $P \geq H_0^8$, we arrive at $D^{(n)} \lesssim H_0^{-7}.$ Consequently, we have the following control:
\begin{equation}
\begin{aligned}
\max_{|\mathbf{k}|= 3}\| D_{\mathbf{k}} \hat{c} \|_\infty &\leq \max_{|\mathbf{k}|= 3}\| D_{\mathbf{k}} (G_0 \ast c) \|_\infty +\max_{|\mathbf{k}|= 3}\| D_{\mathbf{k}} (\hat{c} -G_0 \ast c) \|_\infty\\
&\leq \max_{|\mathbf{k}|= 3}\| D_{\mathbf{k}} c \|_\infty + \|G_0\|_2\left(\overline{C}_1 H_0^4D^{(n)}\right)^\frac{1}{2}\\
&\leq M_3^* +\overline{C}_1 H_0^\frac{7}{2} \left(D^{(n)}\right)^\frac{1}{2},
\end{aligned}
\end{equation}
which yields the desired conclusion.
\end{proof}

\begin{color}{black}
\begin{remark}
\label{remark:blowup}
When the underlying solution remains in the non-blow-up regime, the error analysis shows that, if the algorithmic parameters $P$, $H$, and $N_T=T/\tau$ are increased according to the prescribed coupled scaling, the numerical solution converges to the underlying solution at a quantified rate. Consequently, for sufficiently refined parameters, the method cannot produce an artificial blow-up when the underlying solution remains regular. In this sense, the method provides a one-sided blow-up diagnostic: a numerical blow-up that persists under systematic parameter refinement indicates that the underlying system cannot remain in the non-blow-up regime.

The converse implication need not hold at a fixed finite resolution. Since the dominant numerical effect of a convergent discretization is expected to be typically dissipative rather than aggregative, the method may suppress or delay a genuine blow-up. Therefore, the numerically detected critical mass may be slightly larger than the true critical value. Nevertheless, by progressively increasing $P$, $H$, and $N_T$ according to the prescribed scaling, the numerical dissipation can be reduced and the detected blow-up threshold is expected to approach the true threshold. Thus, the proposed method should be understood as a refinement-based blow-up detection procedure, rather than as a convergence result at the singularity itself.
\end{remark}
\end{color}

\section{Numerical experiments}
\label{sec:Num}
This section presents numerical results to validate the convergence analysis and demonstrate the performance of the SIPF-PIC method. We first establish convergence properties through systematic analysis of particle-grid resolution dependencies, followed by demonstrations of the framework's capacity to resolve finite-time blow-up dynamics.

The section is organized as follows. Section \ref{subsec:rad} quantifies algorithmic convergence through error analysis with respect to the particle number \( P \) and grid resolution \( H \) per dimension, and subsequently compares computational performance using identical Keller-Segel parameters. Building on these foundations, Section \ref{subsec:critical_mass} verifies the theoretical threshold \( M_{\text{critical}} = 8\pi \) in 2D elliptic systems and identifies an empirical threshold in 3D parabolic systems. We then explore blow-up dynamics: Section \ref{subsec:tetra} reveals three mass-dependent phases in tetrahedral symmetric systems, including subcritical dispersion, intermediate unified blowup, and supercritical multistage singularity formation. Finally, Section \ref{subsec:ring_blowup} demonstrates radial concentration into singular ring structures that exhibit post-collapse instabilities. \begin{color}{black}More specifically, the convergence tests in Section~\ref{subsec:rad} do not exhibit blow-up and lie within the scope of the convergence theorem, whereas the blow-up examples in the subsequent subsections are presented as numerical evidence beyond the proven regime.
\end{color}

\subsection{Convergence of the SIPF-PIC method with respect to discretization}\label{subsec:rad}
We first validate the convergence of the algorithm with respect to the particle number $P$ and grid resolution $H$ per dimension. Consider a spherically symmetric initial configuration where $\rho_0(\mathbf{x})$ is uniformly distributed within a unit ball centered at the origin:  
\begin{equation}
\label{radius}
\rho_0(\mathbf{x}) = \begin{cases} 
\frac{60}{\pi} & |\mathbf{x}| \leq 1, \\
0 & \text{otherwise},
\end{cases}
\end{equation}
with total mass $M_0 = 80$. The model parameters in~\eqref{eq:para_system} are set as $\mu = \chi = 1$, $\epsilon = 10^{-4}$, $k = 10^{-1}$, with final time $T = 2 \times 10^{-2}$, timestep $\tau = 10^{-5}$, and total time steps $N_T=2000$.  

By exploiting radial symmetry, we can reduce the 3D system to a 1D radial coordinate $r$ through the following transformation 
\begin{equation}
x_1 = r \cos\theta,\quad x_2 = r \sin\theta\cos\phi,\quad x_3 = r \sin\theta\sin\phi.
\end{equation}
This yields the simplified system:  
\begin{equation}
    \label{eq:1D}
    \begin{cases}
        \displaystyle
        \rho_t = \frac{2\mu}{r}\rho_r + \mu\rho_{rr} - \frac{2\chi}{r}\rho c_r - \chi\rho_r c_r - \chi\rho c_{rr}, \\
        \displaystyle
        \epsilon c_t = \frac{2}{r}c_r + c_{rr} - k^2c + \rho,
    \end{cases}
\end{equation}
together with boundary conditions \(\rho_r|_{r=0}=c_r|_{r=0}=0\).

While the problem defined in~\eqref{radius} is posed in $\mathbb{R}^3$, the particle dynamics remain localized within $|\mathbf{x}| < 2$ during the entire simulation; \begin{color}{black} the region contains at least $99.9\%$ of the mass of the solution (the $0.999$-quantile $r_{0.999}\approx 1.45047$). Thus, the statistics used in the comparison are essentially unaffected by the far-field boundary.

Under radial symmetry, the 3D problem reduces to a 1D radial problem on a ball. By contrast, the fully three-dimensional particle-field calculation is performed on the periodic cube. Although the spherical no-flux and cubic periodic boundary conditions cannot be matched exactly, their discrepancy is reduced by choosing the radial-domain radius $R_{\mathrm{ball}}$ so that the two computational domains have the same volume, namely $\frac{4\pi}{3}R_{\mathrm{ball}}^3=L^3$. For $L=8$, this gives $R_{\mathrm{ball}}\approx4.9628$. A reference solution can then be computed using an implicit finite volume method to solve \eqref{eq:1D} on a high-resolution 1D grid with $N = 4 \times 10^5$ points, utilizing the radial symmetry.\end{color}

Convergence studies for the SIPF-PIC method are conducted by varying $P$ and $H$ while comparing 
the numerical results with this reference solution. We use the cumulative distribution function
\begin{equation}
m(r,t;\rho) = \frac{1}{M_0} \int_{|\mathbf{x}| \leq r} \rho(\mathbf{x},t) \, \mathrm{d}\mathbf{x}
\end{equation}
at the final time $T = 2 \times 10^{-2}$ as the comparative criterion. Since $m(0,T;\rho)=0$ and $m(\infty,T;\rho)=1$ for every mass-preserving density function $\rho$, we compute the $W^1$ distance between the results obtained by the SIPF-PIC method and the reference solution as the error metric. Specifically, for the density \( \hat{\rho} \) defined by \eqref{eq:def_rho_hat}, we have
\(
m(r,T,\hat{\rho}) = \frac{1}{P} \sum_{p=1}^{P} \mathbf{1}_{\{|\mathbf{X}_p^{(N_T)}| \leq r\}},
\)
and for a radially symmetric density \( \rho(\mathbf{x},t) = \psi(|\mathbf{x}|,t) \), we have
\(
m(r,T,\rho) = \frac{4\pi} {M_0}\int_0^r \psi(s,t) s^2 \, ds.
\)

To quantify the difference between the two distributions, we use a truncated quantile-based Wasserstein discrepancy. Let $n=10^3$ and take $m=\left\lceil \frac{3}{4}n\right\rceil-1$. We compute the quantiles $q_{j/n}$ and $\hat q_{j/n}$ for $j=0,1,\ldots,m$, and define
\begin{equation}
\mathcal{E}
:=
\frac{1}{m+1}
\sum_{j=0}^{m}
\left|q_{j/n}-\hat q_{j/n}\right|.
\end{equation}
This quantity is a truncated analogue of the one-dimensional $W_1$ distance. The upper quartile is excluded because it is more sensitive to the discrepancy between the spherical no-flux boundary used for the radial reference solution and the periodic boundary used for the three-dimensional particle simulation. The truncated discrepancy therefore emphasizes agreement in the major part of the distribution, where boundary effects are negligible.

We test all combinations of $P \in \{2^{11}, \ldots, 2^{20}\}$ and $H \in \{8, 16, 32, 64, 128, 256\}$. \begin{color}{black}In these simulations, we set $H_{0}=\left\lceil 8H^{8/13}L^{5/13}\right\rceil$, where $H$ is the spatial grid resolution. Accordingly, the Gaussian spectral filter $\Phi_{\mathbf q}=\exp\left(-2\pi^2|\mathbf q|^2/H_{0}^2\right)$ is applied to the particle density after the particle-to-grid projection and before the field solve. This filter corresponds to convolution with a Gaussian having standard deviation $L/H_{0}$ in each spatial direction. For every parameter configuration, we perform $N_{\mathrm{sim}}=30$ statistically independent simulations, using mutually independent random seeds for both the initial particle sampling and the Brownian increments. The truncated quantile-based Wasserstein discrepancy is evaluated separately for each independent run and then averaged over all runs, namely, \(\overline{\mathcal{E}}:=N_{\mathrm{sim}}^{-1}\sum_{\ell=1}^{N_{\mathrm{sim}}}\mathcal{E}_{\ell}\).\end{color}

\begin{figure}[htbp]
    \centering
    \begin{subfigure}[b]{\textwidth}
        \centering
        \begin{subfigure}[b]{0.42\textwidth}
            \includegraphics[width=\textwidth]{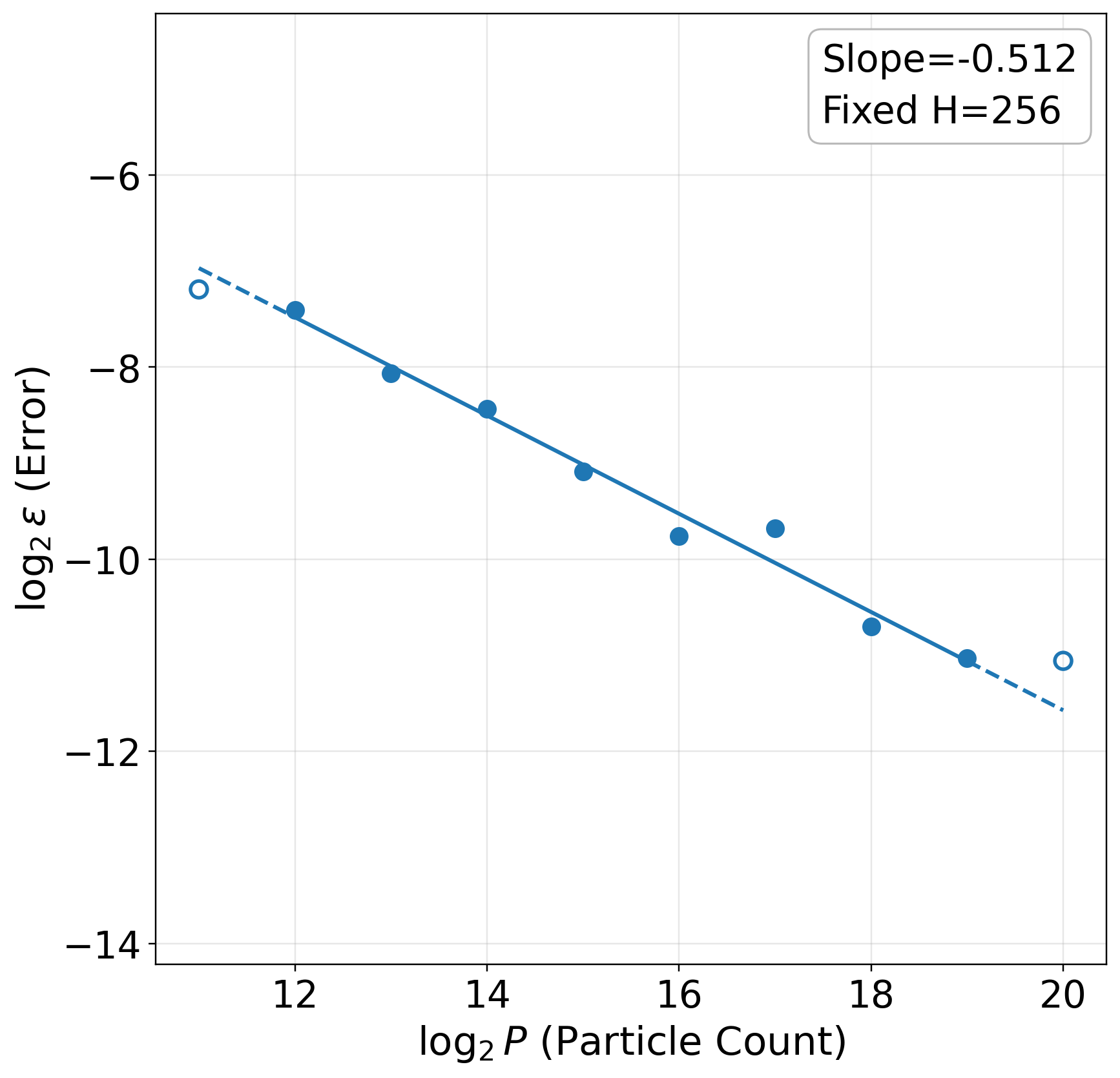}
            \caption{Error convergence against particle number ($P$) with fixed grid resolution.}
            \label{fig:convergence_4_2_P}
        \end{subfigure}
        \hfill
        \begin{subfigure}[b]{0.42\textwidth}
            \includegraphics[width=\textwidth]{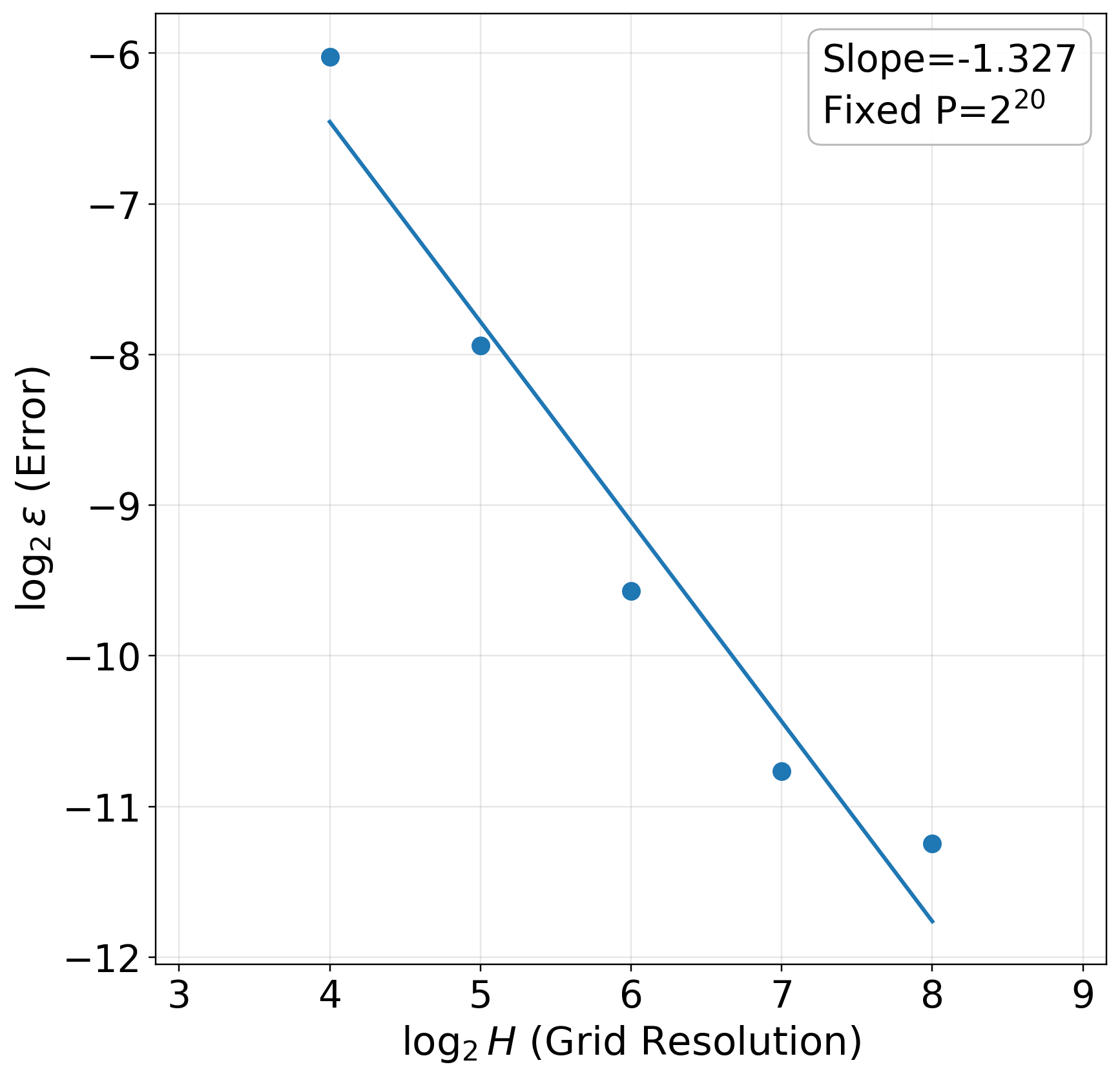}
            \caption{Error convergence against grid resolution ($H$) with fixed particle number.}
            \label{fig:convergence_4_2_H}
        \end{subfigure}
    \end{subfigure}
    \caption{Convergence characteristics under the dual parameters of particle number and grid resolution.}
    \label{fig:convergence2}
\end{figure}

Figure~\ref{fig:convergence2} demonstrates the convergence behavior of the SIPF-PIC method using fourth-order particle-to-grid and second-order grid-to-particle interpolations, quantifying convergence rates with respect to the individual parameters. Two complementary studies isolate the parametric dependencies: fixing the grid resolution at $H_{\max}$ yields the particle number convergence rate via the slope of $\log_2 \mathcal{E}$ versus $\log_2 P$, with a measured value of $-0.512$, which closely matches the theoretical order of $1/2$ derived in Section~\ref{sec:Error}; fixing $P_{\max}$ determines the spatial convergence order through the slope of $\log_2 \mathcal{E}$ versus $\log_2 H$, with a measured value of $-1.327$, which closely matches the theoretical order of $16/13$. \begin{color}{black}For each parameter configuration, the standard deviation across independent runs is less than one half of the mean error, and more than \(90\%\) of the random trajectories have errors less than twice the mean error.\end{color}

\begin{figure}[htbp]
    \centering
    \begin{subfigure}[b]{\textwidth}
        \centering
        \begin{subfigure}[b]{0.42\textwidth}
            \includegraphics[width=\textwidth]{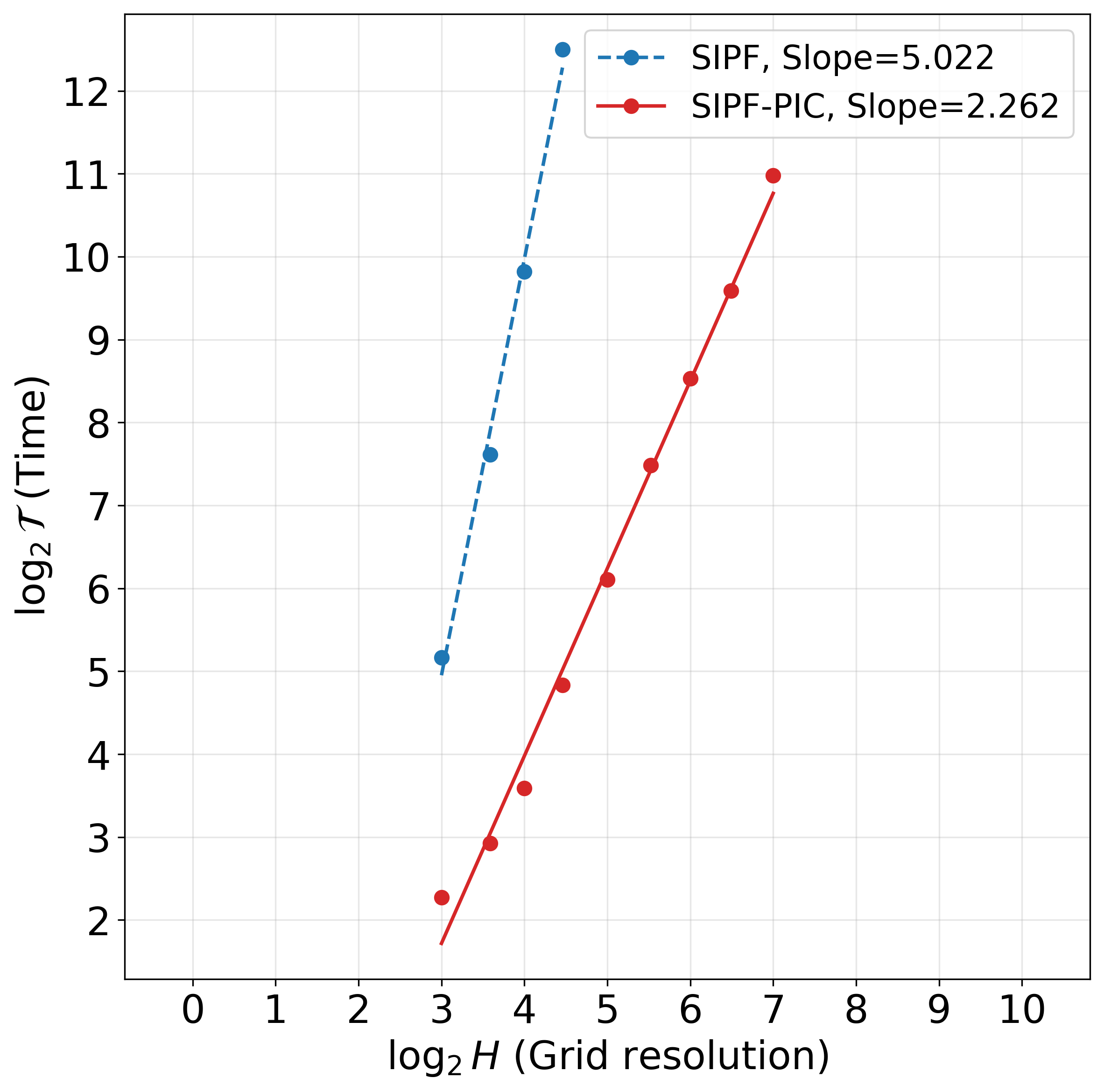}
            \caption{Running time convergence against grid resolution ($H$) with coupled particle number $P=H^3$.}
            \label{fig:cost_error1}
        \end{subfigure}
        \hfill
        \begin{subfigure}[b]{0.42\textwidth}
            \includegraphics[width=\textwidth]{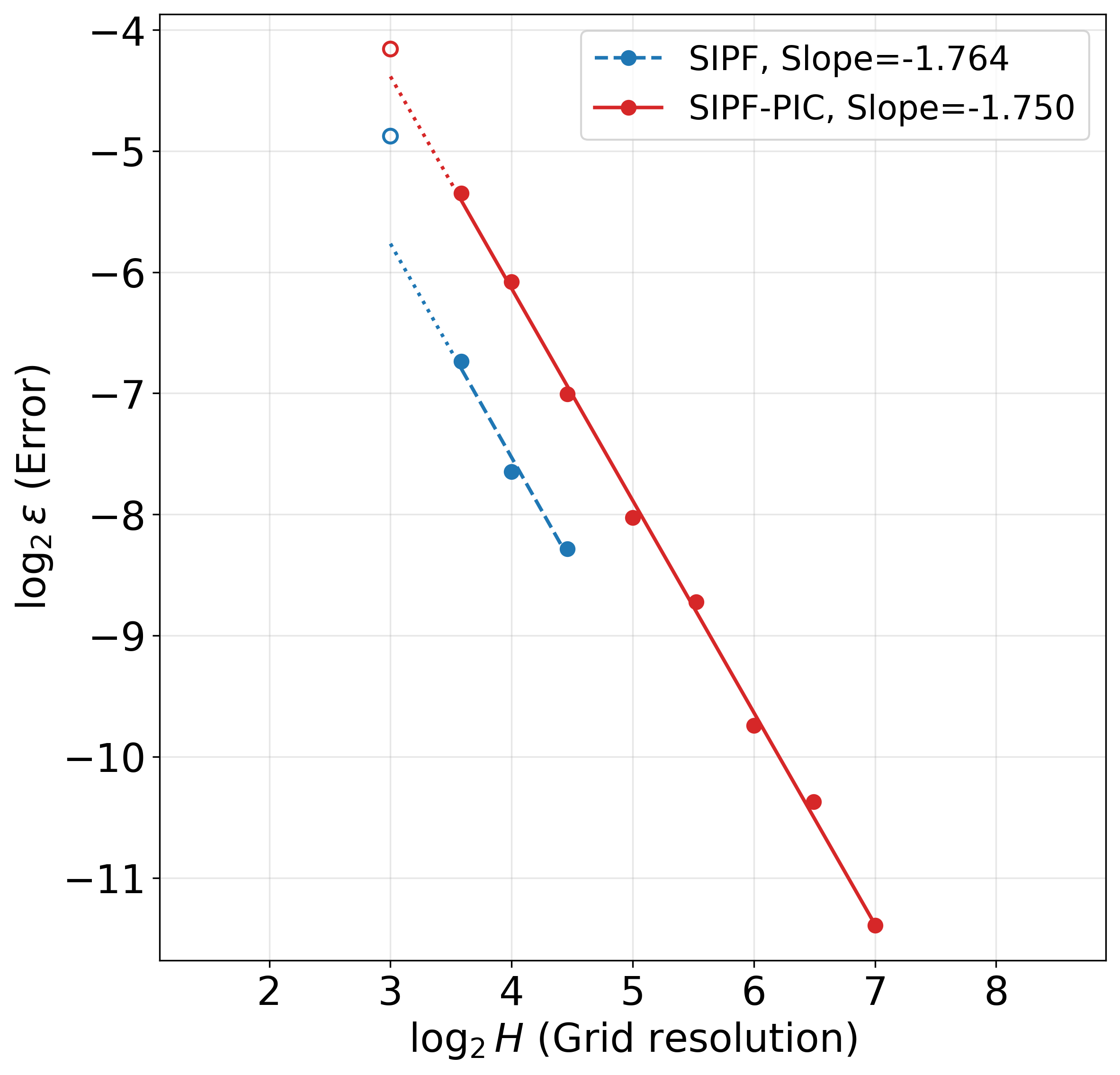}
            \caption{Error convergence against grid resolution ($H$) with coupled particle number $P=H^3$.}
            \label{fig:cost_error2}
        \end{subfigure}
    \end{subfigure}
    \caption{Time-error profiles of SIPF and SIPF-PIC. The coupled particle and grid resolutions satisfy \(P=H^3\).}
    \label{fig:cost_error}
\end{figure}

In addition to the convergence results presented in Figure~\ref{fig:convergence2}, we investigate a 
performance comparison between the original SIPF method and our new SIPF-PIC method through timing and error measurements for the KS system, using identical temporal parameters and initial conditions to those in Figure~\ref{fig:convergence2}. The test cases employ grid resolutions $H \in \{8, 12, 16, 22, 32, 46, 64, 90, 128\}$ and particle numbers $P =H^3$, with results shown in Figure~\ref{fig:cost_error}. The Gaussian filter is set as $H_{0}=\left\lceil 8H^{8/13}L^{5/13}\right\rceil$ (only for SIPF-PIC). The error is measured using the same metric as in the preceding experiment. \begin{color}{black}The reported computational time is the mean of the $N_{\mathrm{sim}}=30$ individual running times and is measured on a personal computer equipped with an Intel Core i9 CPU. The missing entries indicate that the original SIPF method could not complete the computation because of either excessive running time or insufficient memory.\end{color}

The complexity reduction from $O(P H^3)$ to $O(P + H^3 \log H)$ through the SIPF-PIC method's localized particle-grid interpolations and spectral field solutions fundamentally decouples particle-field interactions, achieving accelerated parameter exploration while preserving numerical accuracy and avoiding the reintroduction of $P H^3$ bottlenecks. This complexity separation permits high-resolution simulations at extreme resolutions (e.g., $H\geq256$ and  $P\geq2^{20}$) with controlled error growth, overcoming prior SIPF limitations via adaptive accuracy-efficiency tradeoffs.

\subsection{Critical mass in 2D and 3D systems}  
\label{subsec:critical_mass}  

Recall the 2D parabolic-elliptic KS system with parameters \((k, \epsilon, \mu, \chi) = (0, 0, 1, 1)\), where the critical mass threshold \(M_{\text{crit}} = 8\pi\) governs finite-time blow-up (see Section~\ref{subsec:review}). Numerical simulations initialize with uniform distributions on the unit disk \(\rho_0(\mathbf{x}) = \frac{M_0}{\pi}\mathbf{1}_{|\mathbf{x}| \leq 1}\), for masses \(M_0 \in \{22.0, 23.0, \ldots, 27.0\}\), discretized in the spatial domain \(\Omega = [-4,4]^2\) ($L=8$) with grid resolution \(H = 2048\) on each axis and \(P = 2^{22}\) particles, utilizing fourth-order particle-to-grid and second-order grid-to-particle interpolations, coupled with the Gaussian filter $H_{0}=\left\lceil 8H^{8/13}L^{5/13}\right\rceil$, and the time step \(\tau = 5 \times 10^{-4}\) up to final time \(T = 2.0\), totaling \( N_T = 4000 \) steps.

\begin{figure}[htbp]
  \centering
  \subfloat[$M_0 = 22.0$]{\includegraphics[width=0.3\textwidth]{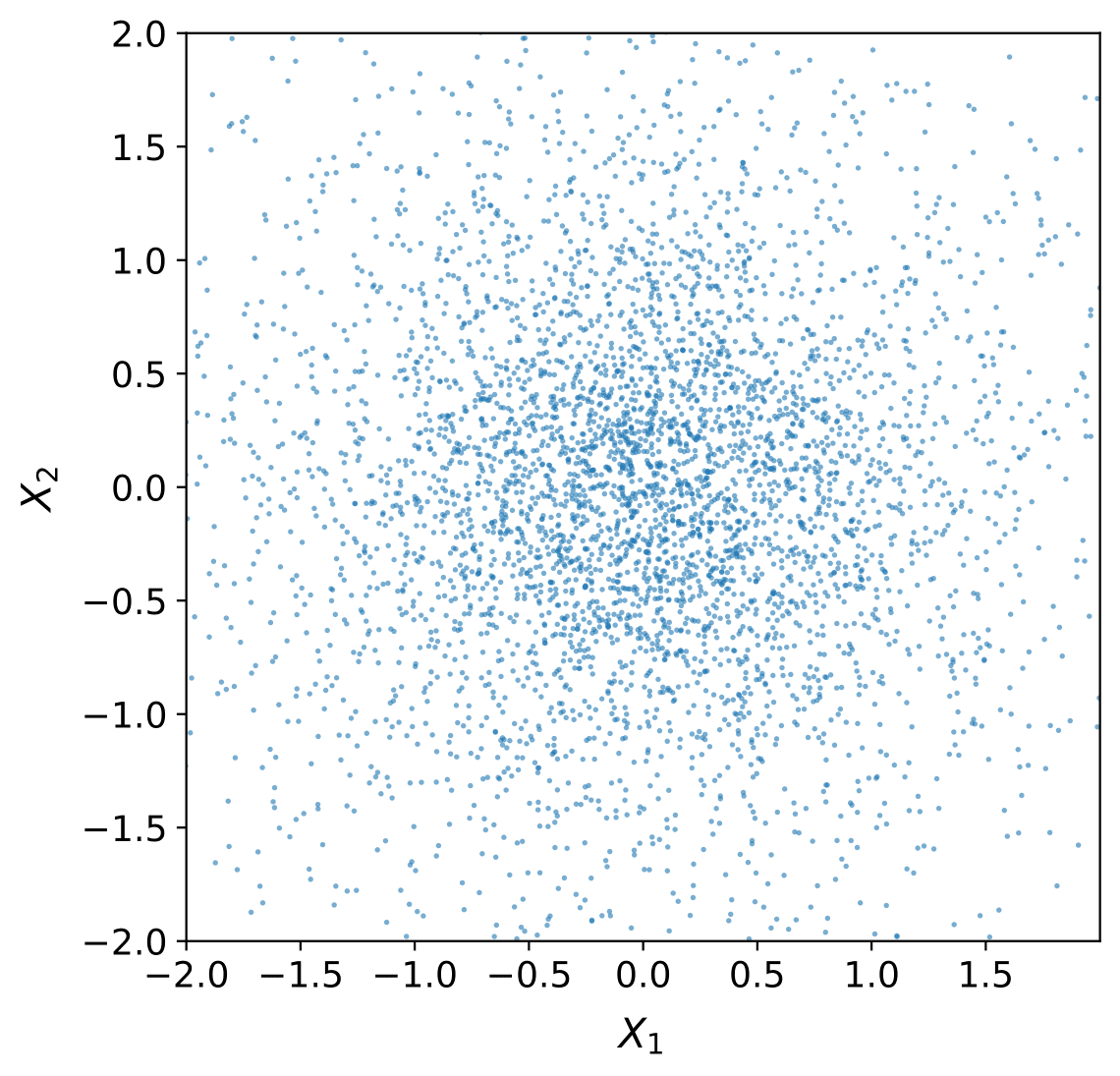}\label{fig:M22_0}}
  \hfill
  \subfloat[$M_0 = 23.0$]{\includegraphics[width=0.3\textwidth]{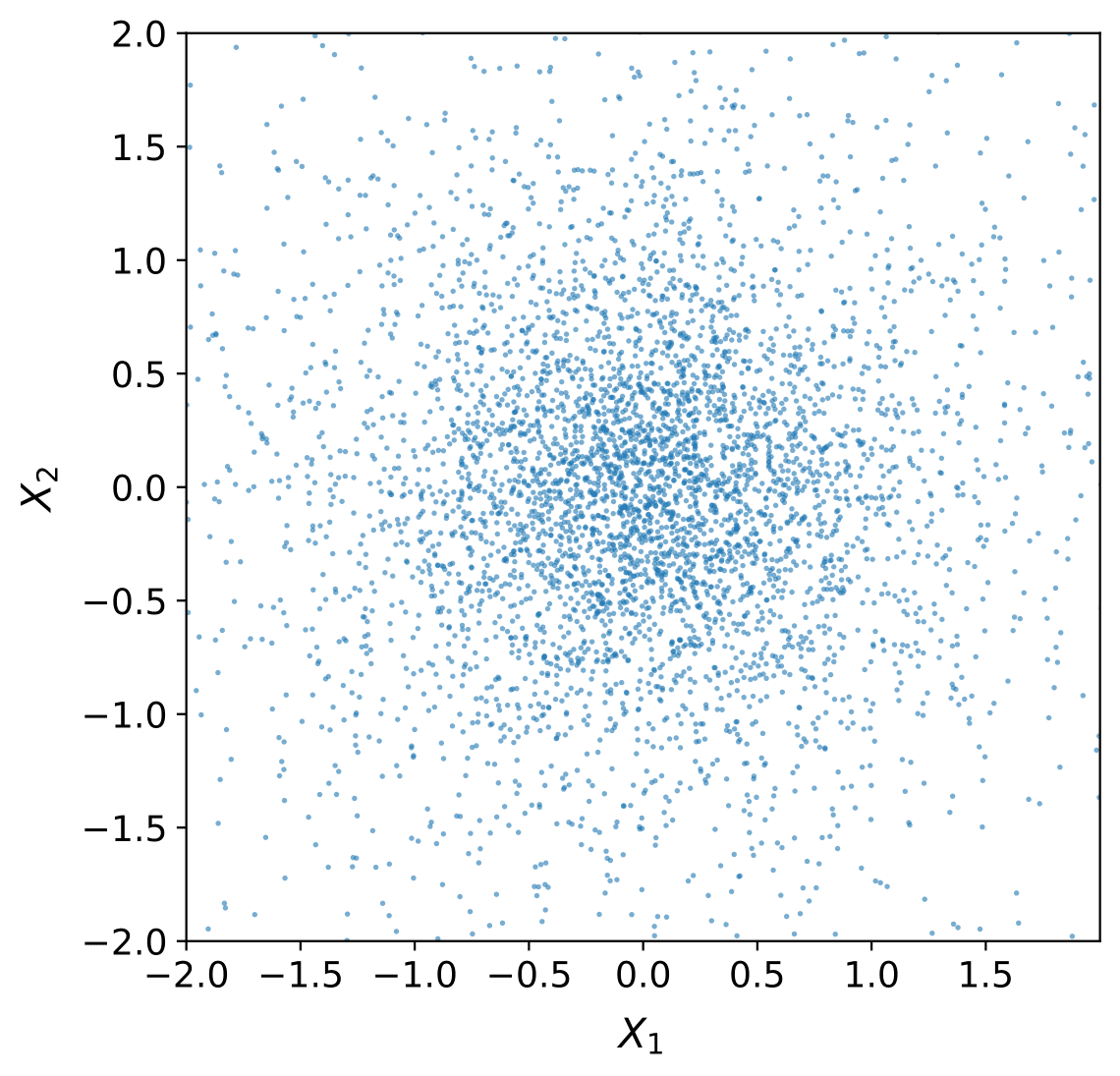}\label{fig:M23_0}}
  \hfill
  \subfloat[$M_0 = 24.0$]{\includegraphics[width=0.3\textwidth]{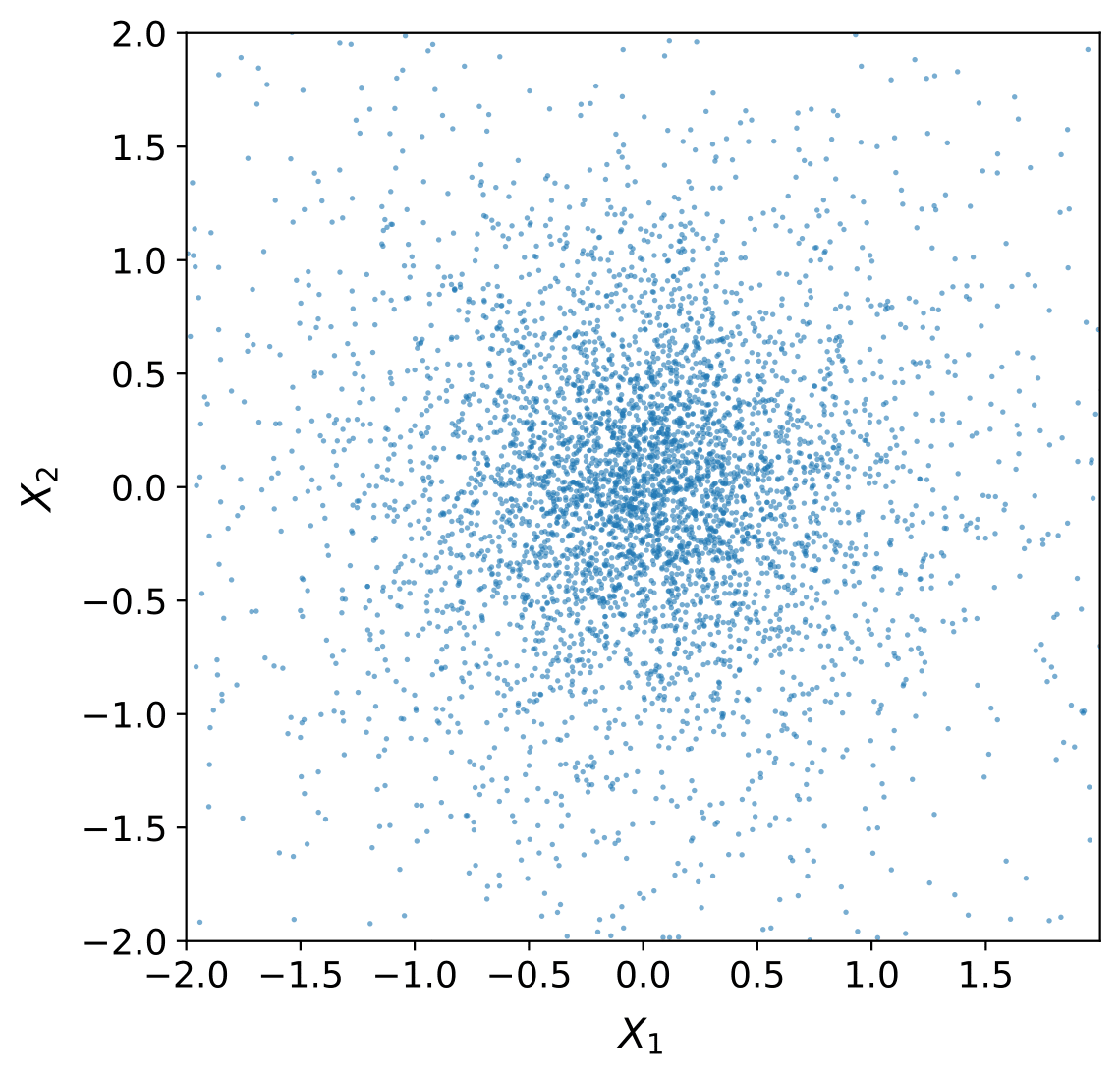}\label{fig:M24_0}}
  
  \vspace{0.5cm}
  \centering
  \subfloat[$M_0 = 25.0$]{\includegraphics[width=0.3\textwidth]{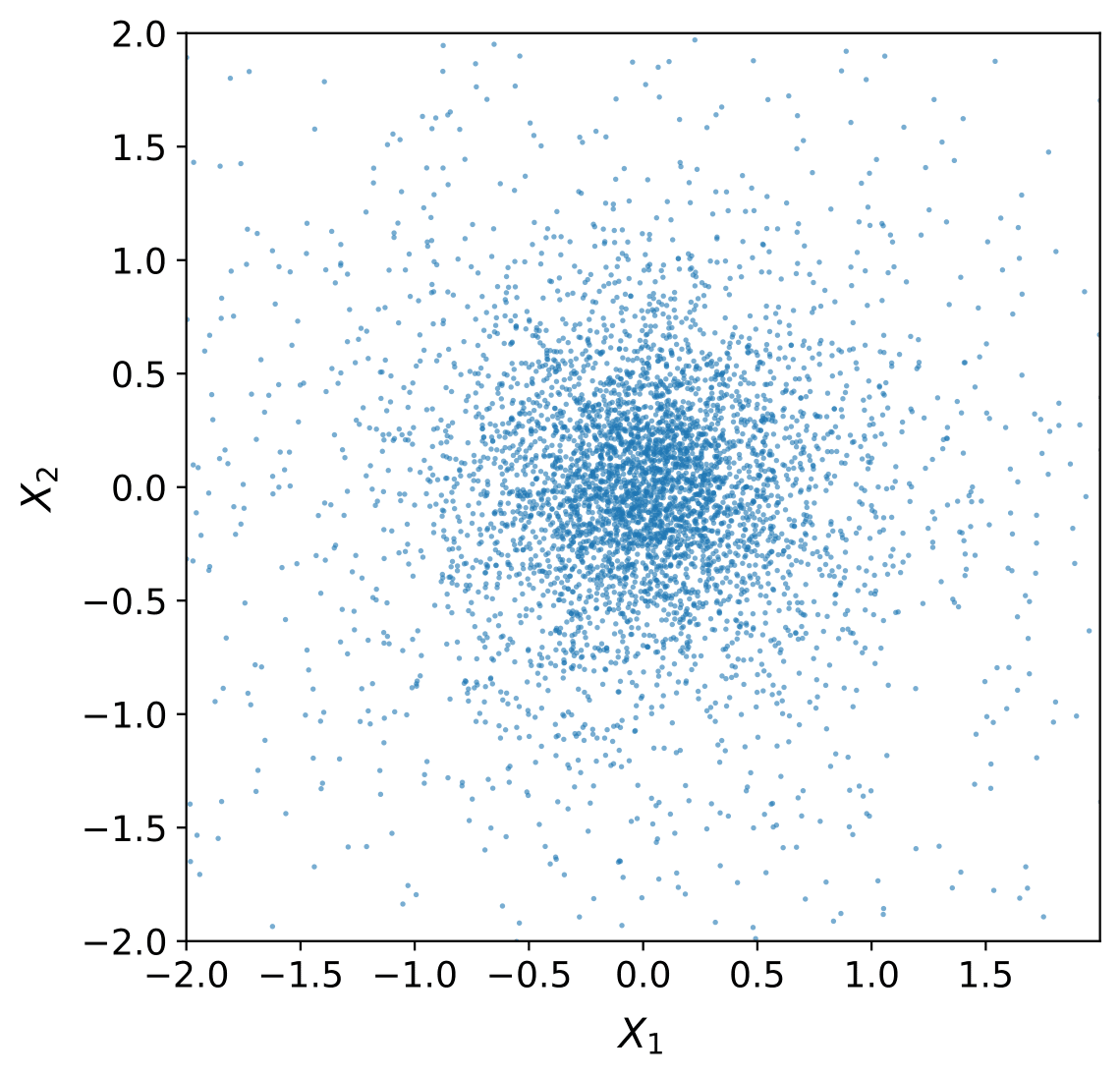}\label{fig:M25_0}}
  \hfill
  \subfloat[$M_0 = 26.0$]{\includegraphics[width=0.3\textwidth]{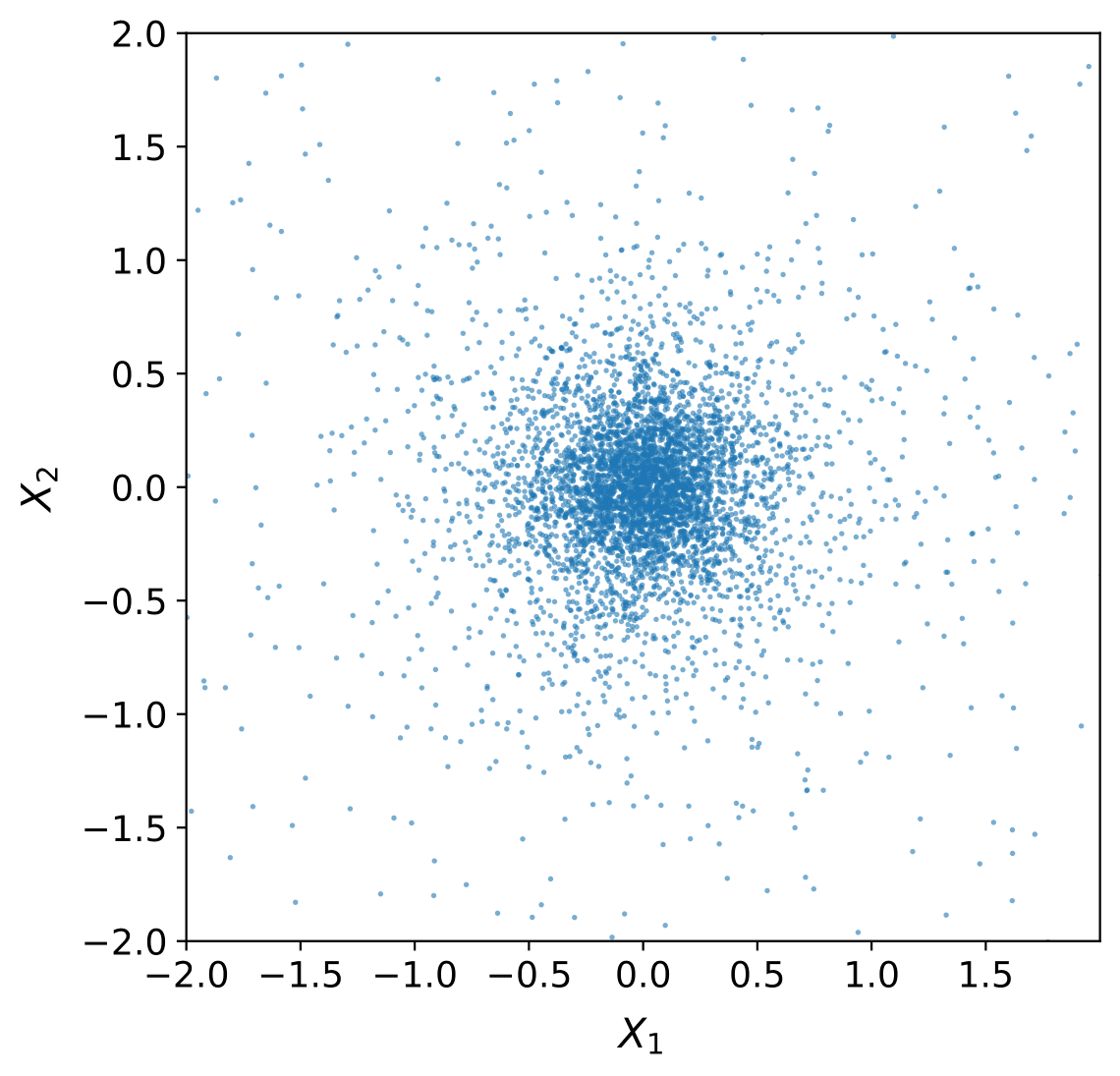}\label{fig:M26_0}}
  \hfill
  \subfloat[$M_0 = 27.0$]{\includegraphics[width=0.3\textwidth]{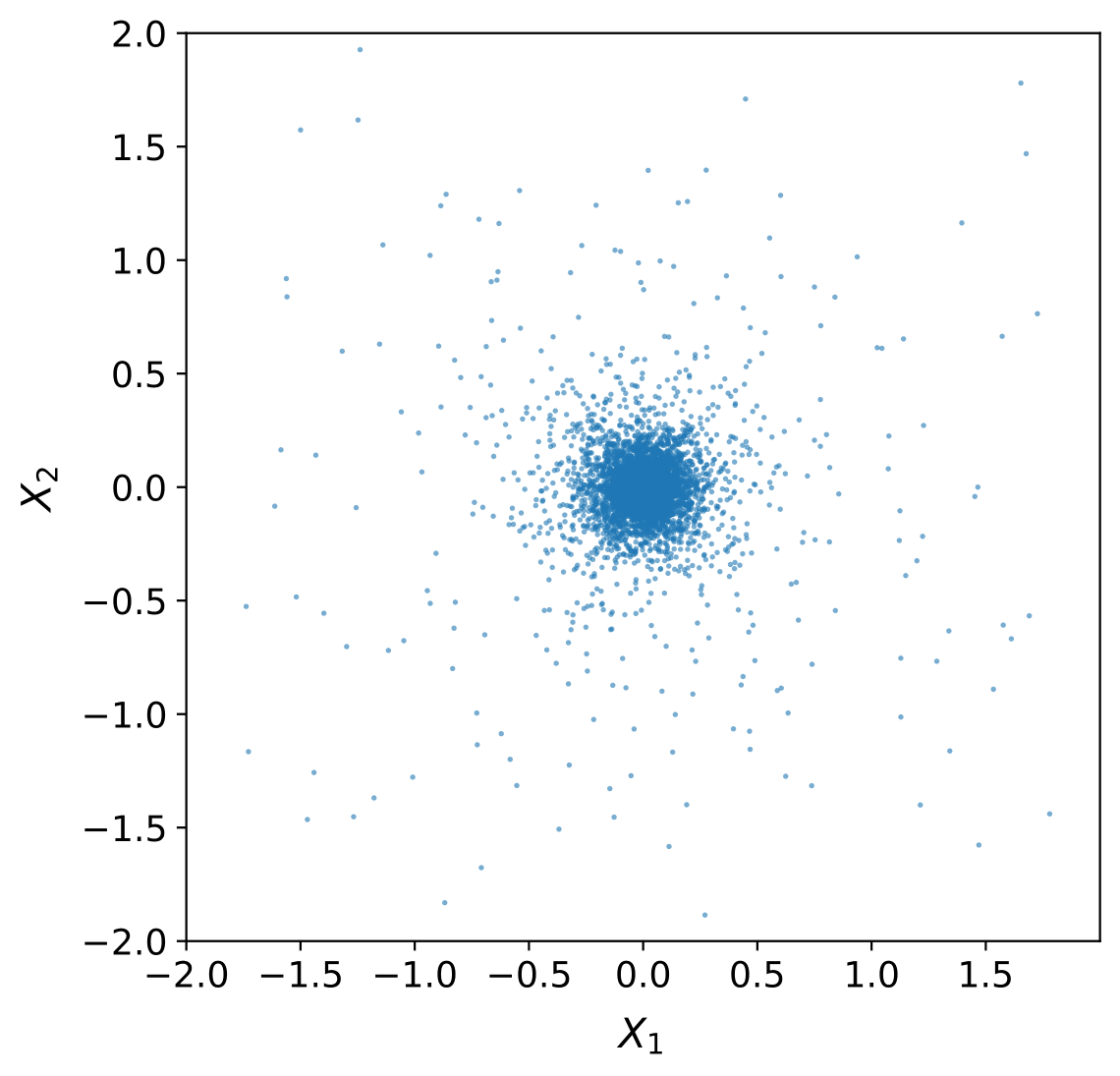}\label{fig:M27_0}}
  
  \caption{Two-dimensional particle distributions at $T = 2.0$ for different initial masses. (a)-(d): Subcritical to near-critical regimes ($M_0 = 22.0$ - $25.0$); (e)-(f): Supercritical regimes ($M_0 = 26.0$ - $27.0$). Parameters $P=2^{22},H=2048,L=8,H_{0}=\left\lceil 8H^{8/13}L^{5/13}\right\rceil,\tau=5\times 10^{-4}$.}
  \label{fig:2d_elliptic}
\end{figure}

Figure~\ref{fig:2d_elliptic} demonstrates the particle distributions at time \( T = 2.0 \) for different initial masses \( M_0 = 22.0, 23.0, \ldots, 27.0 \), exhibiting quantitative agreement with theoretical predictions: subcritical masses (\( M_0 < 8\pi \)) disperse radially, near-critical masses (\( M_0 \approx 8\pi \)) maintain quasi-stationary profiles, and supercritical masses (\( M_0 > 8\pi \)) undergo concentration toward the origin. \begin{color}{black}In all scatter plots presented in this section, $P_{\mathrm{plot}}=6000$ particles are randomly sampled for visualization.\end{color}

\begin{figure}[htbp]
    \centering
    \begin{subfigure}[t]{0.42\textwidth}
        \centering
        \includegraphics[width=\textwidth]{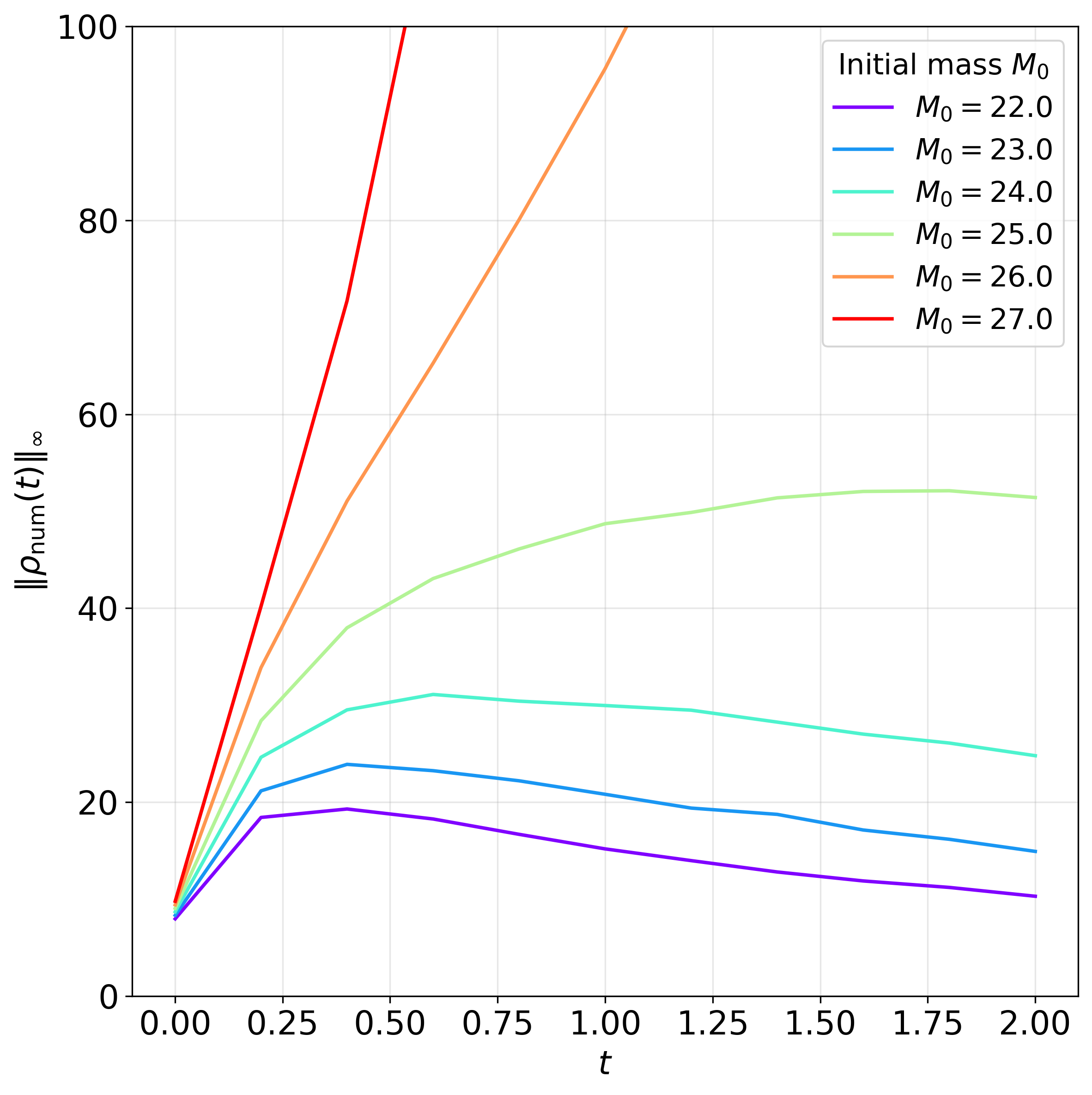}
        \caption{Evolution of $\|\rho_{\mathrm{num}}\|_{\infty}$.}
        \label{fig:2d_elliptic2}
    \end{subfigure}
    \hfill
    \begin{subfigure}[t]{0.42\textwidth}
        \centering
        \includegraphics[width=\textwidth]{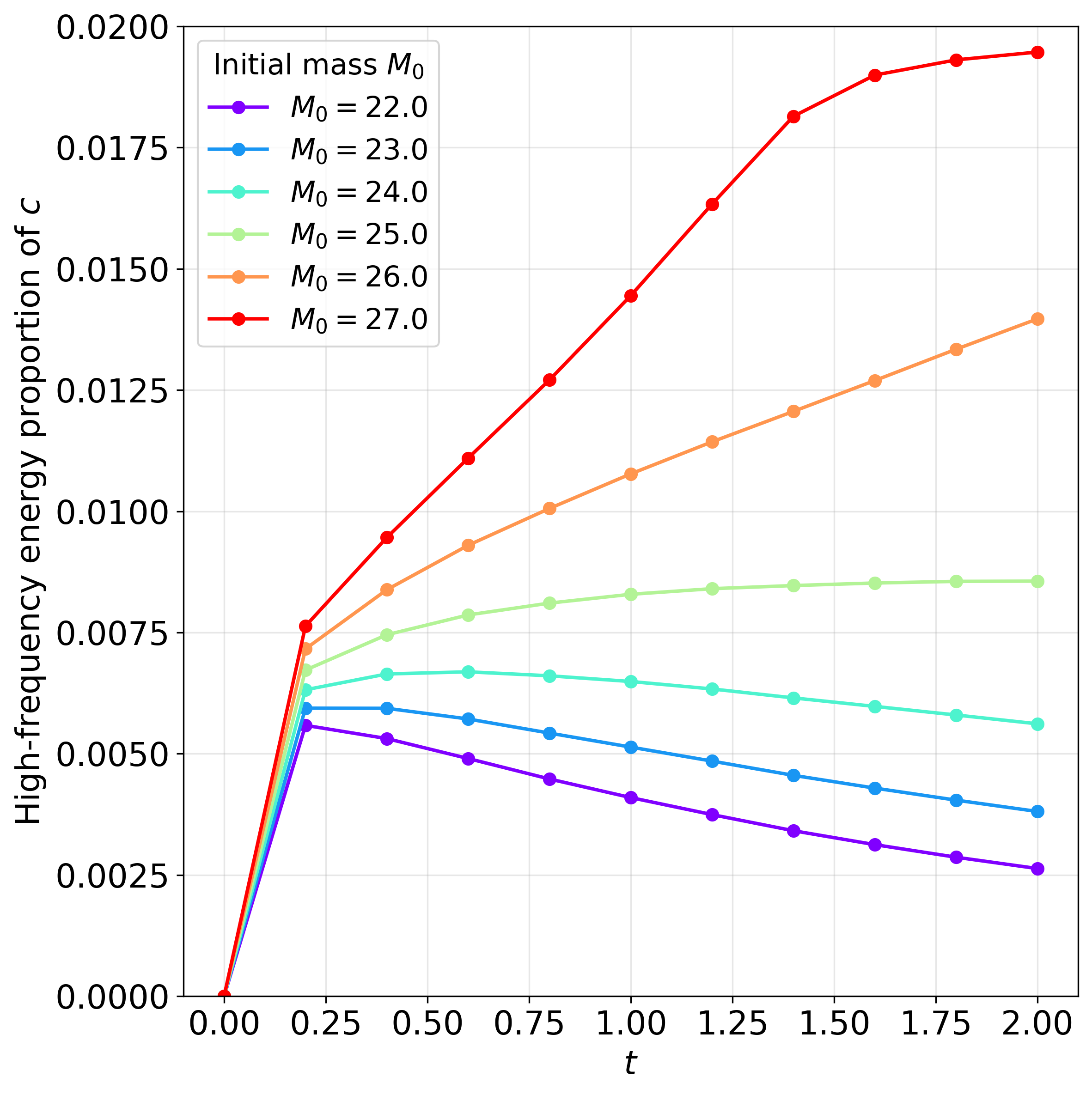}
        \caption{Proportion of high-frequency energy.}
        \label{fig:2d_elliptic3}
    \end{subfigure}

    \caption{Quantitative indicators for the two-dimensional elliptic problem.}
    \label{fig:2d_elliptic_indicators}
\end{figure}

\begin{color}{black}As a quantitative complement, Figure~\ref{fig:2d_elliptic2} shows the evolution of $\|\rho_{\mathrm{num}}(t)\|_{\infty}$ for several initial masses. The horizontal axis represents time, and each curve corresponds to a different value of $M_0$. Here, $\widehat{\rho}_P:=\sum_{p=1}^{P}w_p\delta(\cdot -\mathbf X_p)$ denotes the empirical particle measure, with equal particle weights $w_p=M_0/P$, and $\rho_{\mathrm{num}}:=\widehat{\rho}_P*G_0$ is the associated regularized numerical density. The kernel $G_0$ is an isotropic Gaussian with standard deviation $L/H_0$ in each spatial direction. The quantity $\|\rho_{\mathrm{num}}(t)\|_{\infty}$ is evaluated as the maximum of this regularized density reconstructed on the physical grid.

The original SIPF method \cite{SIPF1} demonstrates that, as blow-up occurs, the proportion of the energy of \(\widehat{c}\) contained in the high-frequency modes ($\|\mathbf{q}\|\geq 4$) increases substantially. Figure~\ref{fig:2d_elliptic3} presents the evolution of this indicator, where the horizontal axis represents time and each curve corresponds to a different value of $M_0$. The observed transition identifies critical masses consistent with those inferred from the particle distributions in Figures~\ref{fig:2d_elliptic} and~\ref{fig:2d_elliptic2}.
\end{color}

The following example demonstrates that, when the periodic domain size $L$ is sufficiently large and the grid spacing $L/H$ is sufficiently small, the influence of the discrepancy between the periodic and unbounded boundary conditions can be progressively reduced, thereby enabling an accurate approximation of the corresponding problem on the whole space. We consider a series of near-critical initial masses \(M_0 = 24.6, 24.8, \ldots, 25.6\), which are close to the theoretical critical threshold \(8\pi \approx 25.1\). For an unbounded domain, the exact second moment at final time \(T=2\) is given by  
\begin{equation}
m_2(t=2) = \int_{\mathbb{R}^2} |\mathbf{x}|^2 \rho(\mathbf{x}, t=2) \, \mathrm{d}\mathbf{x} = \frac{17}{2} - \frac{M_0}{\pi}    
\end{equation}
according to equation \eqref{eq:second_moment}.

To approximate the unbounded domain, we simulate the system on periodic square domains of increasing sizes \(L = 10, 20, 40\), while keeping a fixed spatial mesh size \(L/H = 5/512\). The numerical second moment is computed from the particles as  
\begin{equation}
\widehat{m}_2(t=2) = \frac{1}{P} \sum_{j=1}^P |\widehat{\mathbf{X}}_j(t=2)|^2, 
\end{equation}
where \(P = 2^{20}\) particles are used.

\begin{figure}[htbp]
    \centering \includegraphics[width=0.42\linewidth]{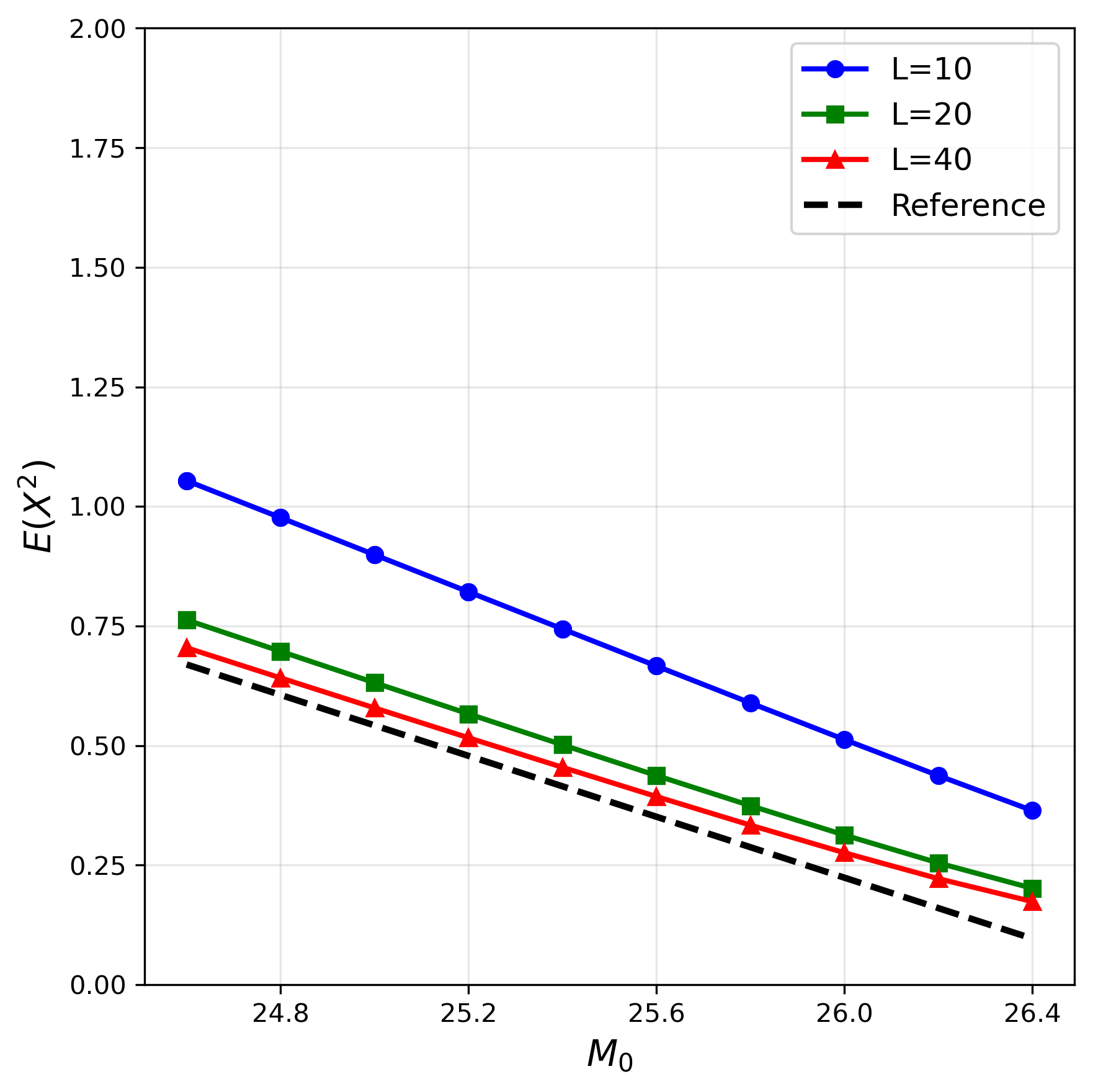}
    \caption{Convergence of the numerical second moment to the theoretical value in unbounded domain. Each polyline corresponds to a different period \(L\), while the dashed line represents the theoretical value. Parameters $P=2^{20},L/H=5/512,H_{0}=\left\lceil 8H^{8/13}L^{5/13}\right\rceil,\tau=5\times 10^{-4}$.}
    \label{fig:second_moment_convergence}
\end{figure}
As shown in Figure~\ref{fig:second_moment_convergence}, the numerical second moment converges to the theoretical prediction as \(L\) grows, confirming that periodic boundary effects diminish and the simulation approaches the unbounded-domain result. This agreement validates the accuracy of our numerical scheme, even near the critical mass where the dynamics are sensitive.

\begin{figure}[htbp]
  \centering
  \subfloat[$M_0 = 45.0$]{\includegraphics[width=0.3\textwidth]{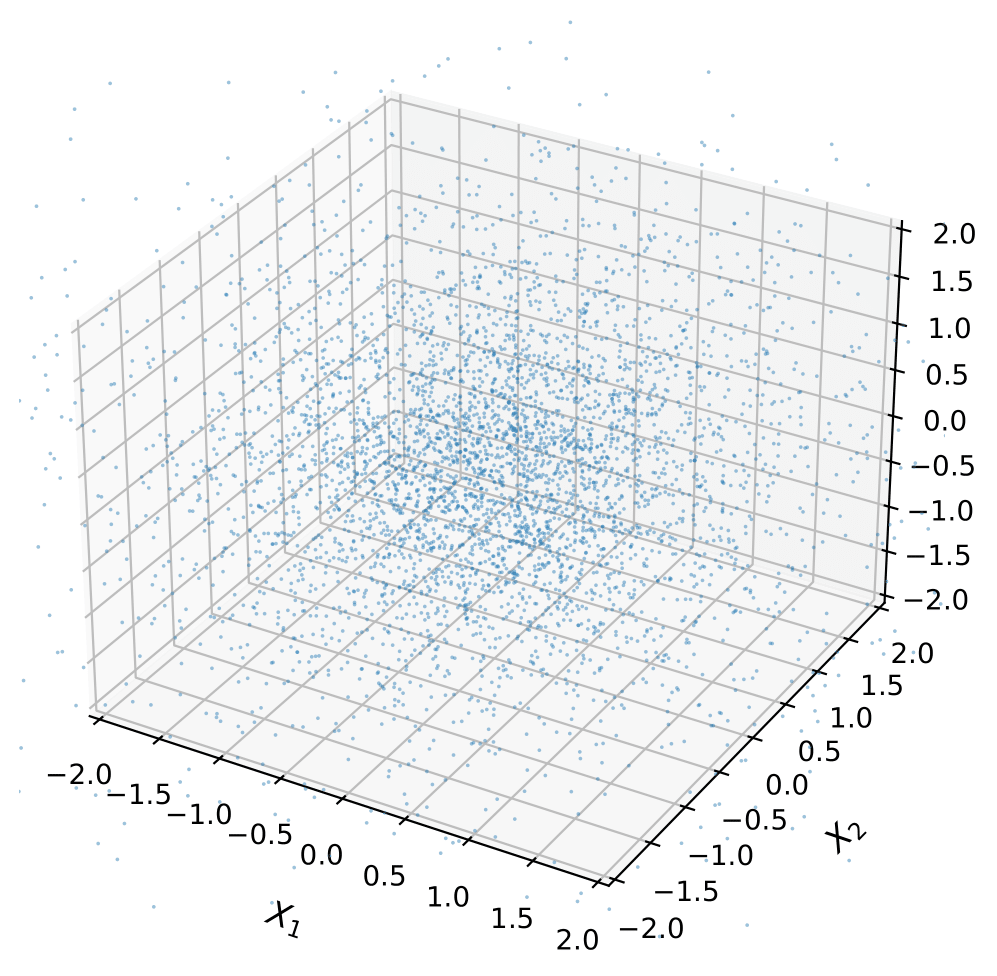}\label{fig:M0_45}}
  \hfill
  \subfloat[$M_0 = 46.0$]{\includegraphics[width=0.3\textwidth]{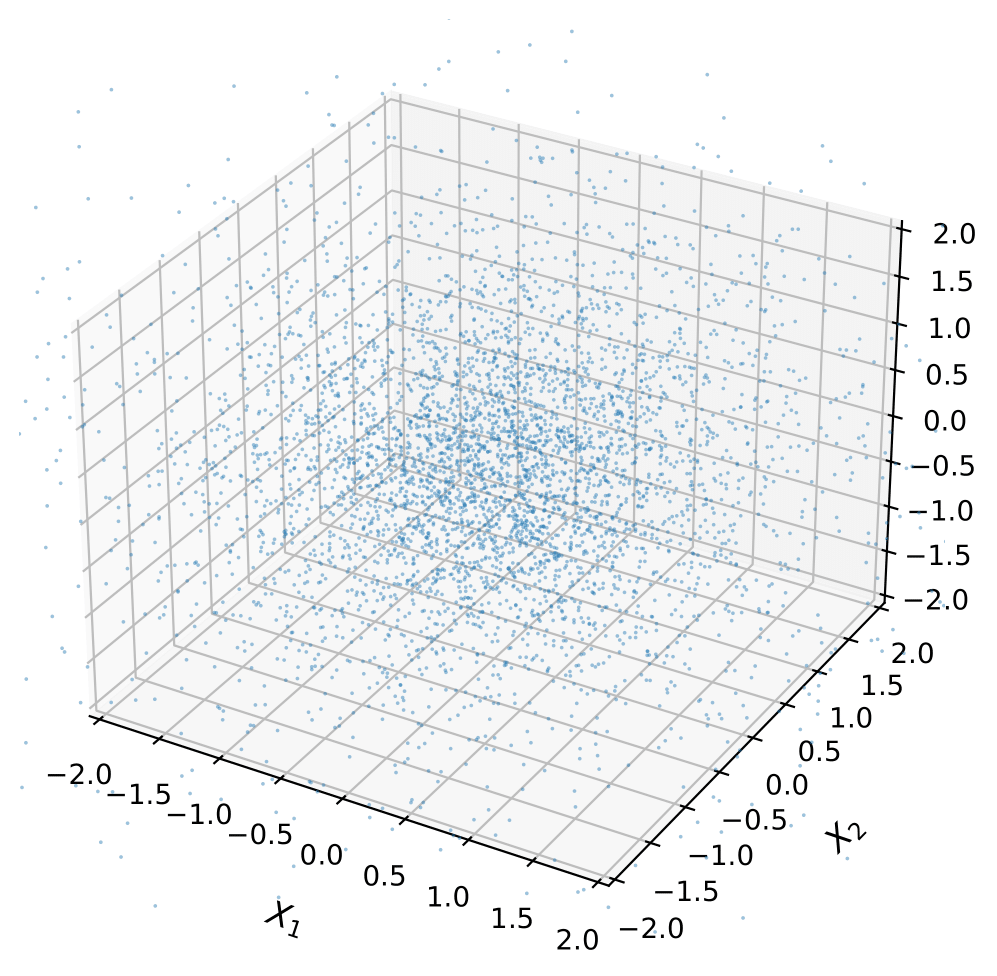}\label{fig:M0_46}}
  \hfill
  \subfloat[$M_0 = 47.0$]{\includegraphics[width=0.3\textwidth]{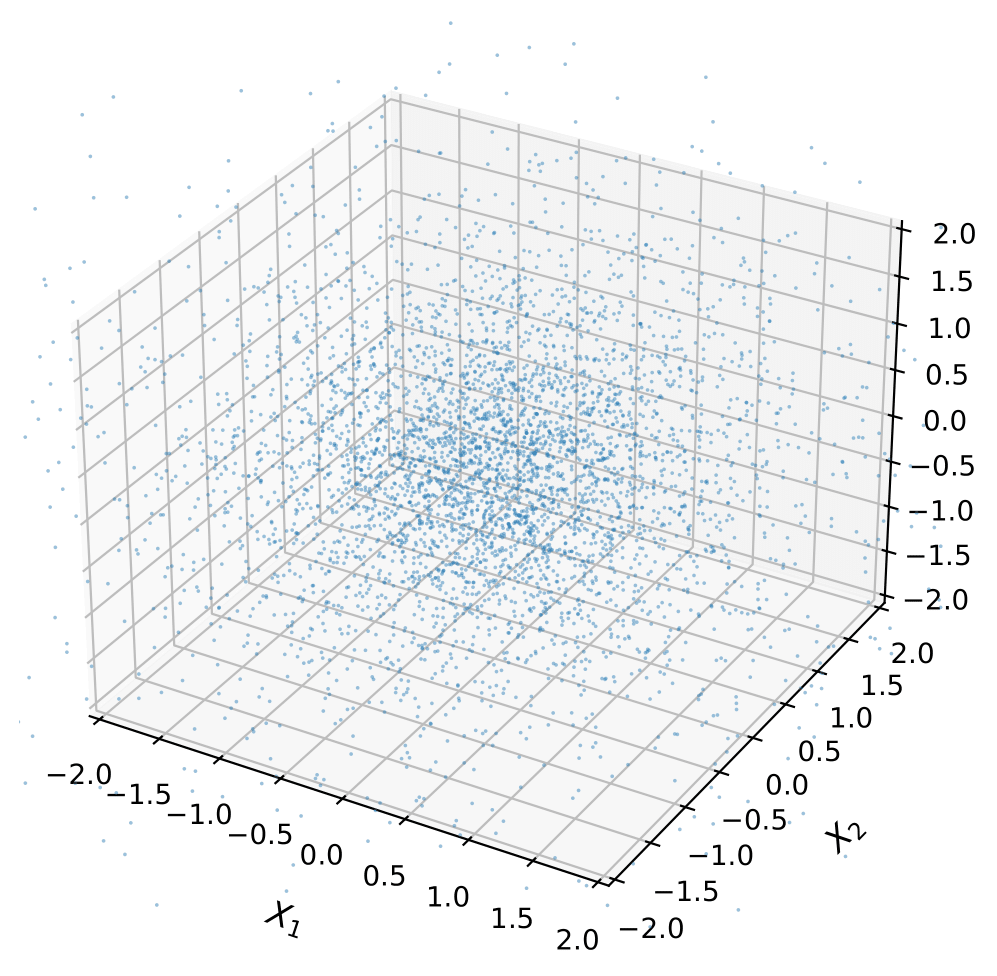}\label{fig:M0_47}}
  
  \vspace{0.5cm}
  \centering
  \subfloat[$M_0 = 48.0$]{\includegraphics[width=0.3\textwidth]{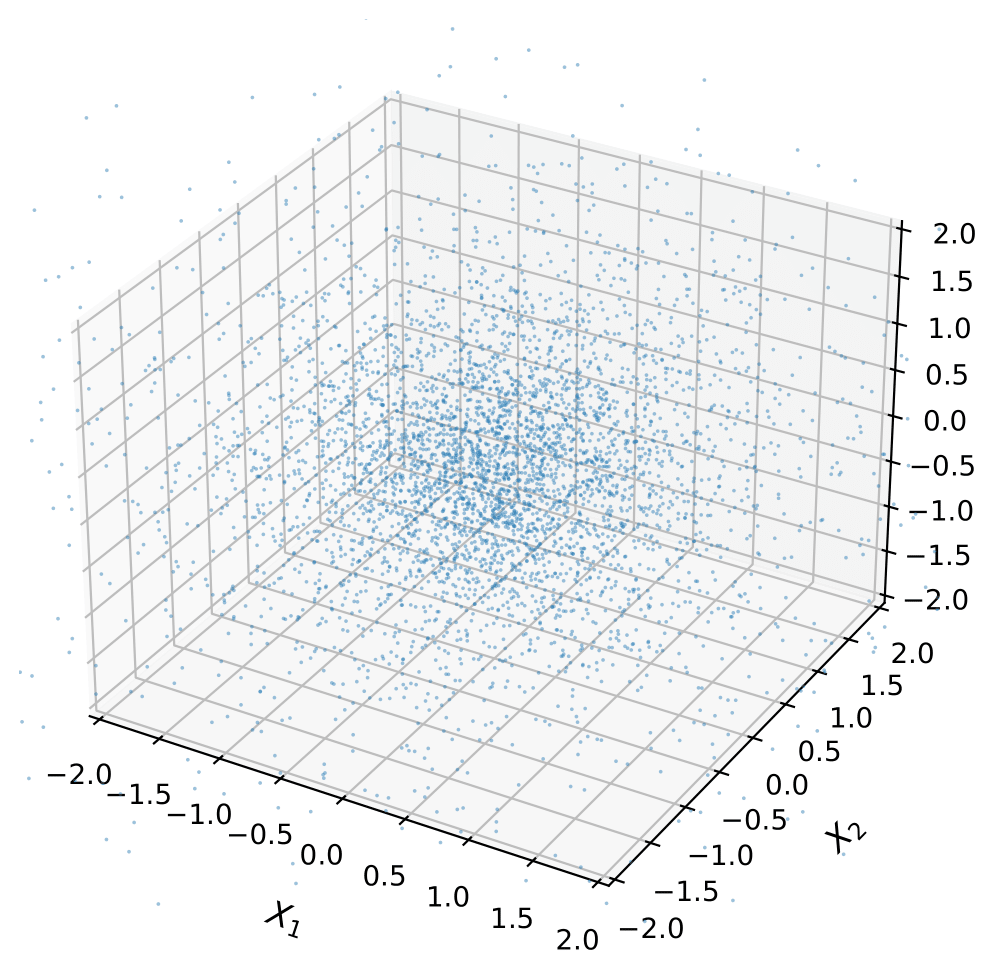}\label{fig:M0_48}}
  \hfill
  \subfloat[$M_0 = 49.0$]{\includegraphics[width=0.3\textwidth]{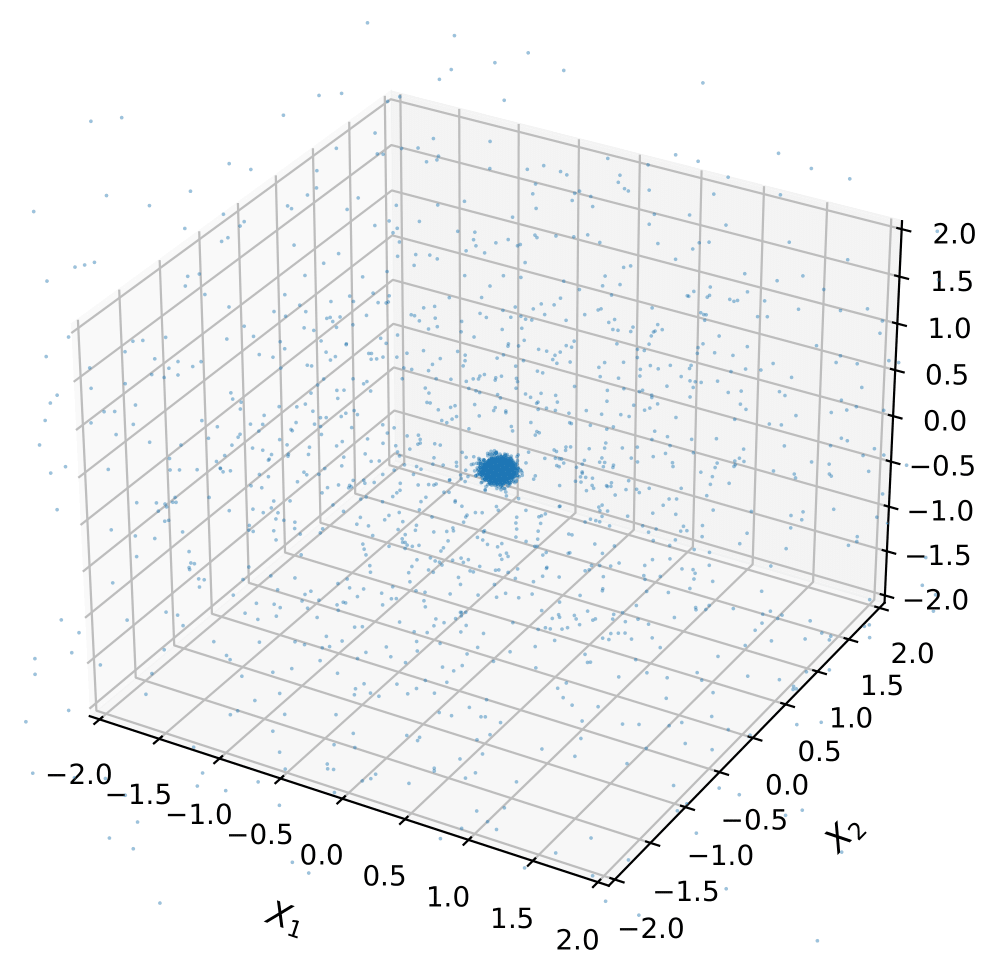}\label{fig:M0_49}}
  \hfill
  \subfloat[$M_0 = 50.0$]{\includegraphics[width=0.3\textwidth]{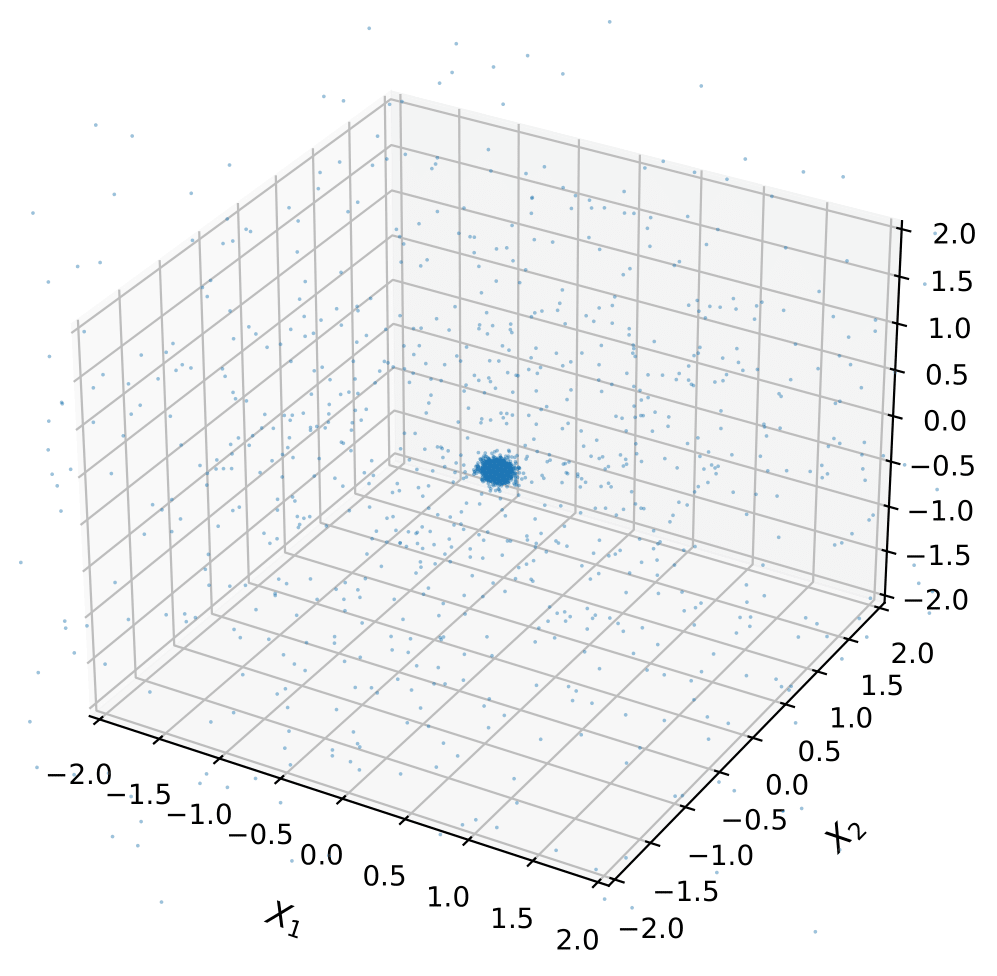}\label{fig:M0_50}}
  
  \caption{Three-dimensional particle distributions at $T = 1.0$ for different initial masses. 
           (a)-(d): Subcritical to near-critical regimes ($M_0 = 45.0$ - $48.0$); 
           (e)-(f): Supercritical regimes ($M_0 = 49.0$ - $50.0$). Parameters $P=2^{20},H=256,L=8,H_{0}=\left\lceil 8H^{8/13}L^{5/13}\right\rceil,\tau=5\times 10^{-4}$.}
  \label{fig:3d_blowup}
\end{figure}

\begin{figure}[htbp]
    \centering
    \begin{subfigure}[t]{0.42\textwidth}
        \centering
        \includegraphics[width=\textwidth]{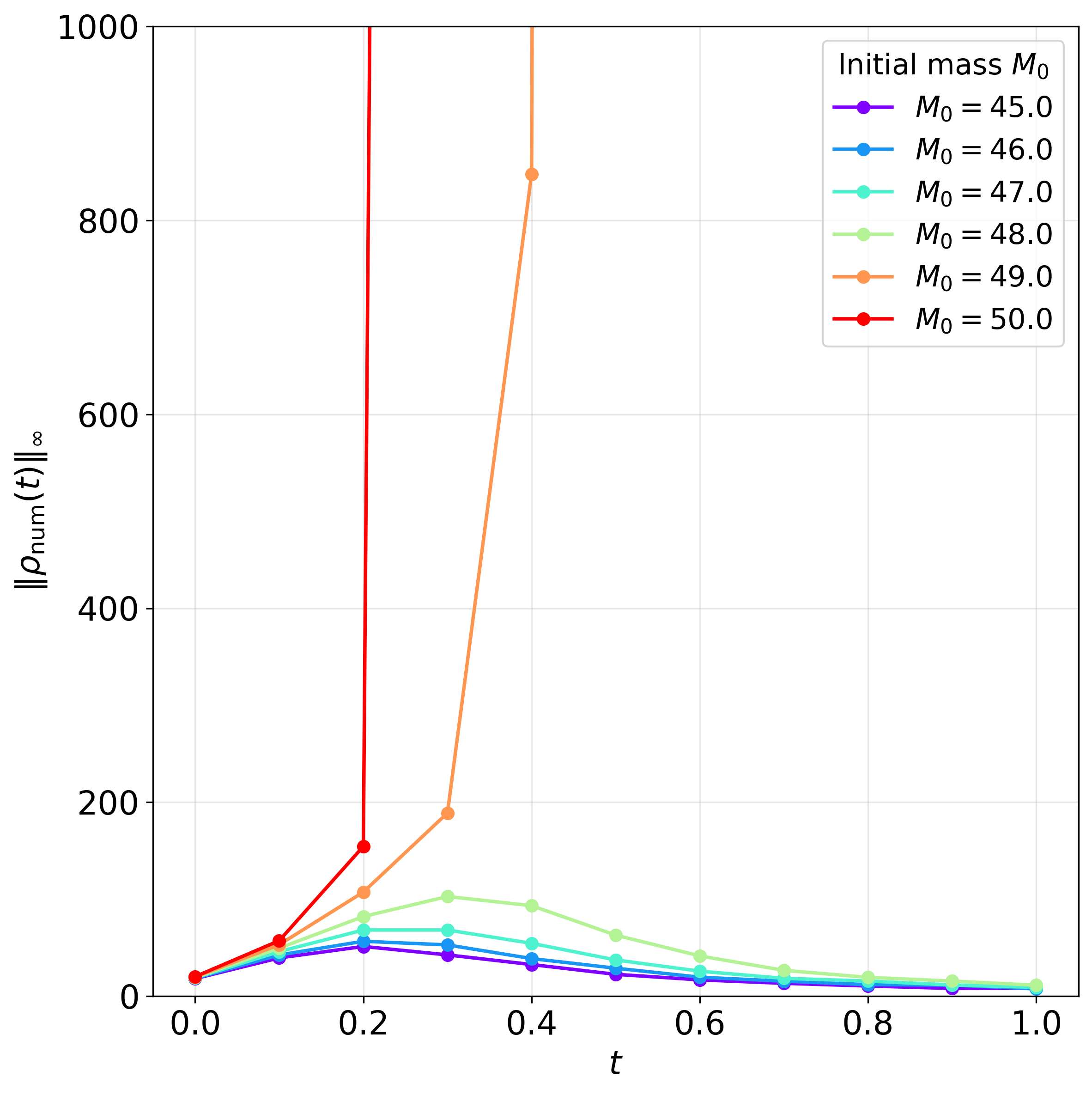}
        \caption{Evolution of $\|\rho_{\mathrm{num}}\|_{\infty}$.}
        \label{fig:3d_blowup2}
    \end{subfigure}
    \hfill
    \begin{subfigure}[t]{0.42\textwidth}
        \centering
        \includegraphics[width=\textwidth]{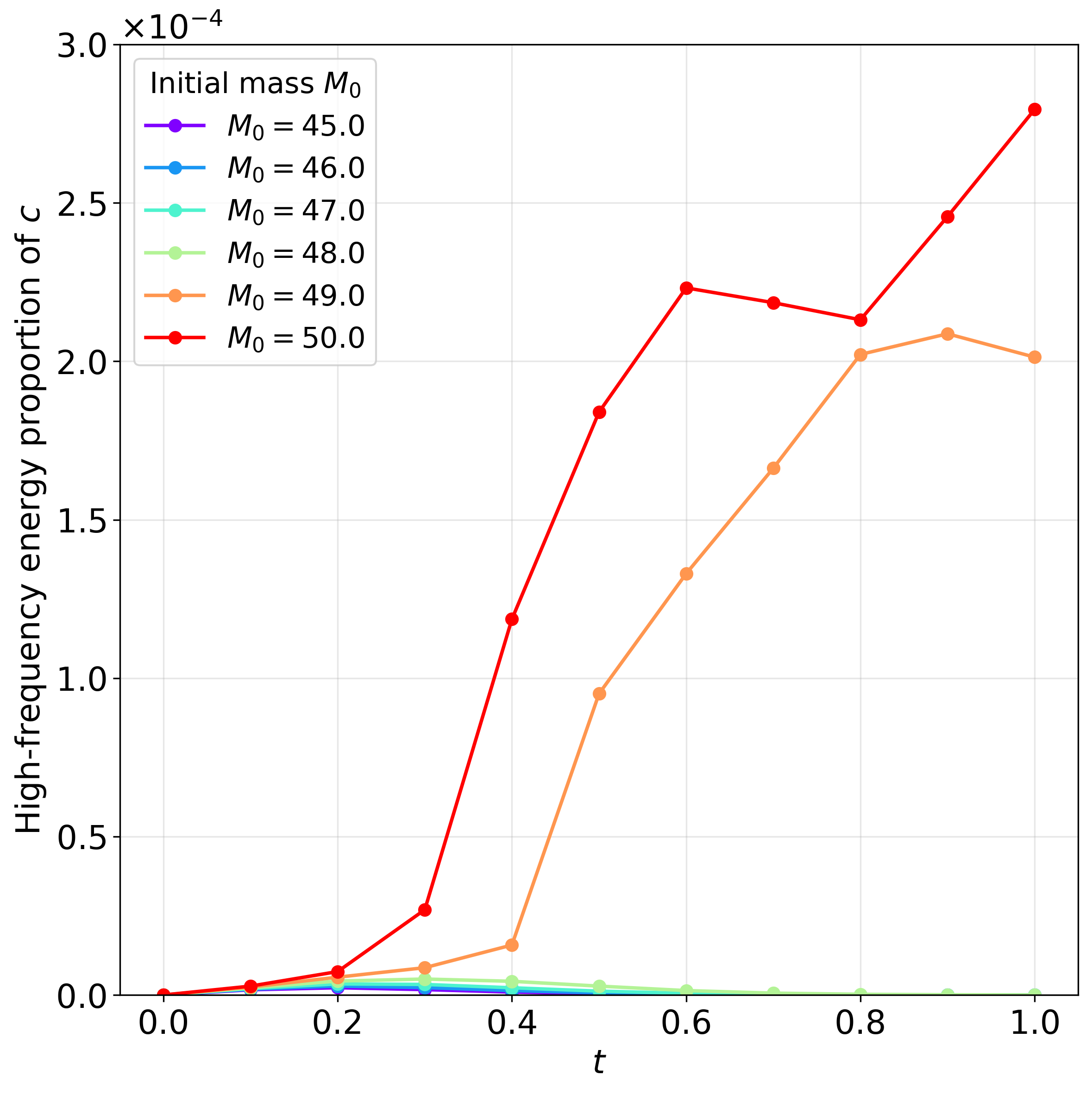}
        \caption{Proportion of high-frequency energy.}
        \label{fig:3d_blowup3}
    \end{subfigure}

    \caption{Quantitative indicators for the three-dimensional parabolic problem.}
    \label{fig:3d_blowup_indicators}
\end{figure}

For the 3D system \eqref{eq:para_system} with parameters $(k, \epsilon, \mu, \chi) = (10^{-1}, 10^{-4}, 1, 1)$, no explicit critical mass threshold is theoretically established. Numerical experiments are conducted with initial masses \(M_0 \in \{45, 46, \ldots, 50\}\) uniformly distributed within the unit ball \(\rho_0(\mathbf{x}) = \frac{3M_0}{4\pi}\mathbf{1}_{\|\mathbf{x}\| \leq 1}\). The spatial domain \(\Omega = [-4,4]^3\) is discretized using \(H = 256\) grid points per axis and \(P = 2^{20}\) particles, coupled with time step \(\tau = 5 \times 10^{-4}\) up to the final time \(T = 1.0\), totaling \( N_T = 2000 \) steps. The particle distributions at the final time corresponding to different initial masses are shown in Figure~\ref{fig:3d_blowup}. \begin{color}{black}The corresponding quantitative indicators, namely the maximum density and the proportion of high-frequency energy, are presented in Figures~\ref{fig:3d_blowup2} and~\ref{fig:3d_blowup3}, respectively. In both figures, the horizontal axis represents time, and each curve corresponds to a different initial mass.\end{color}

Figures \ref{fig:3d_blowup},~\ref{fig:3d_blowup2} and~\ref{fig:3d_blowup3} indicate a critical mass threshold between \(M_0 = 48.0\) and \(M_0 = 49.0\), consistent with the results of finite-volume solutions to the spherically symmetric 1D problem.

\subsection{Blowup dynamics in tetrahedral configuration}  
\label{subsec:tetra}  
\begin{figure}[htbp]
  \centering
  \includegraphics[width=0.42\textwidth]{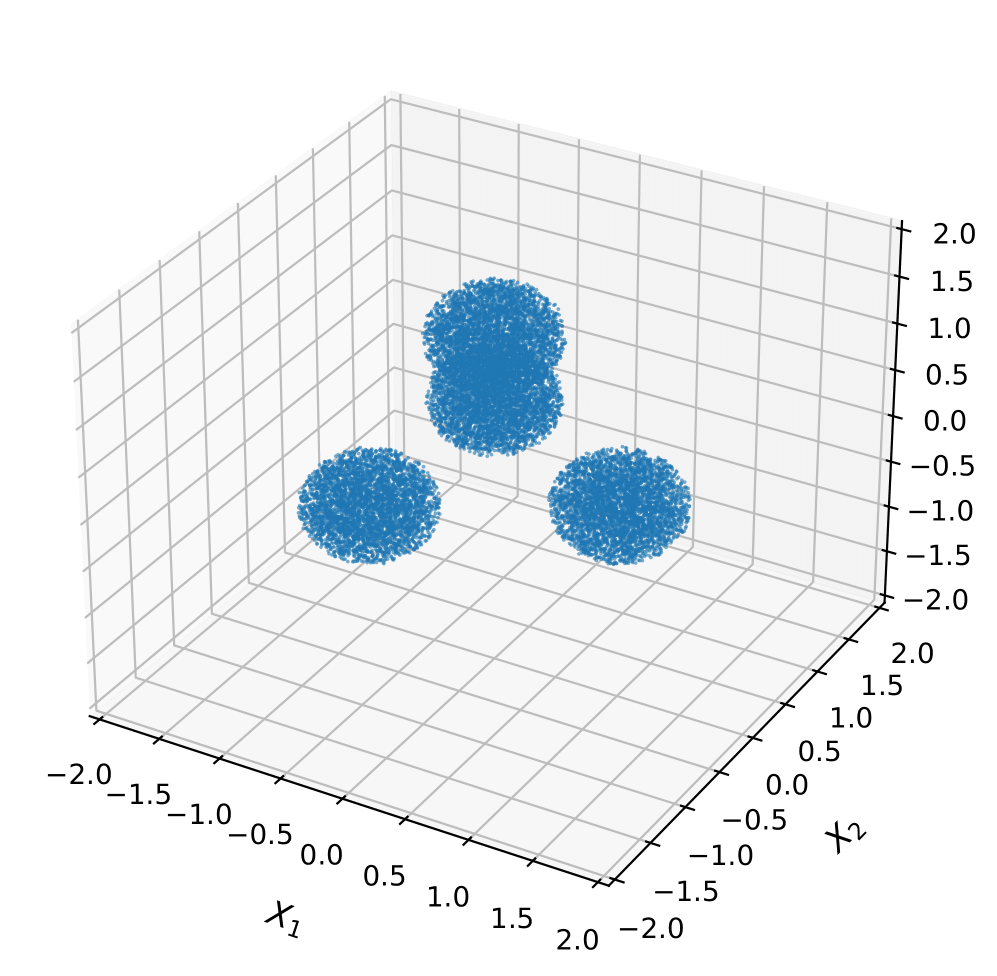}
  \caption{Initial tetrahedral density configuration.}
  \label{fig:tetra_initial}
\end{figure}

We investigate singularity formation using a non-radially symmetric initial configuration: four uniform density balls centered at the vertices of a regular tetrahedron (Fig.~\ref{fig:tetra_initial}). Each ball has mass $\frac{M_0}{4}$, radius $\frac{1}{2}$ and centers at:  
\[
\mathbf{x}_1 = \begin{pmatrix} 1 \\ 0 \\ 0 \end{pmatrix},\;  
\mathbf{x}_2 = \begin{pmatrix} -\frac{1}{2} \\ \frac{\sqrt{3}}{2} \\ 0 \end{pmatrix},\;  
\mathbf{x}_3 = \begin{pmatrix} -\frac{1}{2} \\ -\frac{\sqrt{3}}{2} \\ 0 \end{pmatrix},\;  
\mathbf{x}_4 = \begin{pmatrix} 0 \\ 0 \\ \sqrt{2} \end{pmatrix}
\]  
preserving tetrahedral symmetry with edge length $\sqrt{3}$. The parameters of the system \eqref{eq:para_system} remain $\mu = \chi = 1$, $\epsilon = 10^{-4}$, and $k = 10^{-1}$ as in Section \ref{subsec:rad}. 

\begin{figure}[htbp]
    \centering
    \begin{subfigure}[b]{0.3\textwidth}
        \includegraphics[width=\textwidth]{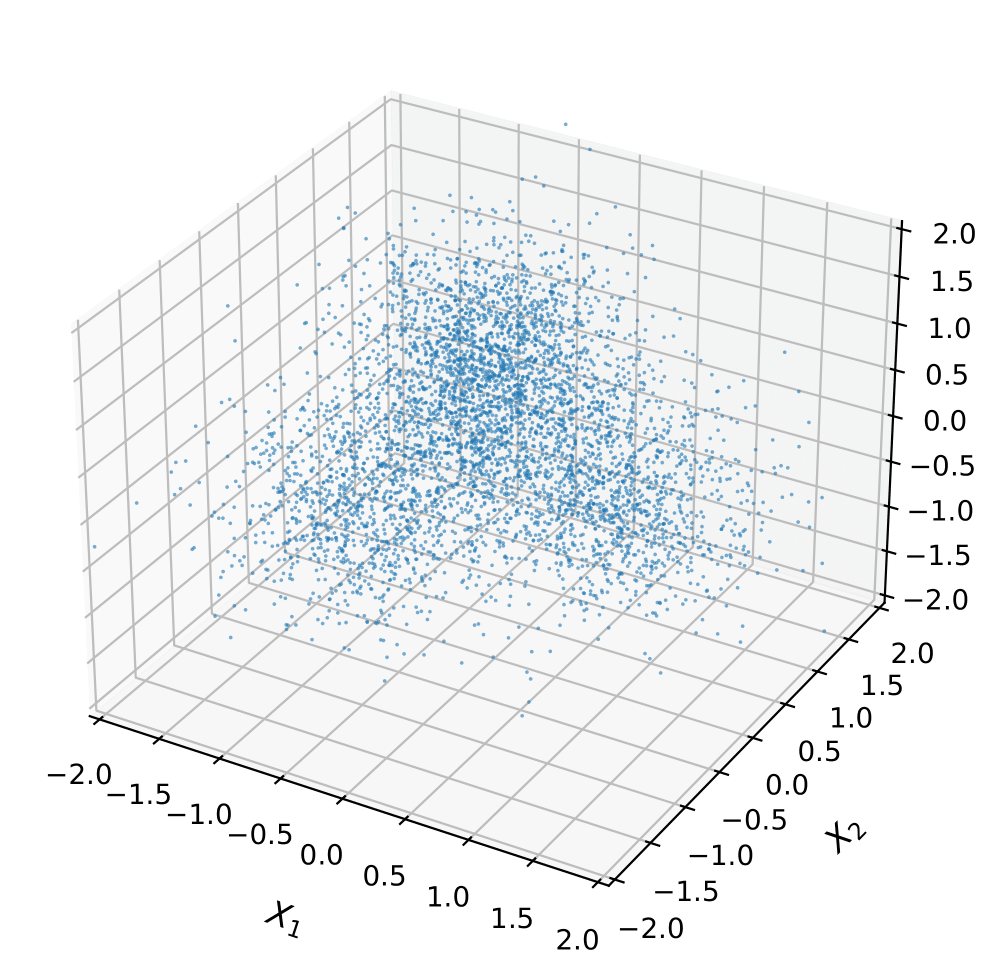}
        \caption*{$M_0 = 20,\,T = 0.1$}
    \end{subfigure}
    \hfill
    \begin{subfigure}[b]{0.3\textwidth}
        \includegraphics[width=\textwidth]{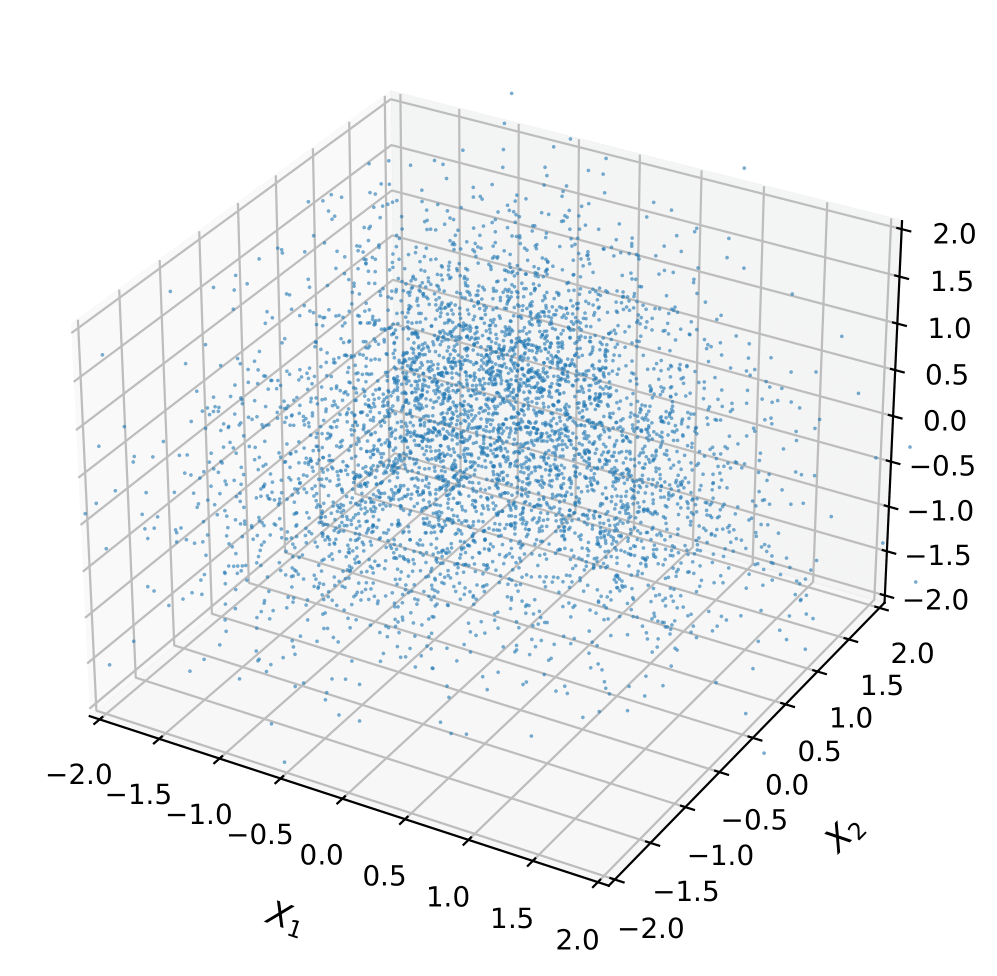}
        \caption*{$M_0 = 20,\,T = 0.2$}
    \end{subfigure}
    \hfill
    \begin{subfigure}[b]{0.3\textwidth}
        \includegraphics[width=\textwidth]{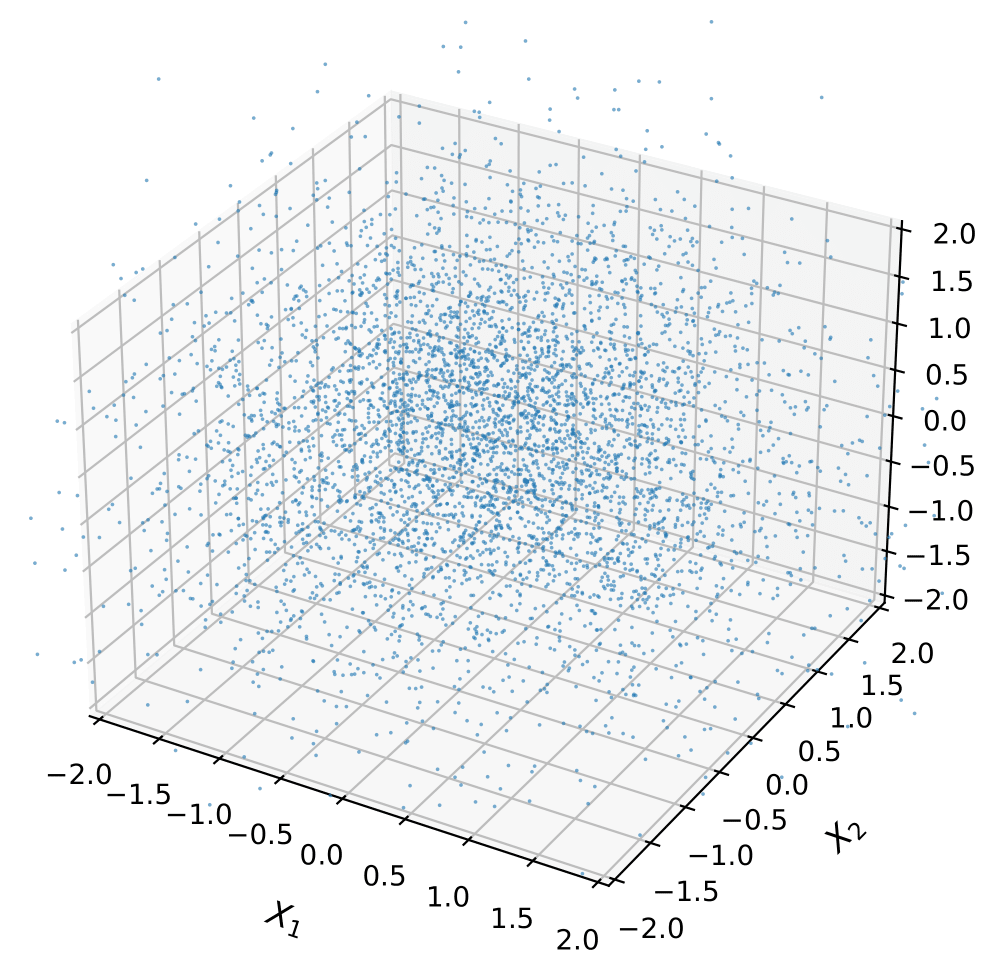}
        \caption*{$M_0 = 20,\,T = 0.4$}
    \end{subfigure}

    \begin{subfigure}[b]{0.3\textwidth}
        \includegraphics[width=\textwidth]{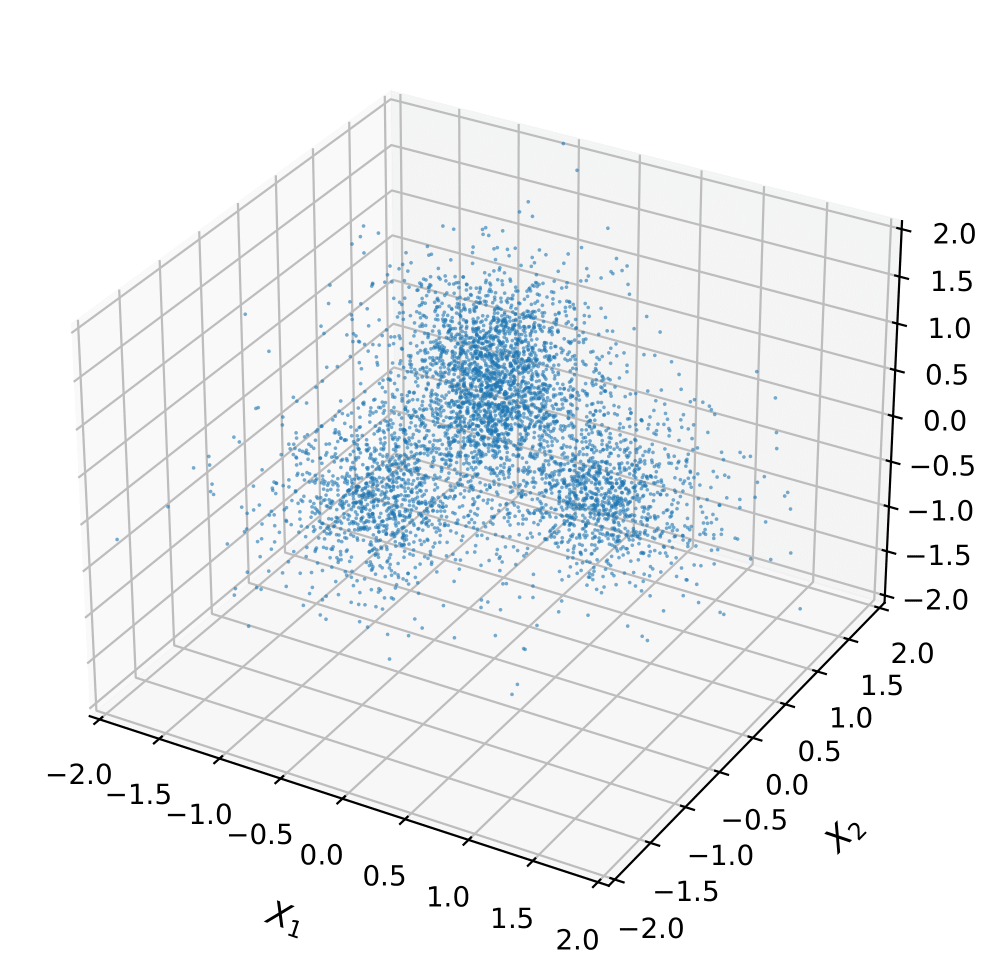}
        \caption*{$M_0 = 80,\,T = 0.1$}
    \end{subfigure}
    \hfill
    \begin{subfigure}[b]{0.3\textwidth}
        \includegraphics[width=\textwidth]{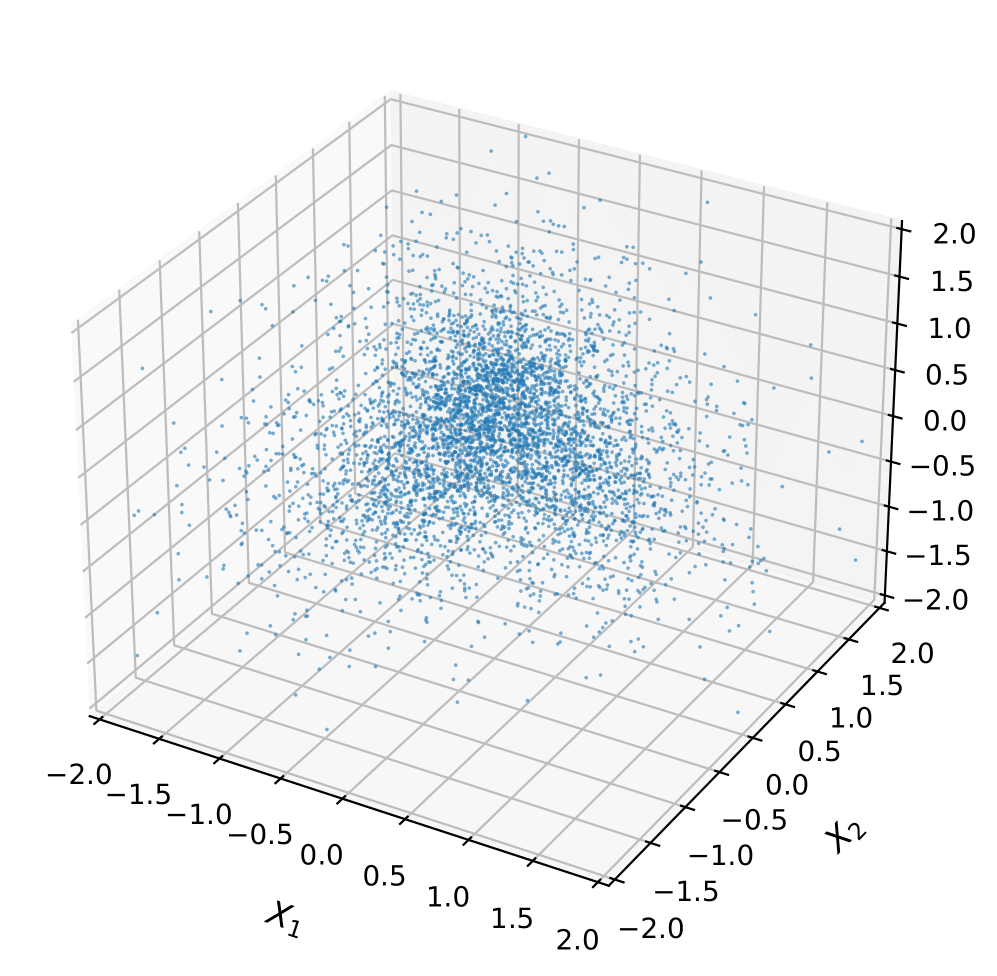}
        \caption*{$M_0 = 80,\,T = 0.2$}
    \end{subfigure}
    \hfill
    \begin{subfigure}[b]{0.3\textwidth}
        \includegraphics[width=\textwidth]{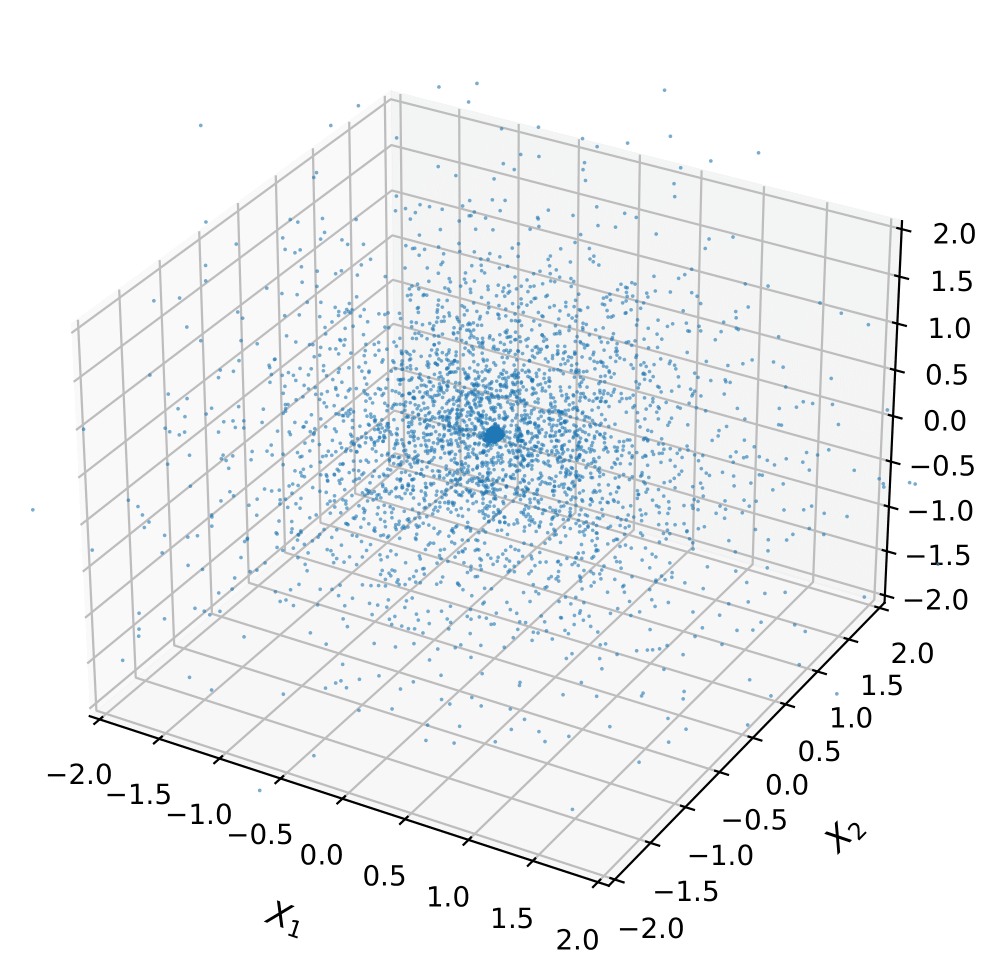}
        \caption*{$M_0 = 80,\,T = 0.4$}
    \end{subfigure}

    \begin{subfigure}[b]{0.3\textwidth}
        \includegraphics[width=\textwidth]{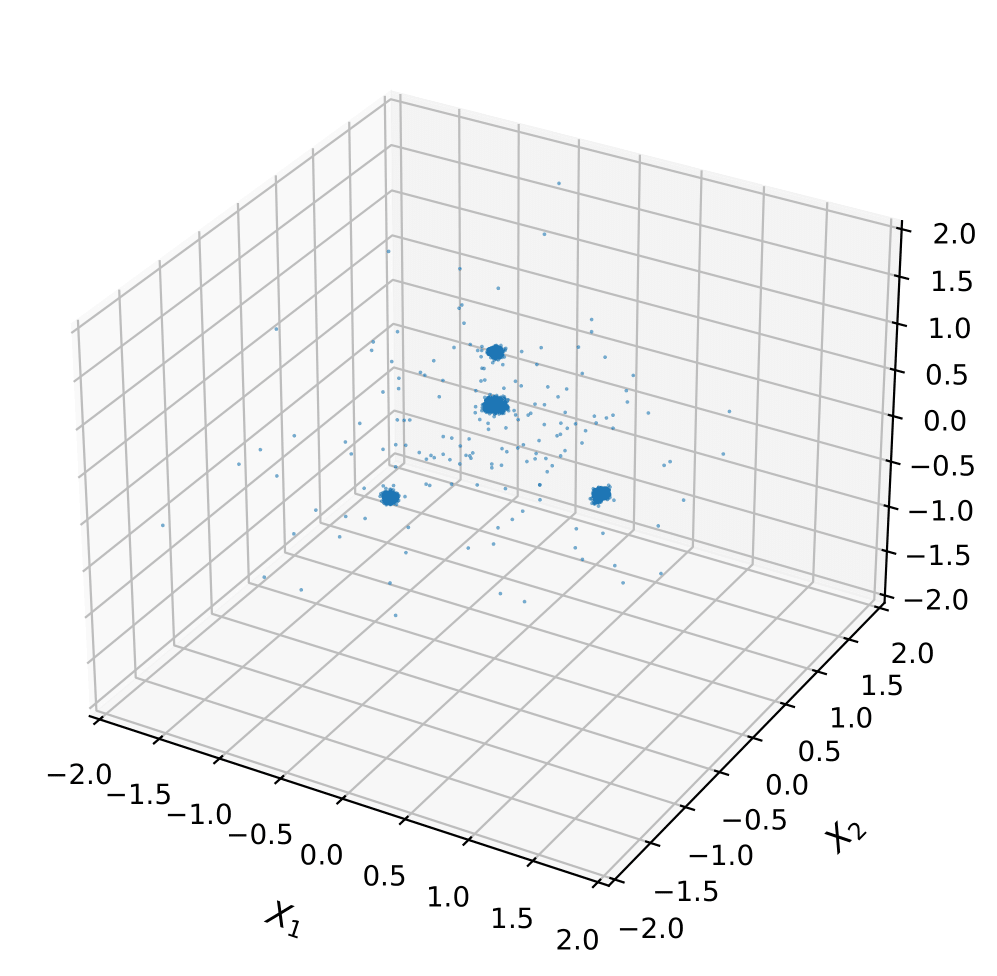}
        \caption*{$M_0 = 160,\,T = 0.1$}
    \end{subfigure}
    \hfill
    \begin{subfigure}[b]{0.3\textwidth}
        \includegraphics[width=\textwidth]{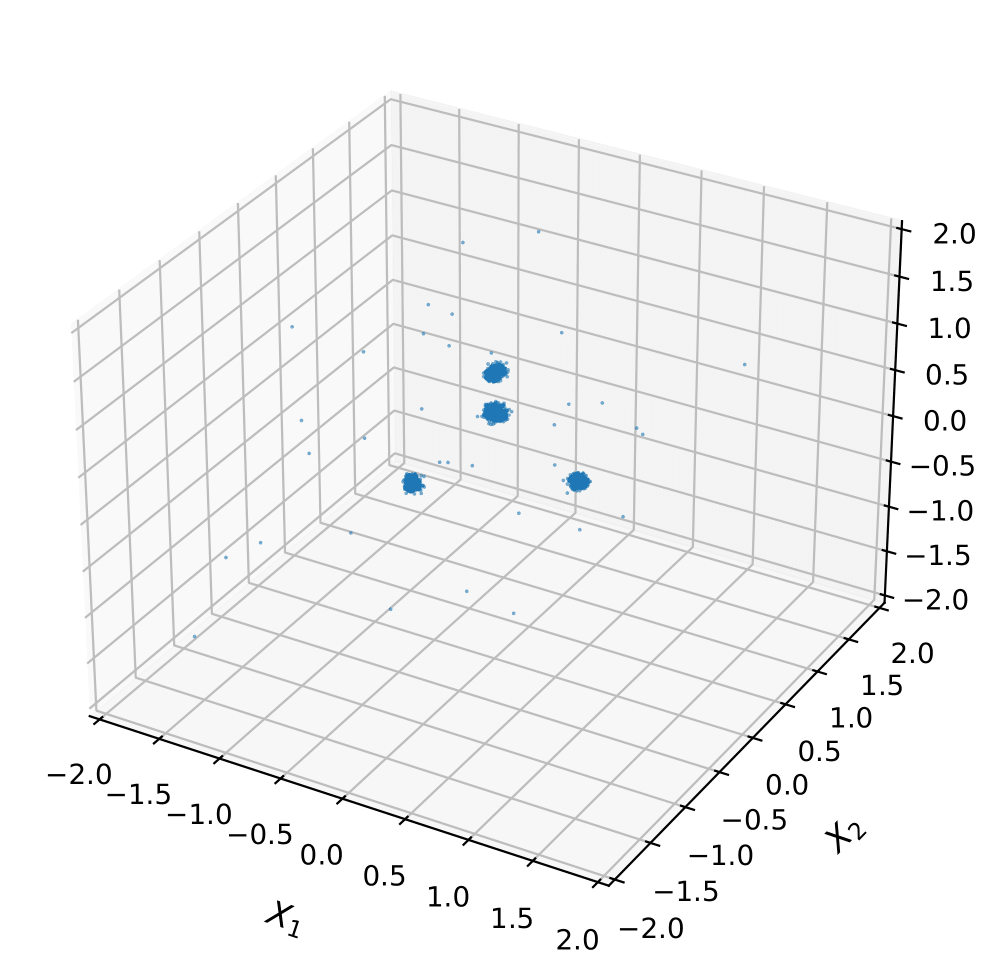}
        \caption*{$M_0 = 160,\,T = 0.2$}
    \end{subfigure}
    \hfill
    \begin{subfigure}[b]{0.3\textwidth}
        \includegraphics[width=\textwidth]{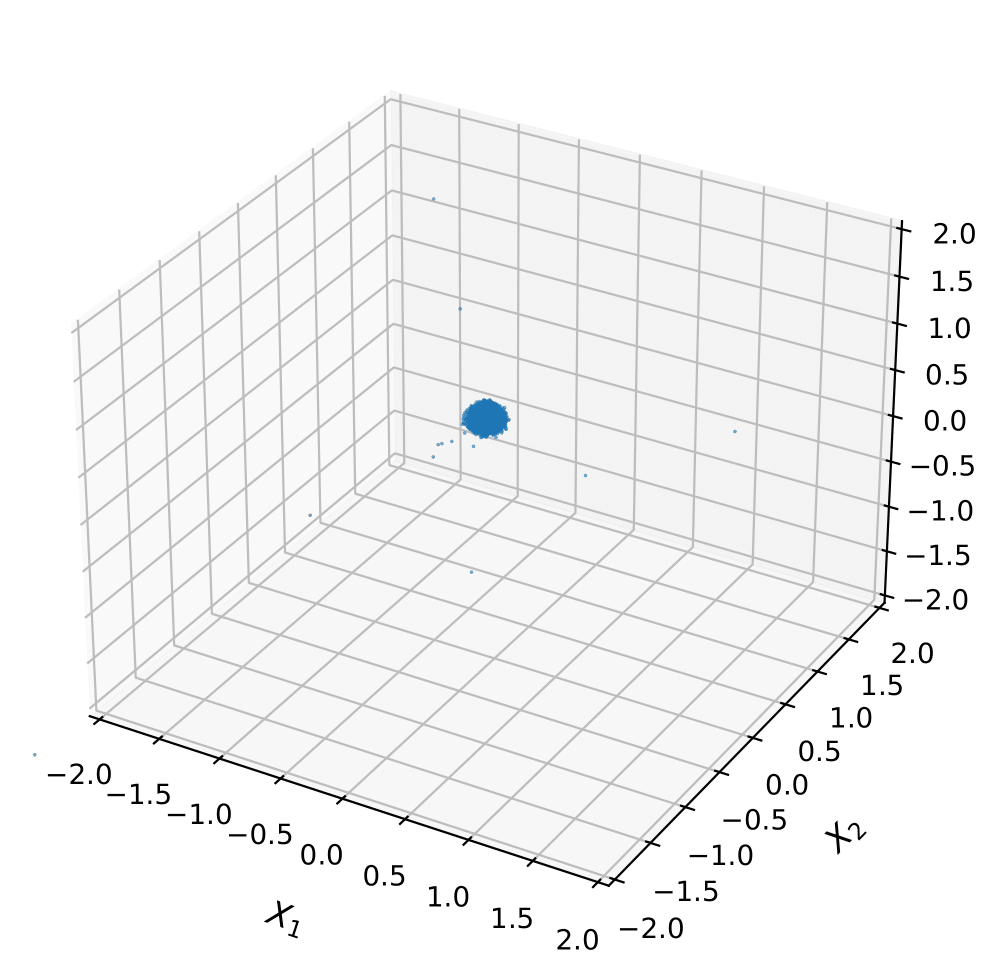}
        \caption*{$M_0 = 160,\,T = 0.4$}
    \end{subfigure}
    
    \caption{
        Blowup patterns for different $M_0$. Parameters $P=2^{20},H=256,L=8,H_{0}=\left\lceil 8H^{8/13}L^{5/13}\right\rceil,\tau= 10^{-4}$.
    }
    \label{fig:tetra_blowup}
\end{figure}

\begin{figure}[htbp]
    \centering
    \begin{subfigure}[t]{0.42\textwidth}
        \centering
        \includegraphics[width=\textwidth]{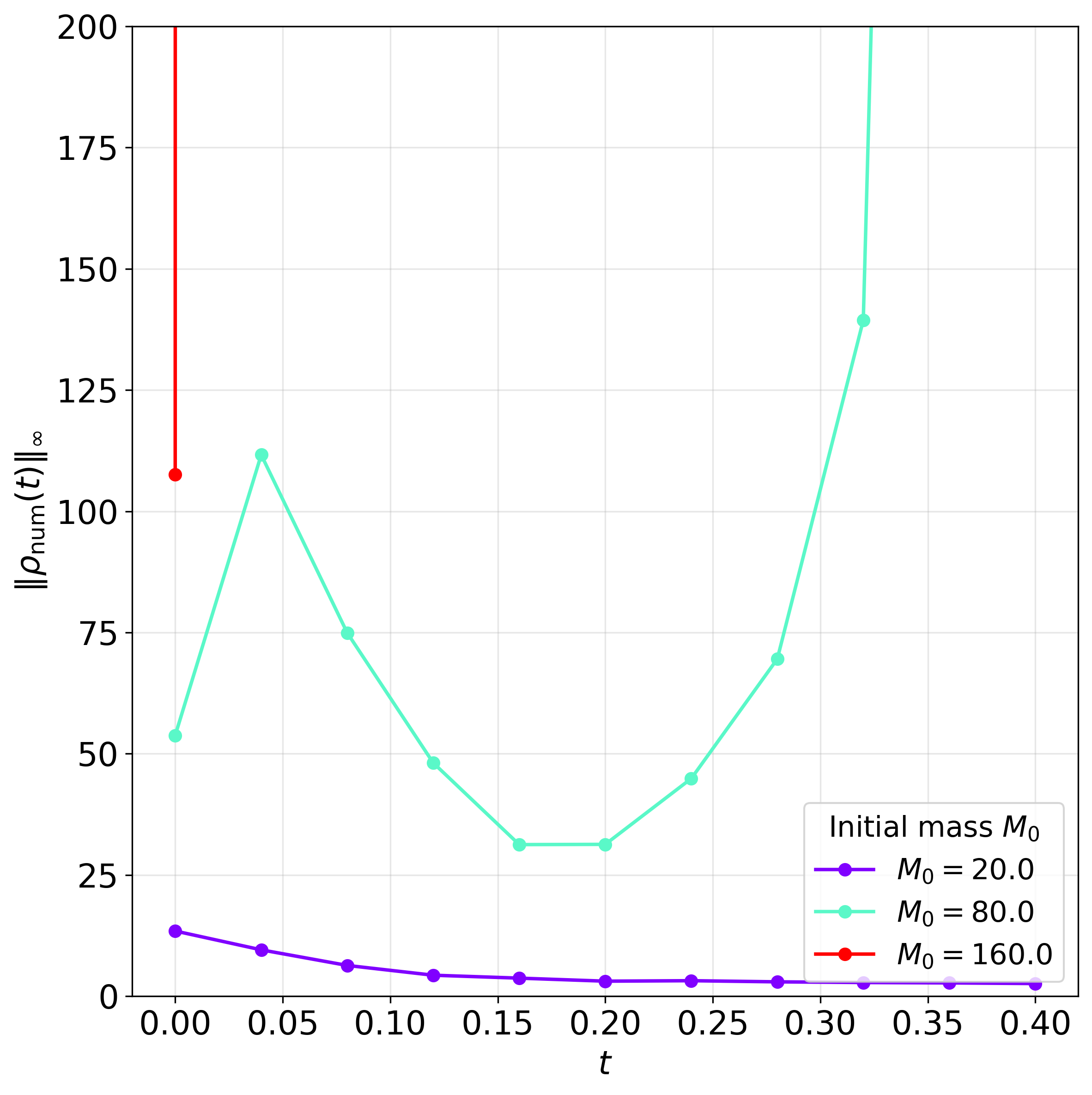}
        \caption{Evolution of $\|\rho_{\mathrm{num}}\|_{\infty}$.}
        \label{fig:tetra_blowup2}
    \end{subfigure}
    \hfill
    \begin{subfigure}[t]{0.42\textwidth}
        \centering
        \includegraphics[width=\textwidth]{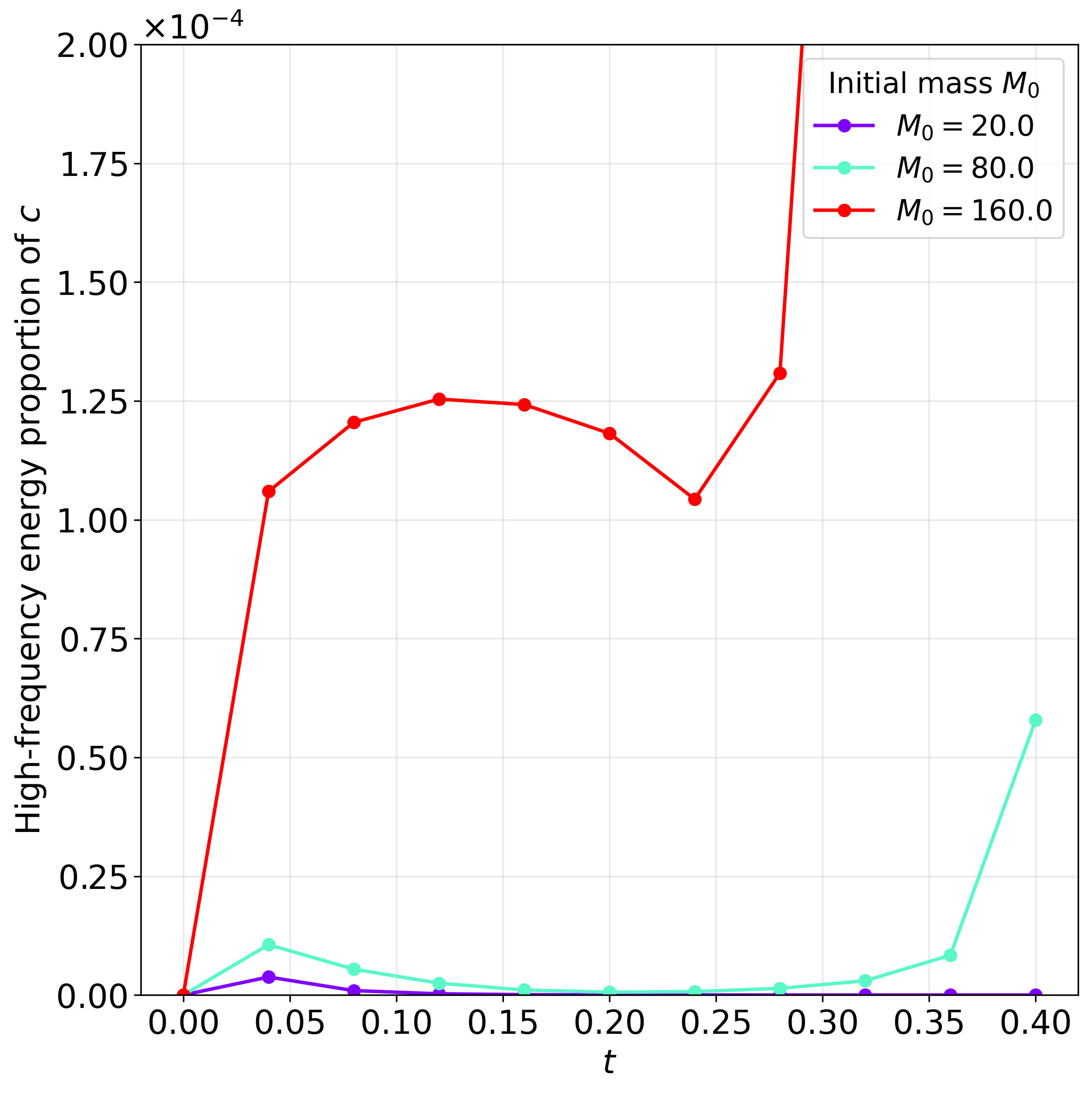}
        \caption{Proportion of high-frequency energy.}
        \label{fig:tetra_blowup3}
    \end{subfigure}

    \caption{Quantitative indicators for the three-dimensional tetrahedral problem.}
    \label{fig:tetra_blowup_indicators}
\end{figure}

We performed simulations for $M_0 \in \{20,28,40,56,80,112,160\}$ with total simulation time $T = 0.4$ and a time step $\tau = 10^{-4}$ (\(N_T = 4000 \)). The grid resolution is fixed at \( H = 256 \), the Gaussian filter $H_{0}=\left\lceil 8H^{8/13}L^{5/13}\right\rceil$. Particle configurations at $T = 0.1$, $0.2$, and $0.4$ are shown in Fig.~\ref{fig:tetra_blowup} for three regimes: $M_0 = 20$ (non-blowup), 80 (central blowup), and 160 (two-stage blowup), with other values demonstrating analogous behavior within their respective categories. \begin{color}{black}The corresponding quantitative indicators, namely the maximum density and the proportion of high-frequency energy, are presented in Figures~\ref{fig:tetra_blowup2} and~\ref{fig:tetra_blowup3}, respectively. In both figures, the horizontal axis represents time, and each curve corresponds to a different initial mass.\end{color}

Three distinct dynamical regimes are observed. For small $M_0$, no blowup phenomenon occurs. At intermediate $M_0$, the four clusters merge into a single aggregate prior to a single unified blowup at the center. For large $M_0$, the system exhibits a two-stage process: simultaneous blowups initially develop at the original tetrahedral positions, followed by coalescence into a central singularity.

\subsection{Collapsing ring blowup dynamics}\label{subsec:ring_blowup}
In the last experiment, we investigate a special blow-up pattern where the mass concentrates on an evolving ring \(x_1^2 + x_2^2 = r(t)^2\), with the radius \(r(t)\) decaying as \(t \to T\).

Numerical simulations initialize with a uniform distribution inside the toroidal volume defined by \((1.0 - \sqrt{x_1^2 + x_2^2})^2 + x_3^2 \leq 0.4^2\), with parameters \((k, \epsilon, \mu, \chi) = (10^{-1}, 10^{-4}, 1, 1)\) in system \eqref{eq:para_system}, an initial mass \(M_0 = 180\), and a final computational time \(T = 6 \times 10^{-2}\). The spatial domain \(\Omega = [-4,4]^3\) is discretized with \(H = 256\) grid points per axis, using \(P = 2^{20}\) particles, the Gaussian filter $H_{0}=\left\lceil 8H^{8/13}L^{5/13}\right\rceil$, time step \(\tau = 5 \times 10^{-5}\), and total steps \( N_T = 1200 \).

\begin{figure}[htbp]
  \centering
  \subfloat[\(t = 0.01\)]{\includegraphics[width=0.3\textwidth]{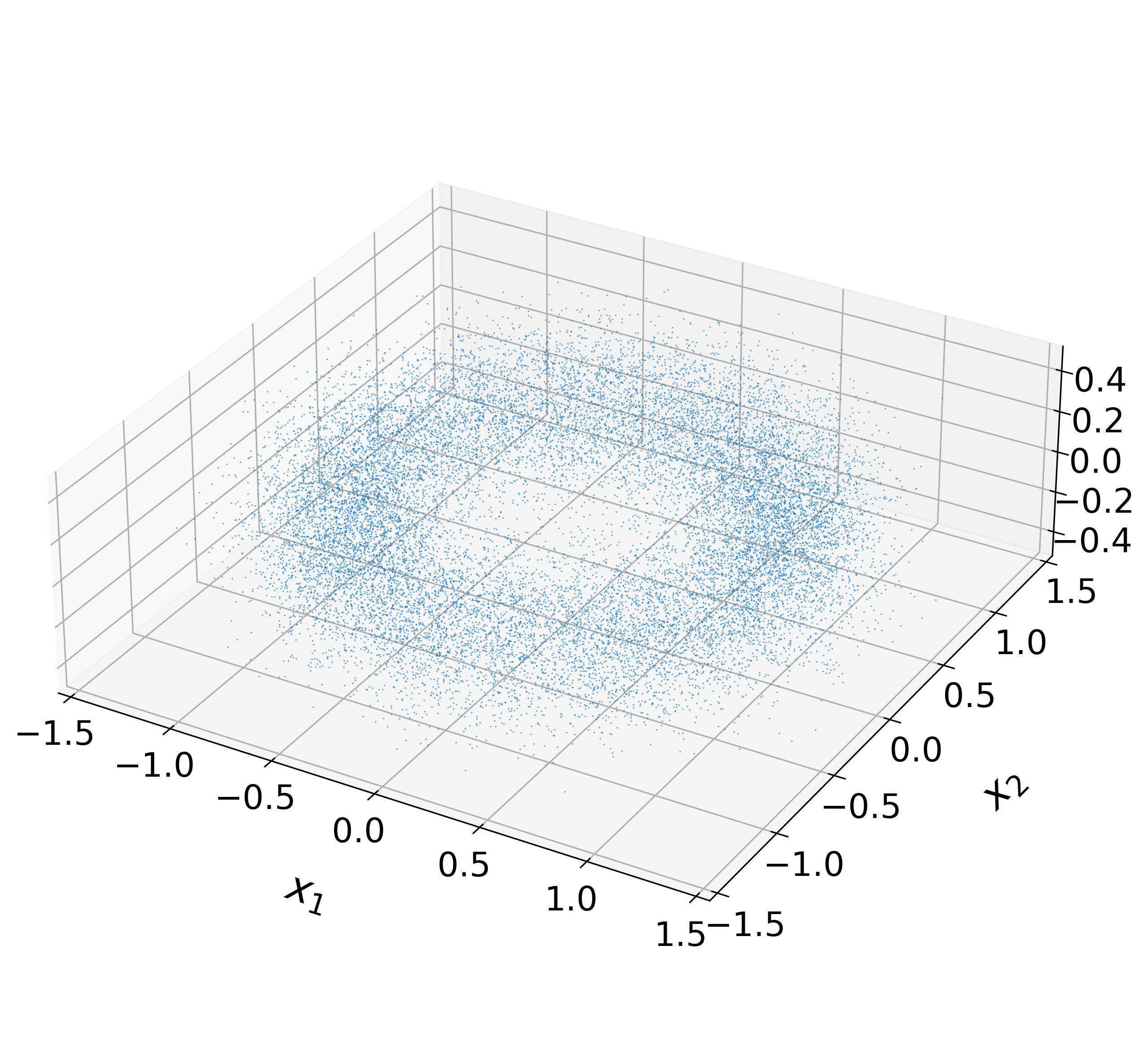}}
  \hfill
  \subfloat[\(t = 0.02\)]{\includegraphics[width=0.3\textwidth]{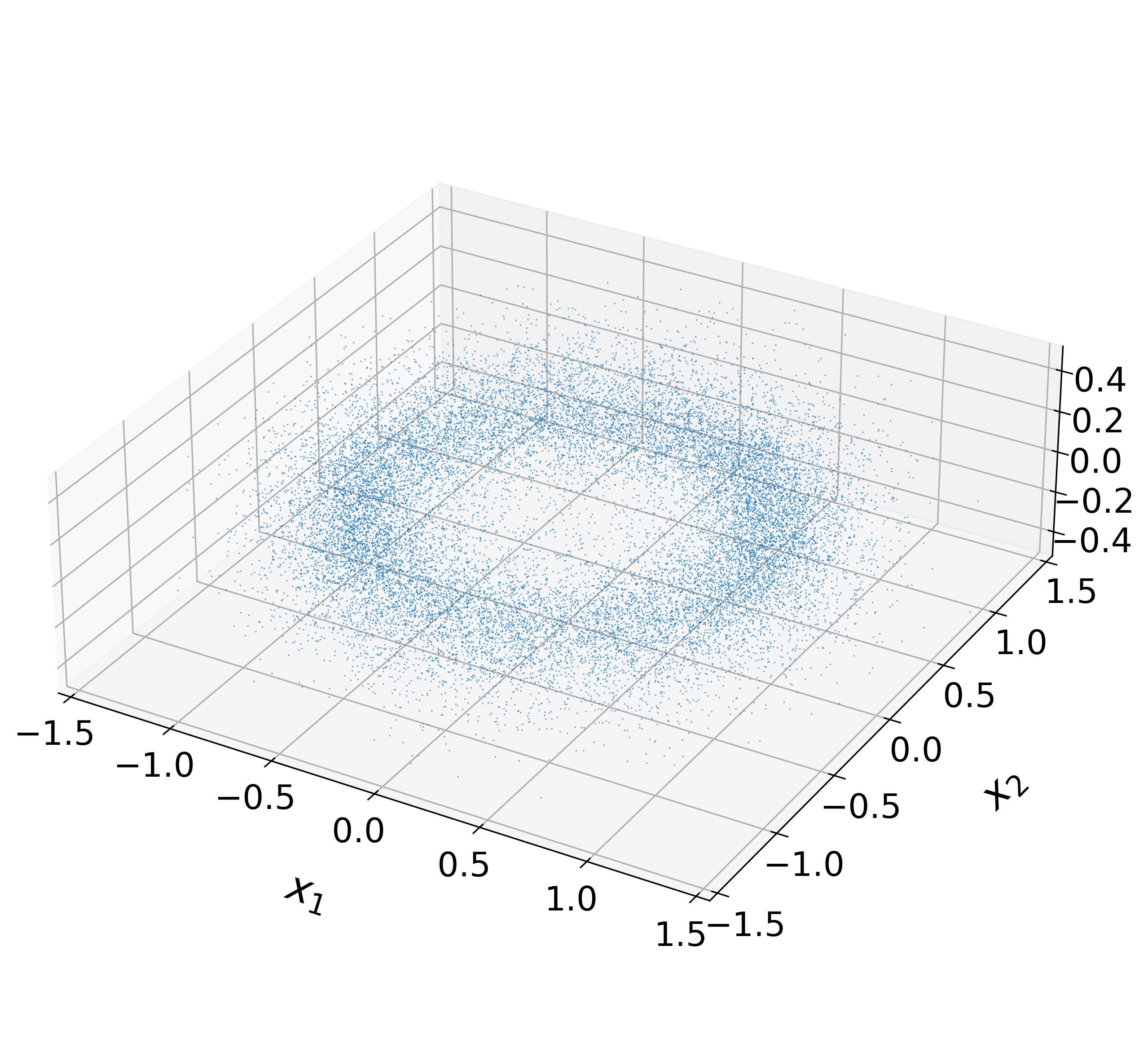}}
  \hfill
  \subfloat[\(t = 0.03\)]{\includegraphics[width=0.3\textwidth]{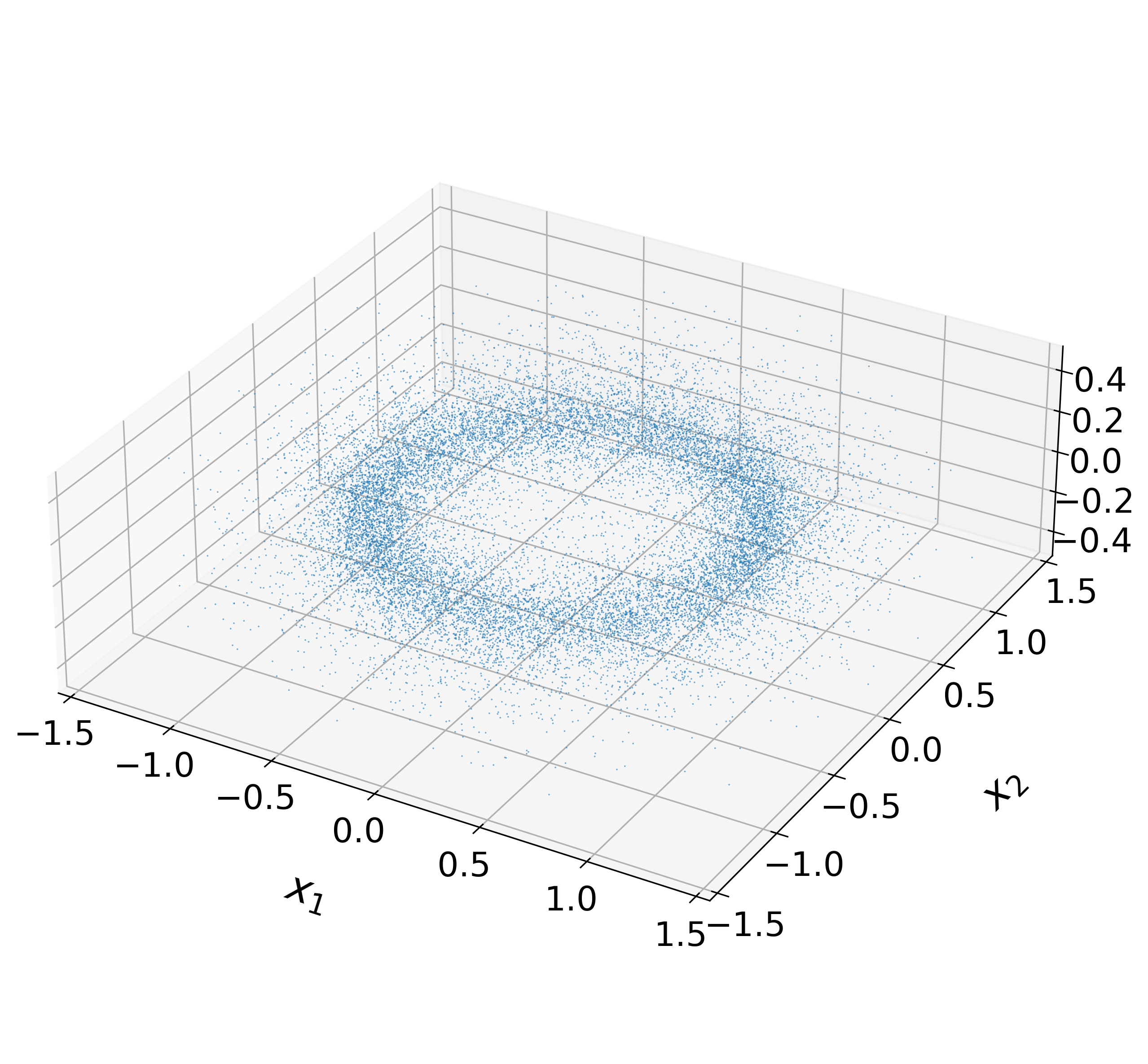}}
  
  \vspace{0.3cm}
  \subfloat[\(t = 0.04\)]{\includegraphics[width=0.3\textwidth]{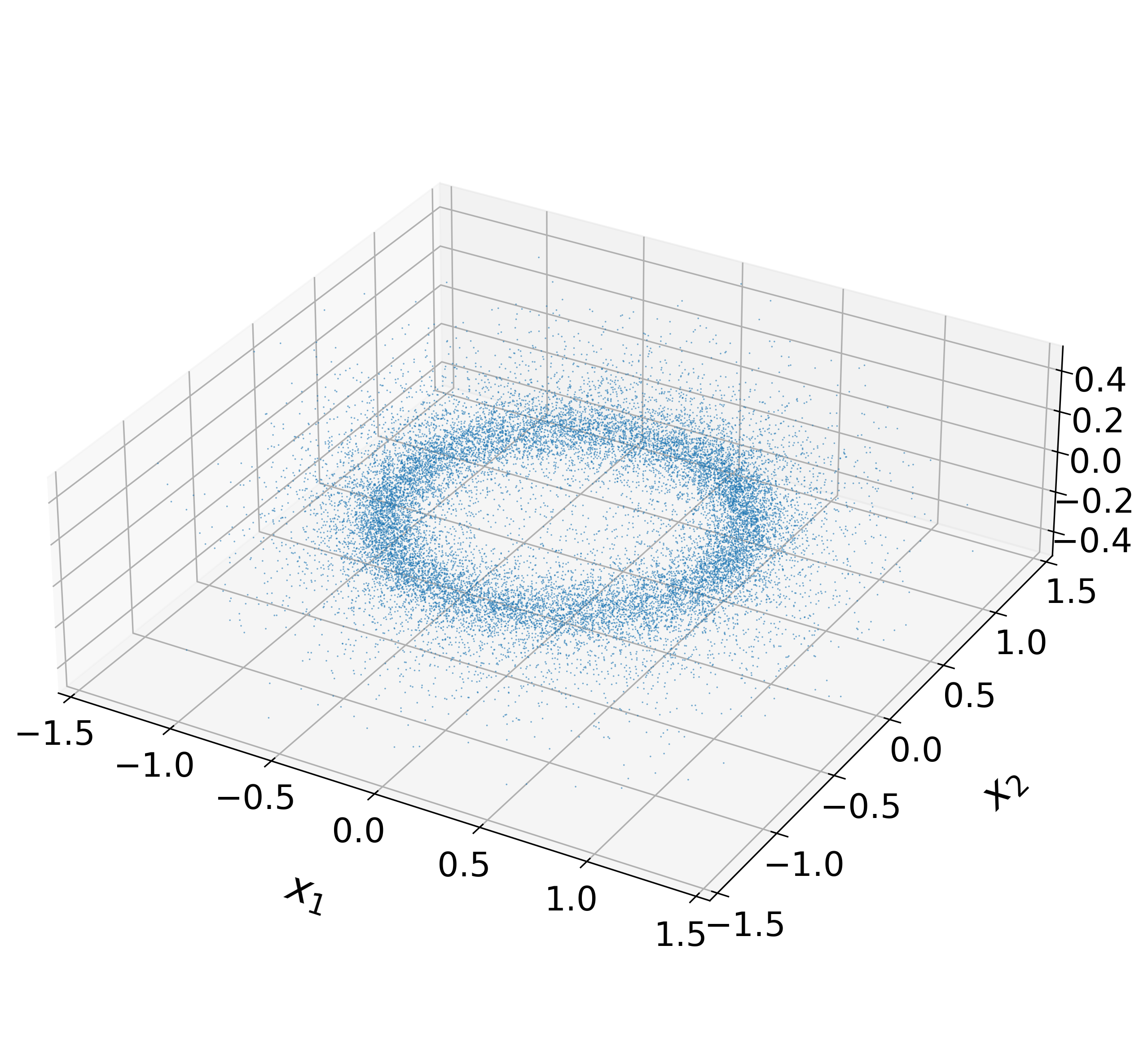}}
  \hfill
  \subfloat[\(t = 0.05\)]{\includegraphics[width=0.3\textwidth]{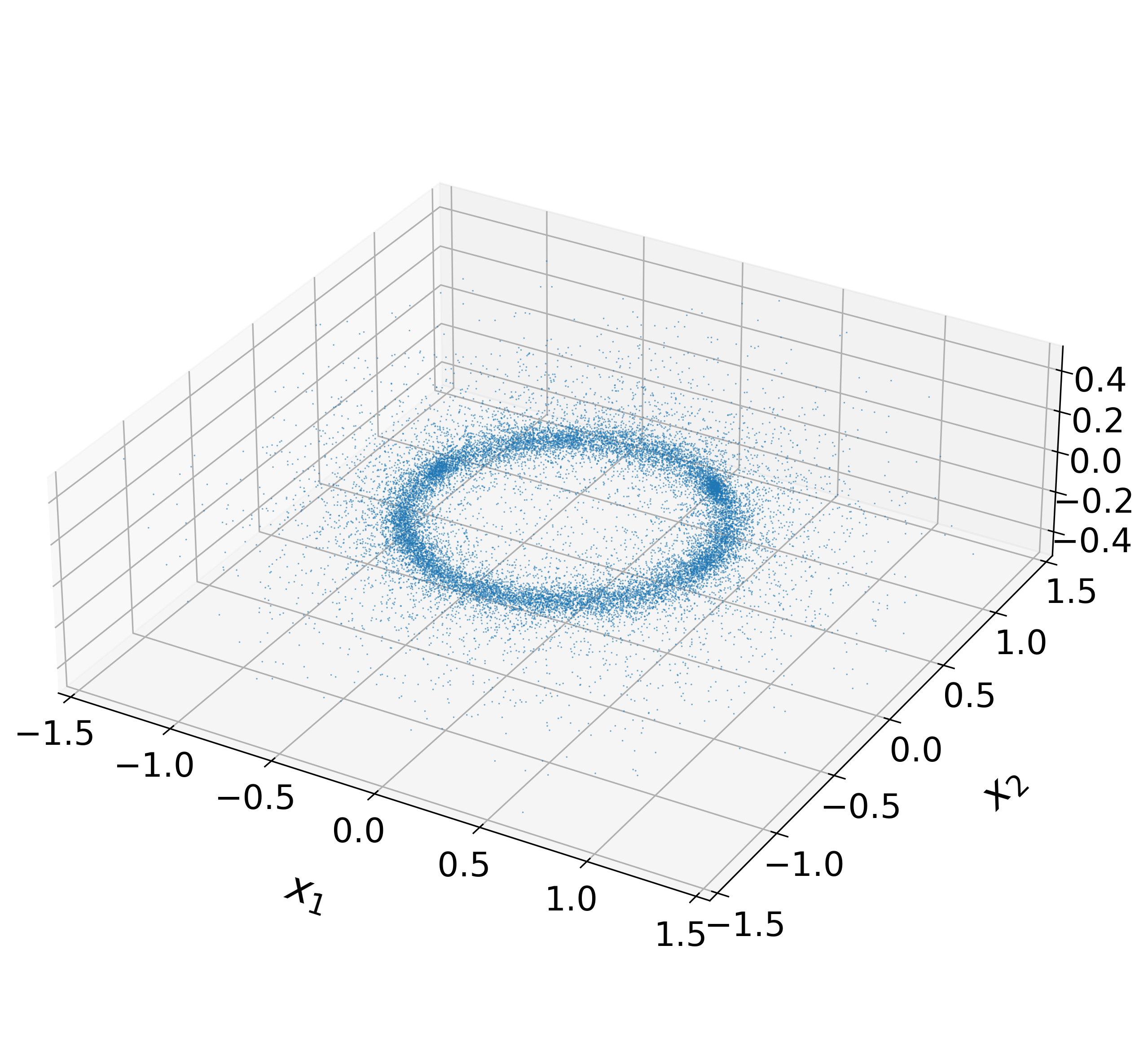}}
  \hfill
  \subfloat[\(t = 0.06\)]{\includegraphics[width=0.3\textwidth]{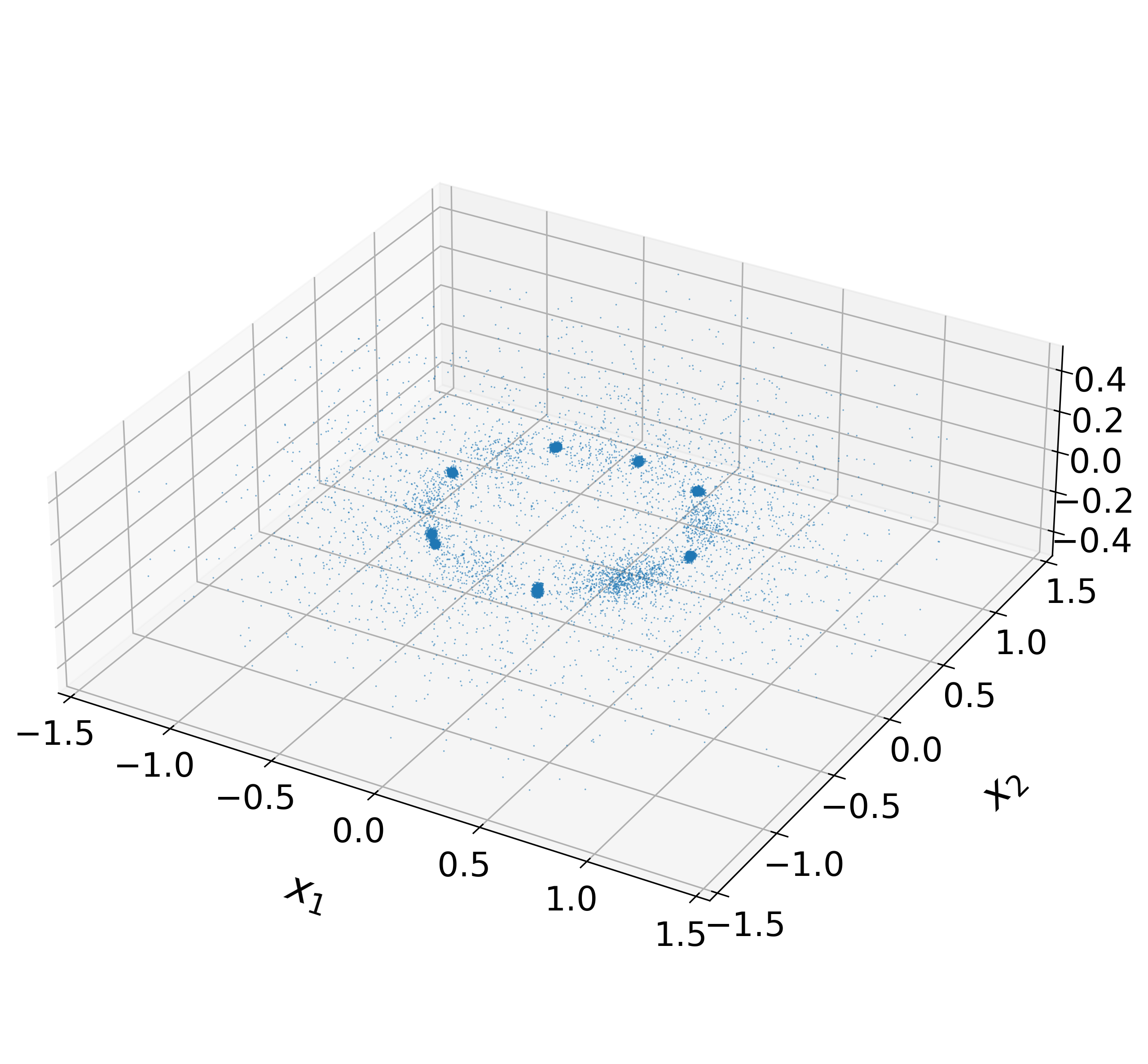}}
  
  \caption{Toroidal collapse evolution (Scatter Plot). Parameters $P=2^{20},H=256,L=8,H_{0}=\left\lceil 8H^{8/13}L^{5/13}\right\rceil,\tau=5\times 10^{-5}$.}
  \label{fig:torus_1}
\end{figure}

\begin{figure}[htbp]
  \centering
  \subfloat[\(t = 0.01\)]{\includegraphics[width=0.3\textwidth]{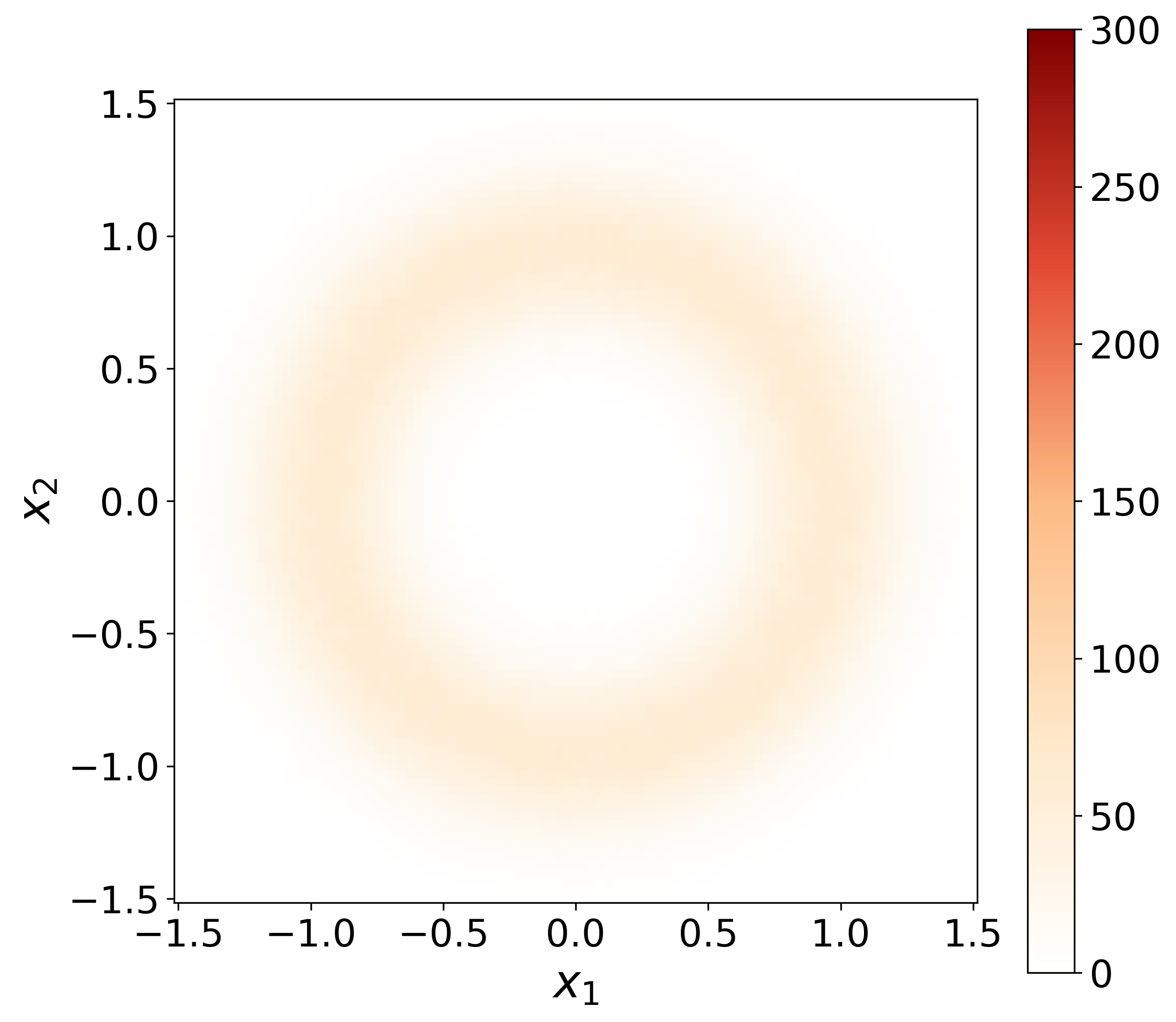}}
  \hfill
  \subfloat[\(t = 0.02\)]{\includegraphics[width=0.3\textwidth]{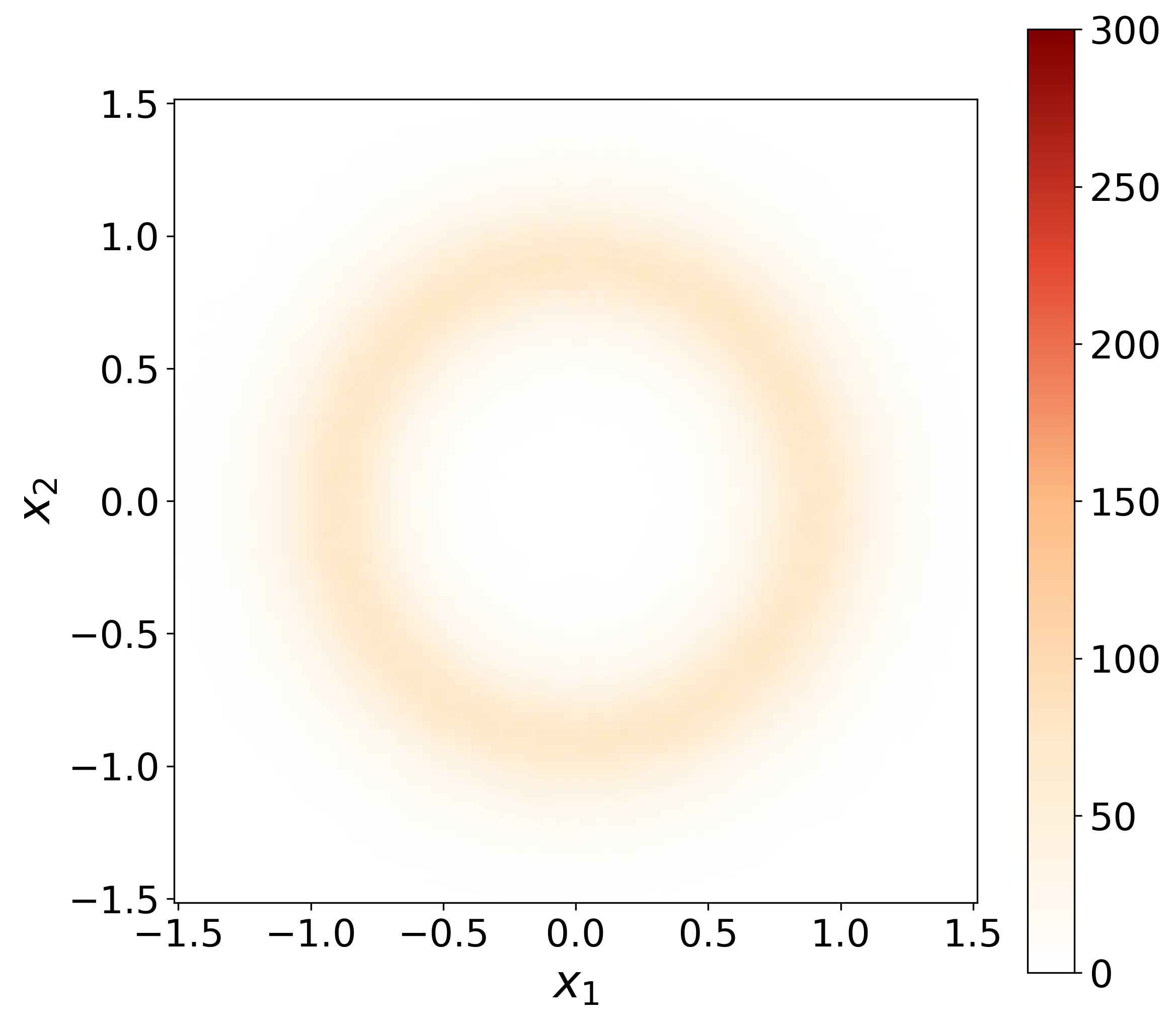}}
  \hfill
  \subfloat[\(t = 0.03\)]{\includegraphics[width=0.3\textwidth]{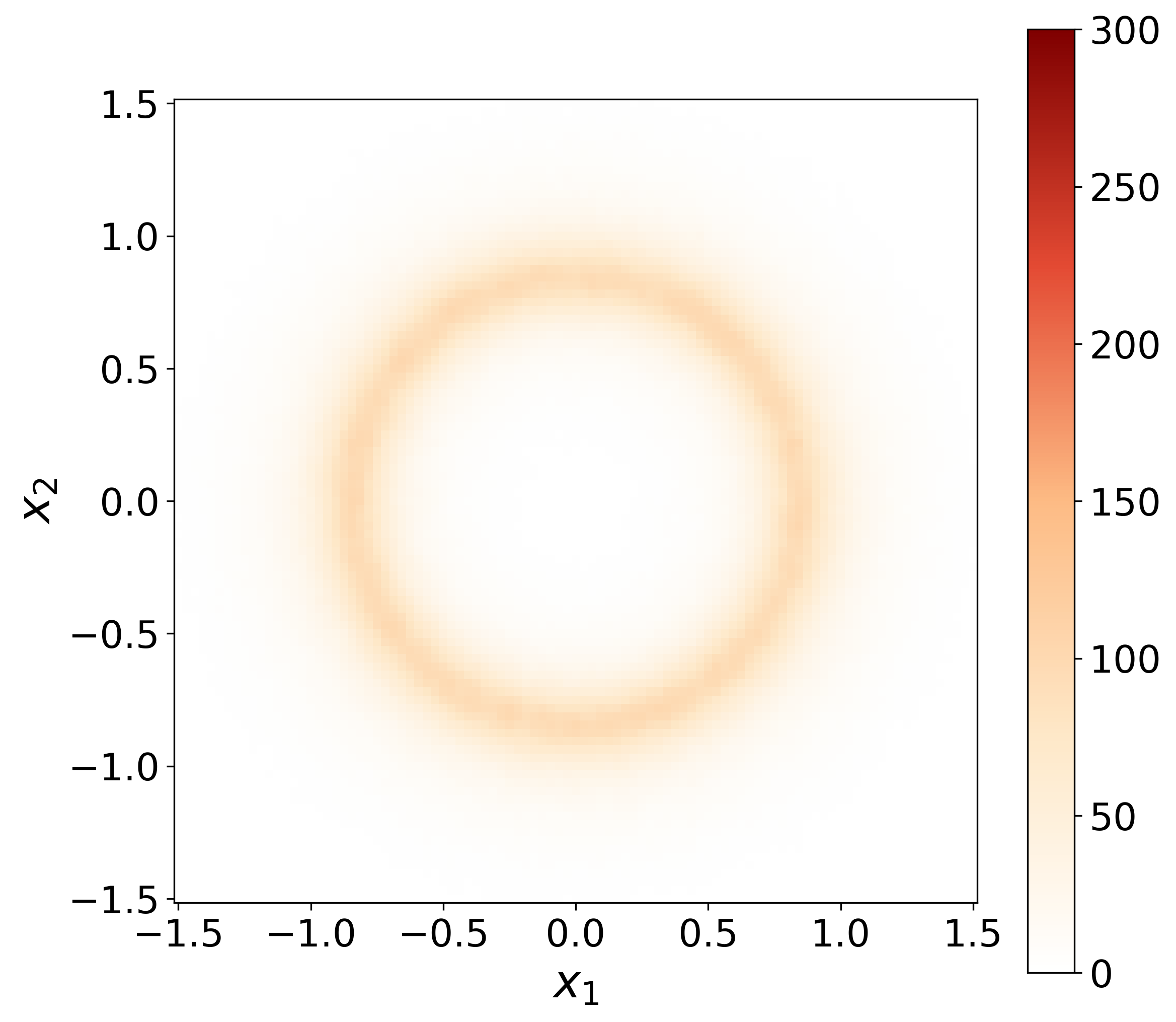}}

  \vspace{0.3cm}
  \subfloat[\(t = 0.04\)]{\includegraphics[width=0.3\textwidth]{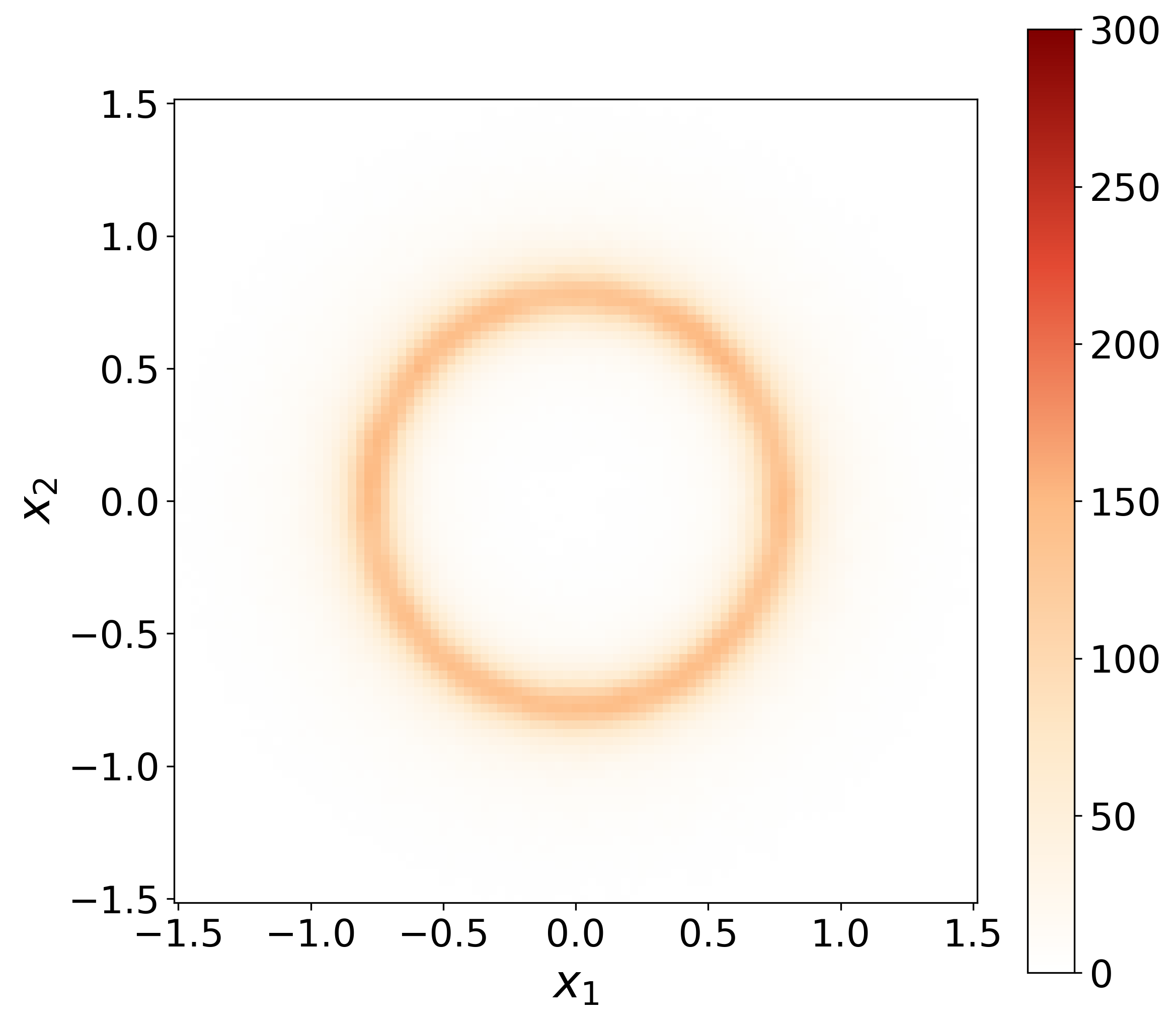}}
  \hfill
  \subfloat[\(t = 0.05\)]{\includegraphics[width=0.3\textwidth]{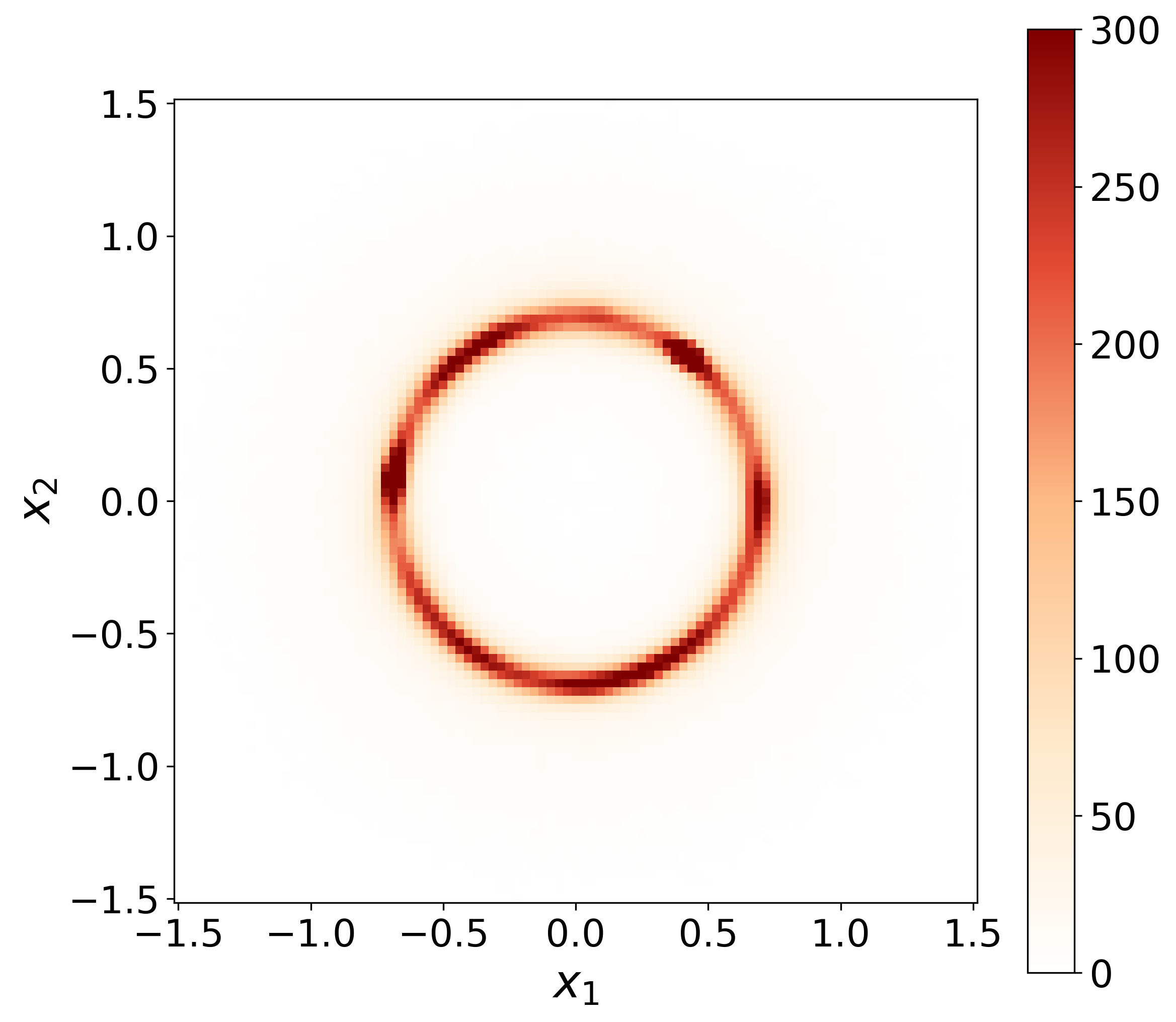}}
  \hfill
  \subfloat[\(t = 0.06\)]{\includegraphics[width=0.3\textwidth]{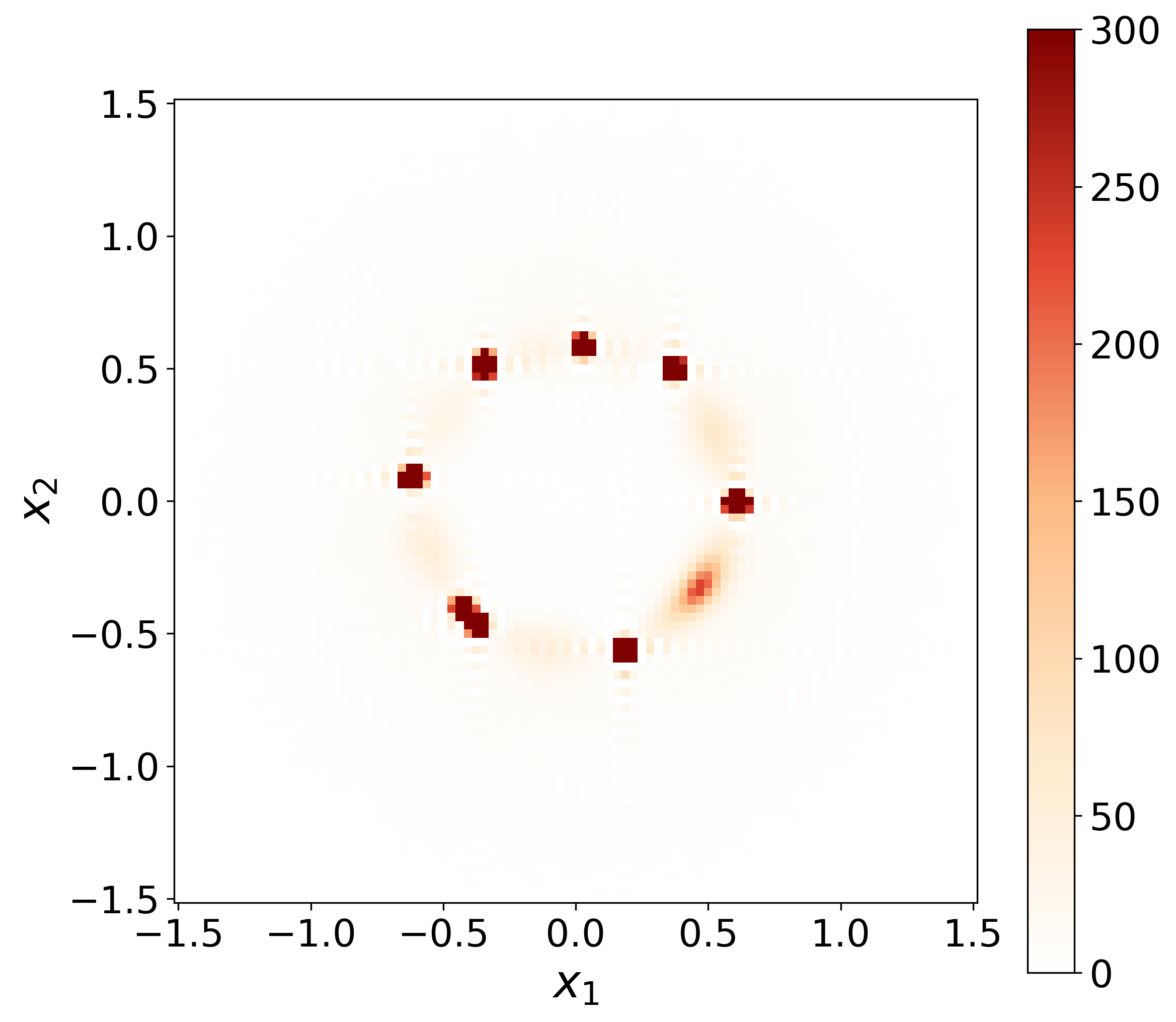}}

  \caption{Toroidal collapse evolution (Density Plot).}
  \label{fig:torus_2}
\end{figure}

\begin{figure}[htbp]
  \centering
  \subfloat[Evolution of \(\|\rho_{\mathrm{num}}\|_{\infty}\).]
  {
    \includegraphics[width=0.42\textwidth]{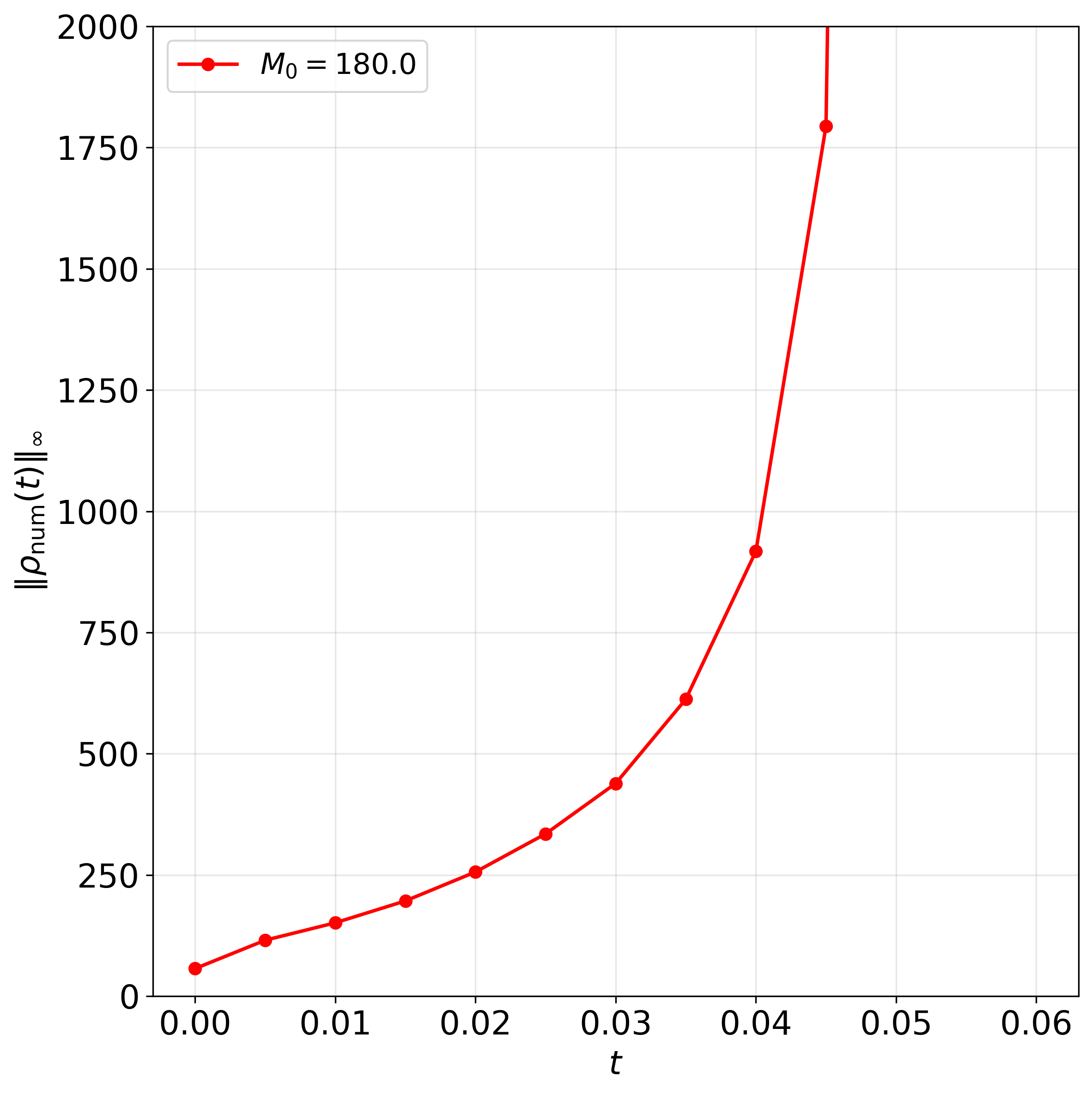}
    \label{fig:torus_3}
  }
  \hfill
  \subfloat[Evolution of the numerical radius \(R_\text{num}\).]
  {
    \includegraphics[width=0.42\textwidth]{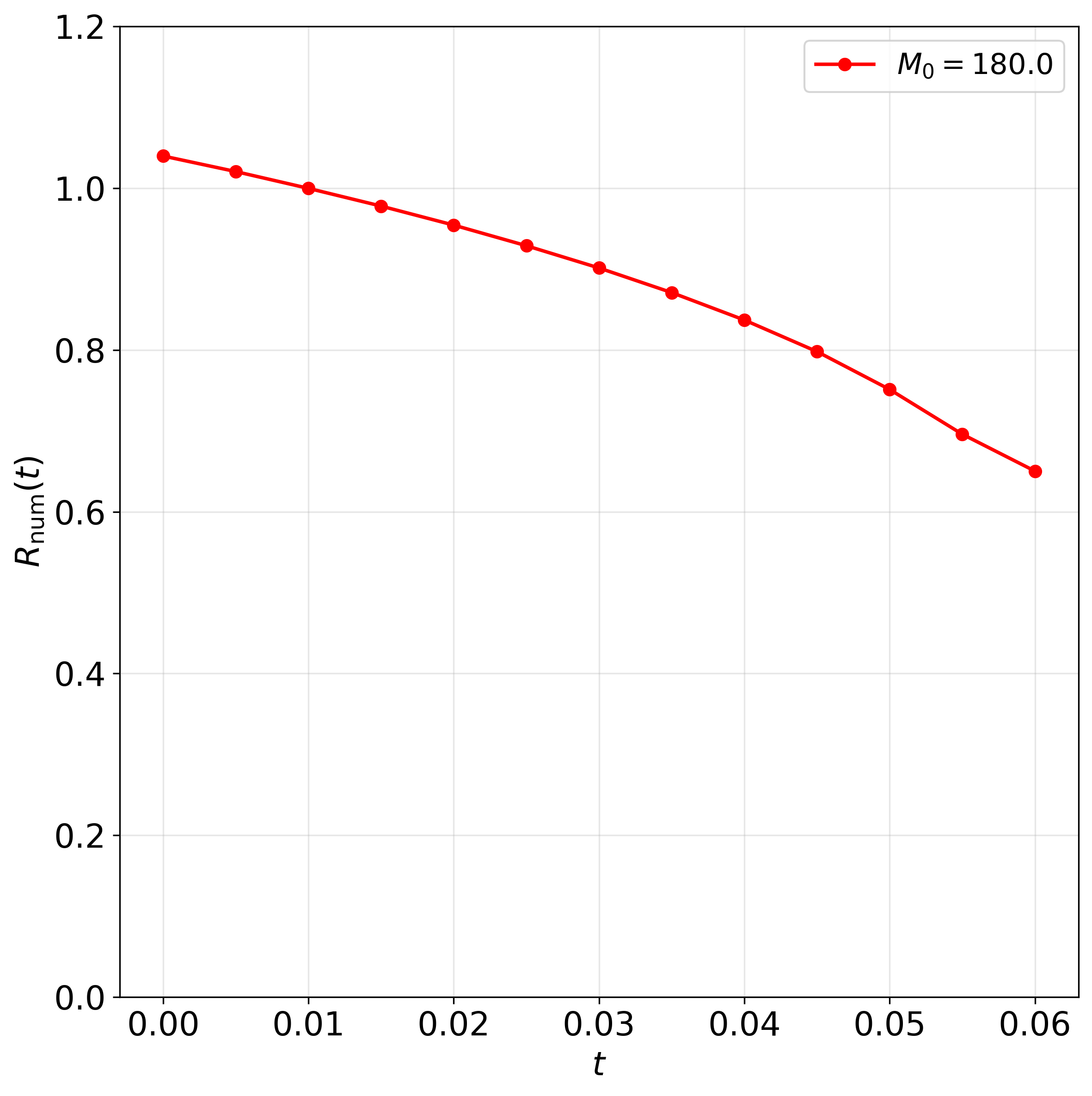}
    \label{fig:torus_4}
  }

  \caption{Quantitative indicators for the three-dimensional toroidal problem.}
  \label{fig:torus_quantitative}
\end{figure}

\begin{color}{black}Figure~\ref{fig:torus_1} shows the scatter plots of the particle distribution at different time steps. Figure~\ref{fig:torus_2} shows the corresponding numerical density, projected onto the \(x_1\)-\(x_2\) plane at the same time steps. The numerical density is obtained by convolving the empirical full-particle measure with the Gaussian kernel \(G_0\). Figure~\ref{fig:torus_3} presents the temporal evolution of the density maximum \(\|\rho_{\mathrm{num}}\|_{\infty}\), while Figure~\ref{fig:torus_4} reports the numerical radius $R_{\text{num}}(t)=1/P\sum_{p=1}^P \sqrt{\left(\widehat{X}_{p,1}^{(t/\tau)}\right)^2+\left(\widehat{X}_{p,2}^{(t/\tau)}\right)^2}$.\end{color} The mass progressively collapses into a thin ring with an asymptotically linear density profile as \(t \to T\). After complete collapse, numerical stability degrades severely, revealing an extreme sensitivity to asymmetric perturbations. This instability drives particles to aggregate into point-like clusters (singularities) within localized regions, each with a radius of approximately one grid spacing. High particle number and grid resolution are necessary to resolve the ring structure. The original SIPF method ($P=10^4, H=32$) is not sufficient to resolve the ring structure, as the particles undergo a blow-up at the origin rather than aggregating into a ring.


\section{Conclusions} \label{sec:conclusion}
In this paper, we develop the SIPF-PIC method, a hybrid particle-spectral particle-in-cell method that combines localized particle-grid interpolations with global spectral field solvers accelerated by the FFT method. 
This new method achieves a computational complexity of $O(P + H^3 \log H)$, enabling efficient simulations of 3D Keller-Segel (KS) systems with large particle numbers and fine grids. A unified error analysis framework is established to rigorously quantify approximation errors arising from particle-grid interpolations and spectral field evaluations, thus providing convergence guarantees for our SIPF-PIC method in computing non-blow-up solutions of parabolic-parabolic KS systems.

The SIPF-PIC method successfully resolves complex blow-up dynamics beyond classical single-point singularities, such as ring-shaped blow-up patterns in which mass concentrates on an evolving ring structure. Numerical experiments, including convergence tests and 3D visualizations, demonstrate the enhanced precision and scalability of the method in capturing aggregation and near-singular phenomena.  

Future work will focus on extending the SIPF-PIC method to other related systems with collective particle-field interactions. Furthermore, we aim to develop adaptive strategies to improve the robustness of the method in handling near-blow-up cases, where high-resolution tracking of density singularities remains a critical challenge.

\section{Acknowledgements}
 \noindent  ZW was partially supported by NTU SUG-023162-00001 and MOE AcRF Tier 1 Grant RG17/24. JX was partially supported by NSF Grant DMS-2309520, the Swedish Research Council
Grant No. 2021-06594 at the Institut Mittag-Leffler in Djursholm, Sweden, and the E. Schrödinger Institute in Vienna, Austria, both during his stay in the fall of 2025. ZZ was partially supported by the National Natural Science Foundation of China (Projects 92470103 and 12171406), the Hong Kong RGC Grant (Projects 17304324 and 17300325), the Seed Funding Programme for Basic Research (HKU), and the Hong Kong RGC Research Fellow Scheme 2025. The computations were performed at the research computing facilities provided by the Greenplanet Cluster at UC Irvine and the Information Technology Services, The University of Hong Kong.


\printbibliography


\end{document}